\documentclass[11pt]{article}
\usepackage{amsmath}
\usepackage{amssymb}
\usepackage{amscd}
\usepackage{xspace}
\usepackage{verbatim}
\newcommand{\mysection}[1]{
\section{#1}\setcounter{equation}{0}}
\title{\bf Capacitary estimates of solutions \\ of semilinear parabolic 
equations }
\author{ {\bf Moshe Marcus}\\
{\small Department of Mathematics,}\\
 {\small  Technion, Haifa, ISRAEL}
\and {\bf Laurent Veron}\\
{\small Department of Mathematics,}\\
 {\small  Univ. of Tours,  FRANCE}
}

\date{}
\begin{document}
\maketitle

\newcommand{\txt}[1]{\;\text{ #1 }\;}
\newcommand{\tbf}{\textbf}
\newcommand{\tit}{\textit}
\newcommand{\tsc}{\textsc}
\newcommand{\trm}{\textrm}
\newcommand{\mbf}{\mathbf}
\newcommand{\mrm}{\mathrm}
\newcommand{\bsym}{\boldsymbol}
\newcommand{\scs}{\scriptstyle}
\newcommand{\sss}{\scriptscriptstyle}
\newcommand{\txts}{\textstyle}
\newcommand{\dsps}{\displaystyle}
\newcommand{\fnz}{\footnotesize}
\newcommand{\scz}{\scriptsize}
\newcommand{\be}{
\begin{equation}
}
\newcommand{\bel}[1]{
\begin{equation}
\label{#1}}
\newcommand{\ee}{
\end{equation}
}
\newcommand{\eqnl}[2]{
\begin{equation}
\label{#1}{#2}
\end{equation}
}
\newtheorem{subn}{\name}
\renewcommand{\thesubn}{}
\newcommand{\bsn}[1]{\def\name{#1}
\begin{subn}}
\newcommand{\esn}{
\end{subn}}
\newtheorem{sub}{\name}[section]
\newcommand{\dn}[1]{\def\name{#1}}   
\newcommand{\bs}{
\begin{sub}}
\newcommand{\es}{
\end{sub}}
\newcommand{\bsl}[1]{
\begin{sub}\label{#1}}
\newcommand{\bth}[1]{\def\name{Theorem}
\begin{sub}\label{t:#1}}
\newcommand{\blemma}[1]{\def\name{Lemma}
\begin{sub}\label{l:#1}}
\newcommand{\bcor}[1]{\def\name{Corollary}
\begin{sub}\label{c:#1}}
\newcommand{\bdef}[1]{\def\name{Definition}
\begin{sub}\label{d:#1}}
\newcommand{\bprop}[1]{\def\name{Proposition}
\begin{sub}\label{p:#1}}
\newcommand{\R}{\eqref}
\newcommand{\rth}[1]{Theorem~\ref{t:#1}}
\newcommand{\rlemma}[1]{Lemma~\ref{l:#1}}
\newcommand{\rcor}[1]{Corollary~\ref{c:#1}}
\newcommand{\rdef}[1]{Definition~\ref{d:#1}}
\newcommand{\rprop}[1]{Proposition~\ref{p:#1}}
\newcommand{\BA}{
\begin{array}}
\newcommand{\EA}{
\end{array}}
\newcommand{\BAN}{\renewcommand{\arraystretch}{1.2}
\setlength{\arraycolsep}{2pt}
\begin{array}}
\newcommand{\BAV}[2]{\renewcommand{\arraystretch}{#1}
\setlength{\arraycolsep}{#2}
\begin{array}}
\newcommand{\BSA}{
\begin{subarray}}
\newcommand{\ESA}{
\end{subarray}}
\newcommand{\BAL}{
\begin{aligned}}
\newcommand{\EAL}{
\end{aligned}}
\newcommand{\BALG}{
\begin{alignat}}
\newcommand{\EALG}{
\end{alignat}}
\newcommand{\BALGN}{
\begin{alignat*}}
\newcommand{\EALGN}{
\end{alignat*}}
\newcommand{\note}[1]{\textit{#1.}\hspace{2mm}}
\newcommand{\Proof}{\note{Proof}}
\newcommand{\qeda}{\hspace{10mm}\hfill $\square$}
\newcommand{\qed}{\\
${}$ \hfill $\square$}
\newcommand{\Remark}{\note{Remark}}
\newcommand{\modin}{$\,$\\
[-4mm] \indent}
\newcommand{\forevery}{\quad \forall}
\newcommand{\set}[1]{\{#1\}}
\newcommand{\setdef}[2]{\{\,#1:\,#2\,\}}
\newcommand{\setm}[2]{\{\,#1\mid #2\,\}}
\newcommand{\lra}{\longrightarrow}
\newcommand{\lla}{\longleftarrow}
\newcommand{\llra}{\longleftrightarrow}
\newcommand{\Lra}{\Longrightarrow}
\newcommand{\Lla}{\Longleftarrow}
\newcommand{\Llra}{\Longleftrightarrow}
\newcommand{\warrow}{\rightharpoonup}
\newcommand{
\paran}[1]{\left (#1 \right )}
\newcommand{\sqbr}[1]{\left [#1 \right ]}
\newcommand{\curlybr}[1]{\left \{#1 \right \}}
\newcommand{\abs}[1]{\left |#1\right |}
\newcommand{\norm}[1]{\left \|#1\right \|}
\newcommand{
\paranb}[1]{\big (#1 \big )}
\newcommand{\lsqbrb}[1]{\big [#1 \big ]}
\newcommand{\lcurlybrb}[1]{\big \{#1 \big \}}
\newcommand{\absb}[1]{\big |#1\big |}
\newcommand{\normb}[1]{\big \|#1\big \|}
\newcommand{
\paranB}[1]{\Big (#1 \Big )}
\newcommand{\absB}[1]{\Big |#1\Big |}
\newcommand{\normB}[1]{\Big \|#1\Big \|}

\newcommand{\thkl}{\rule[-.5mm]{.3mm}{3mm}}
\newcommand{\thknorm}[1]{\thkl #1 \thkl\,}
\newcommand{\trinorm}[1]{|\!|\!| #1 |\!|\!|\,}
\newcommand{\bang}[1]{\langle #1 \rangle}
\def\angb<#1>{\langle #1 \rangle}
\newcommand{\vstrut}[1]{\rule{0mm}{#1}}
\newcommand{\rec}[1]{\frac{1}{#1}}
\newcommand{\opname}[1]{\mbox{\rm #1}\,}
\newcommand{\supp}{\opname{supp}}
\newcommand{\dist}{\opname{dist}}
\newcommand{\myfrac}[2]{{\displaystyle \frac{#1}{#2} }}
\newcommand{\myint}[2]{{\displaystyle \int_{#1}^{#2}}}
\newcommand{\mysum}[2]{{\displaystyle \sum_{#1}^{#2}}}
\newcommand {\dint}{{\displaystyle \int\!\!\int}}
\newcommand{\q}{\quad}
\newcommand{\qq}{\qquad}
\newcommand{\hsp}[1]{\hspace{#1mm}}
\newcommand{\vsp}[1]{\vspace{#1mm}}
\newcommand{\ity}{\infty}
\newcommand{\prt}{
\partial}
\newcommand{\sms}{\setminus}
\newcommand{\ems}{\emptyset}
\newcommand{\ti}{\times}
\newcommand{\pr}{^\prime}
\newcommand{\ppr}{^{\prime\prime}}
\newcommand{\tl}{\tilde}
\newcommand{\sbs}{\subset}
\newcommand{\sbeq}{\subseteq}
\newcommand{\nind}{\noindent}
\newcommand{\ind}{\indent}
\newcommand{\ovl}{\overline}
\newcommand{\unl}{\underline}
\newcommand{\nin}{\not\in}
\newcommand{\pfrac}[2]{\genfrac{(}{)}{}{}{#1}{#2}}

\def\ga{\alpha}     \def\gb{\beta}       \def\gg{\gamma}
\def\gc{\chi}       \def\gd{\delta}      \def\ge{\epsilon}
\def\gth{\theta}                         \def\vge{\varepsilon}
\def\gf{\phi}       \def\vgf{\varphi}    \def\gh{\eta}
\def\gi{\iota}      \def\gk{\kappa}      \def\gl{\lambda}
\def\gm{\mu}        \def\gn{\nu}         \def\gp{\pi}
\def\vgp{\varpi}    \def\gr{\rho}        \def\vgr{\varrho}
\def\gs{\sigma}     \def\vgs{\varsigma}  \def\gt{\tau}
\def\gu{\upsilon}   \def\gv{\vartheta}   \def\gw{\omega}
\def\gx{\xi}        \def\gy{\psi}        \def\gz{\zeta}
\def\Gg{\Gamma}     \def\Gd{\Delta}      \def\Gf{\Phi}
\def\Gth{\Theta}
\def\Gl{\Lambda}    \def\Gs{\Sigma}      \def\Gp{\Pi}
\def\Gw{\Omega}     \def\Gx{\Xi}         \def\Gy{\Psi}

\def\CS{{\mathcal S}}   \def\CM{{\mathcal M}}   \def\CN{{\mathcal N}}
\def\CR{{\mathcal R}}   \def\CO{{\mathcal O}}   \def\CP{{\mathcal P}}
\def\CA{{\mathcal A}}   \def\CB{{\mathcal B}}   \def\CC{{\mathcal C}}
\def\CD{{\mathcal D}}   \def\CE{{\mathcal E}}   \def\CF{{\mathcal F}}
\def\CG{{\mathcal G}}   \def\CH{{\mathcal H}}   \def\CI{{\mathcal I}}
\def\CJ{{\mathcal J}}   \def\CK{{\mathcal K}}   \def\CL{{\mathcal L}}
\def\CT{{\mathcal T}}   \def\CU{{\mathcal U}}   \def\CV{{\mathcal V}}
\def\CZ{{\mathcal Z}}   \def\CX{{\mathcal X}}   \def\CY{{\mathcal Y}}
\def\CW{{\mathcal W}} \def\CQ{{\mathcal Q}} 
\def\BBA {\mathbb A}   \def\BBb {\mathbb B}    \def\BBC {\mathbb C}
\def\BBD {\mathbb D}   \def\BBE {\mathbb E}    \def\BBF {\mathbb F}
\def\BBG {\mathbb G}   \def\BBH {\mathbb H}    \def\BBI {\mathbb I}
\def\BBJ {\mathbb J}   \def\BBK {\mathbb K}    \def\BBL {\mathbb L}
\def\BBM {\mathbb M}   \def\BBN {\mathbb N}    \def\BBO {\mathbb O}
\def\BBP {\mathbb P}   \def\BBR {\mathbb R}    \def\BBS {\mathbb S}
\def\BBT {\mathbb T}   \def\BBU {\mathbb U}    \def\BBV {\mathbb V}
\def\BBW {\mathbb W}   \def\BBX {\mathbb X}    \def\BBY {\mathbb Y}
\def\BBZ {\mathbb Z}

\def\GTA {\mathfrak A}   \def\GTB {\mathfrak B}    \def\GTC {\mathfrak C}
\def\GTD {\mathfrak D}   \def\GTE {\mathfrak E}    \def\GTF {\mathfrak F}
\def\GTG {\mathfrak G}   \def\GTH {\mathfrak H}    \def\GTI {\mathfrak I}
\def\GTJ {\mathfrak J}   \def\GTK {\mathfrak K}    \def\GTL {\mathfrak L}
\def\GTM {\mathfrak M}   \def\GTN {\mathfrak N}    \def\GTO {\mathfrak O}
\def\GTP {\mathfrak P}   \def\GTR {\mathfrak R}    \def\GTS {\mathfrak S}
\def\GTT {\mathfrak T}   \def\GTU {\mathfrak U}    \def\GTV {\mathfrak V}
\def\GTW {\mathfrak W}   \def\GTX {\mathfrak X}    \def\GTY {\mathfrak Y}
\def\GTZ {\mathfrak Z}   \def\GTQ {\mathfrak Q}

\font\Sym= msam10 
\def\SYM#1{\hbox{\Sym #1}}
\newcommand{\bdw}{\prt\Gw\xspace}
\medskip
\begin{abstract}
We prove that any positive solution of $ \prt_tu-\Delta u+u^q=0$ ($q>1$) in $\BBR^N\ti(0,\infty)$ with initial trace  $(F,0)$, where $F$ is a closed subset of $\BBR^N$ can be represented, up to two universal multiplicative constants, by a series involving the Bessel capacity $C_{2/q,q'}$. As a consequence we prove that there exists a unique positive solution of the equation with such an initial trace. We also characterize the blow-up set of $u(x,t)$ when $t\downarrow 0$ , by using the "density" of $F$ expressed 
in terms of the  $C_{2/q,q'}$-Bessel capacity.
\end{abstract}

\noindent
{\it \footnotesize 2000 Mathematics Subject Classification}. {\scriptsize
35K05;35K55; 31C15; 31B10; 31C40}.\\
{\it \footnotesize Key words}. {\scriptsize Heat equation; singularities; Borel measures; Besov spaces;  real interpolation; Bessel capacities; quasi-additivity; capacitary measures; Wiener type test; initial trace.}
\tableofcontents
\mysection {Introduction}

Let $T\in (0,\infty]$ and $Q_{T}=\BBR^{N}\ti (0,T]$ ($N\geq 1$). If $q>1$ and $u\in
C^{2}(Q_{T})$ is nonnegative and verifies
\begin {equation}
\label {mequ}
\prt_{t}u-\Gd u+u^q=0\quad\mbox {in }Q_{T},
\end {equation}
it has been proven by Marcus and V\'eron \cite {MV2} that there exists 
a unique outer-regular positive Borel measure $\gn$
in $\BBR^{N}$ such that 
\begin {equation}
\label {tr1}
\lim_{t\to 0}u(.,t)=\gn,
\end {equation}
in the sense of Borel measures; the set of such measures is denoted by $\GTB_{_{+}}^{reg}(\BBR^{N})$. To each of its element $\gn$ is associated 
a unique couple $(\CS_{\gn},\gm_{\gn})$ (we write 
$\gn\approx (\CS_{\gn},\gm_\gn)$) where $\CS_\gn$, {\it the singular part} of $\gn$, is a closed subset 
of $\BBR^{N}$ and $\gm_{\gn}$, {\it the regular part} is a nonnegative 
Radon measure on $\CR_{\gn}=\BBR^{N}\setminus \CS_{\gn}$. In this setting,  
relation $(\ref{tr1})$ has the following meaning :
\begin {equation}
\label {tr2}\BA{lcl} (i)\qquad\,\lim_{t\to 0}\int_{\CR_{\gn}}u(.,t)\gz dx =\myint{\CR_{\gn}}{}\gz
d\gm_{\gn},&\qquad\forall\gz\in 
C_{0}(\CR_{\gn}),\\
[2mm] (ii)\phantom{---,}\lim_{t\to 0}\myint{\CO}{}u(.,t) dx=\infty,\quad\;\;& \qquad\forall
\CO\subset\BBR^{N}\mbox { open},
\;\CO\cap \CS_{\gn}\neq\emptyset. 
\EA
\end {equation}
The measure $\gn$ is by definition {\it the initial trace} of $u$ and 
denoted by $Tr_{\BBR^{N}}(u)$. It is wellknown that equation $(\ref{mequ})$ admits a critical exponent
$$
1<q<q_{c}=1+\frac{N}{2}.$$
This is due to the fact, proven by Brezis and Friedman \cite{BF},  that if $q\geq q_c$, isolated singularities of solutions of 
$(\ref{mequ})$ in $\BBR^N\setminus\{0\}$ are removable. Conversely, if $1<q<q_c$, 
it is proven by the same authors that for any $k>0$, equation $(\ref{mequ})$ admits a unique solution $u_{k\gd_0}$ with initial data $k\gd_0$. This existence and uniqueness results extends in a simple way if the initial data $k\gd_0$ is replaced by any 
Radon measure $\gm$ in  $\BBR^N$ (see \cite{Br}). Furthermore, if $k\to\infty$, $u_{k\gd_0}$ increases and converges to a positive, radial and self-similar solution $u_\infty$ of $(\ref{mequ})$. Writing it under the form $u_\infty(x,t)=t^{-\frac{1}{q-1}}f(\abs{x}/\sqrt t)$, $f$ is a positive  solution of
\begin {equation}
\label {bila'}\left\{\BA{ll}
\Gd f+\frac{1}{2}y.Df+\frac{1}{q-1}f-f^q=0\quad\mbox {in }\BBR^N\\[2mm]
\phantom{--,,}\lim_{\abs y\to\infty}\abs y^{\frac{2}{q-1}}f(y)=0.
\EA\right.
\end {equation}
The existence, uniqueness and the expression of the asymptotics of $f$ has been studied thoroughly by Brezis, Peletier and Terman in \cite {BPT}. Later on, Marcus and V\'eron proved in \cite{MV2} that in the same range of exponents, 
for any $\gn\in\GTB_{_{+}}^{reg}(\BBR^{N})$, the 
Cauchy problem
\begin {equation}
\label {CD}\left\{\BA{rll}
\prt_{t}u-\Gd u+u^q=0\;&\quad\mbox {in }Q_{\infty},\\
[2mm] Tr_{\BBR^{N}}(u)=\gn ,&&
\EA\right.
\end {equation}
admits a unique positive solution. This result means that the initial trace establishes a one to one correspondence between the set of positive solutions of $(\ref{mequ})$ and $\GTB_{_{+}}^{reg}(\BBR^{N})$. A key step for proving the uniqueness is the following inequalities
\begin {equation}
\label {bila}
t^{-\frac{1}{q-1}}f(\abs{x-a}/\sqrt t)\leq u(x,t)\leq ((q-1)t)^{-\frac{1}{q-1}}\qquad\qquad\forall (x,t)\in Q_\infty,
\end {equation}
valid for any $a\in \CS_\gn$. 
As a consequence of Brezis and Friedman's result, if $
q\geq q_{c}$, i.e. in the {\it  supercritical range}, Problem $(\ref {CD})$ may admit no solution at all. If $\gn\in\GTB_{_{+}}^{reg}(\BBR^{N})$, $\gn\approx (\CS_\gn,\gm_\gn)$, the necessary and
sufficient conditions for the existence 
of a maximal 
solution $u=\overline u_{\gn}$ to Problem $(\ref {CD})$ are obtained in \cite {MV2}
and expressed in terms of the the Bessel 
capacity 
$C_{2/q,q'}$, 
(with $q'=q/(q-1)$). Furthermore, uniqueness does not hold in general as it was pointed out by Le Gall \cite {LG1}. In the particular case where $\CS_{\gn}=\emptyset$ and $\gn$ is simply the Radon measure $\gm_{\gn}$, the necessary and sufficient condition for solvability is that $\gm_\gn$ does not charge Borel subsets with $C_{2/q,q'}$-capacity zero. This result
was already proven by Baras and Pierre \cite {BP2} in the 
particular case 
of bounded measures and extended by Marcus and V\'eron \cite {MV2} to the general case. We  denote by $\GTM^q_{_{+}}(\BBR^{N})$ the positive cone of  
the space $\GTM^q(\BBR^{N})$ of Radon 
measures which do not charge Borel subsets with zero $C_{2/q,q'}$-capacity.
Notice that $W^{-2/q,q}(\BBR^N)\cap \frak M^b_+(\BBR^N)$ is a subset of $\GTM^q_{_{+}}(\BBR^{N})$ where $\frak M^b_+(\BBR^N)$ is the cone of positive bounded Radon mesures in $\BBR^N$. For such measures, uniqueness always 
holds and we denote $\overline u_{\gm_{\gn}}=u_{\gm_{\gn}}$. \medskip

In view of the already known results concerning the parabolic equation, it is useful to recall the main advanced results previously obtained for the  stationary equation
\begin {equation}
\label {ee}
-\Gd u+u^q=0\quad\mbox {in }\,\Gw,
\end {equation}
in a smooth bounded domain  $\Gw$ of $\BBR^N$. This equation has been intensively studied since 1993, both by probabilists (Le Gall, Dynkin, Kuznetsov) and by analysts (Marcus, V\'eron). The existence of a {\it boundary trace} for positive solutions, in the class of
outer-regular positive Borel measures on $\prt\Gw$, is proven by Le Gall \cite {LG}, \cite {LG1} in the case $q=N=2$, by probabilistic methods, and by Marcus and V\'eron in \cite {MV0}, \cite {MV1} in the general case $q>1$, $N>1$. The existence of a critical exponent $q_e=(N+1)/(N-1)$ is due to Gmira and V\'eron \cite{GV} who shew that, if $q\geq q_e$ boundary isolated singularities of solutions of $(\ref{ee})$ are removable, which is not the case if $1<q<q_e$. In this {\it subcritical case} Le Gall and Marcus and V\'eron proved that the boundary trace establishes a one to one correspondence between positive solutions of $(\ref{ee})$ in $\Gw$ and outer regular positive Borel measures on $\prt\Gw$. This fundamental result does not hold in the {\it supercritical case} $q\geq q_e$.
In \cite {DK2} Dynkin and Kuznetsov introduced the notion of $\gs$-moderate solution which means that $u$ is a positive solution of $(\ref{ee})$ such that there exists an increasing sequence of positive Radon measures on $\prt\Gw$
$\{\gm_n\} $ belonging to $W^{-2/q,q'}(\prt\Gw)$ such that the corresponding solutions $v=v_{\gm_n}$ of 
\begin {equation}
\label {ee1}\left\{\BA {l}
-\Gd v+v^q=0\quad\mbox { in }\,\Gw\\
\phantom{..\Gd u+u}
v=\gm_n\quad\mbox {in }\,\prt\Gw
\EA\right.\end {equation}
converges to $u$ locally uniformly in $\Gw$. This class of solutions plays a fundamental role since Dynkin and Kuznetsov proved that a $\gs$-moderate solution of $(\ref{ee})$ is uniquely determined by its {\it fine trace}, a new notion of trace introduced in order to avoid the non-uniqueness phenomena. Later on,
it is proved by Mselati (if $q=2$)  \cite{Ms}, then by Dynkin (if $q_e\leq q\leq 2$) \cite{Dy} and finally by Marcus with no restriction on $q$ \cite {Mar}, that all the positive solutions of $(\ref{ee})$ are $\gs$-moderate. One of the key-stones element in their proof (partially probabilistic) is the fact that the maximal solution $\overline u_K$ of $(\ref {ee})$ with a boundary trace vanishing outside a compact subset $K\subset\prt\Gw$ is indeed
$\gs$-moderate. This deep result was obtained by a combination of probabilistic and analytic methods by Mselati \cite {Ms} in the case $q=2$ and by purely analytic tools by Marcus and V\'eron \cite{MV5}, \cite{MV6} in the case $q\geq q_e$. Defining $\underline u_K$ as the largest $\gs$-moderate solution of $(\ref{ee})$ with a boundary trace concentrated on $K$, the crucial step in Marcus-V\'eron's proof (non probabilistic) is the bilateral estimate satisfied by $\overline u_K$ and $\underline u_K$
\begin {equation}\label {ee2}
C^{-1}\gr(x)W_K(x)\leq \underline u_K(x)\leq \overline u_K(x)\leq C\gr(x)W_K(x).
\end {equation}
In this expression $C=C(\Gw,q)$, $\gr(x)=\dist(x,\prt\Gw)$ and $W_F(x)$ is the {\it elliptic capacitary potential} of $K$ defined by
\begin {equation}\label {ee3}
W_K(x)=\sum_{-\infty}^\infty 2^{-\frac{m(q+1)}{q-1}}C_{2/q,q'}(2^m K_m(x)),
\end {equation}
where $K_m(x)=K\cap\{z:2^{-m-1}\leq |z-x|\leq 2^{-m}\}$, the Bessel  capacity being relative to $\BBR^{N-1}$. Note that, using a technique introduced in \cite{MV1}, inequality $\overline u_K\leq C^2\underline u_K$ implies $\underline u_K= \overline u_K$.
\medskip


The aim of this article is to initiate the fine study of the complete initial trace problem for positive solutions of $(\ref{mequ})$ in the supercritical case $q\geq q_c$ and to give in particular the parabolic counterparts of the results of \cite {Ms}, \cite{MV5} and \cite{MV6}.
Extending Dynkin's ideas to the parabolic case,  we introduce the following notion
\bdef {sigmoder} A positive solution $u$ of $(\ref{mequ})$ is called $\gs$-moderate if their exists an increasing sequence $\{\gm_n\}\subset 
W^{-2/q,q}(\BBR^N)\cap \frak M^b_+(\BBR^N)$ such that the corresponding solution $u:=u_{\gm_n}$ of 
\begin {equation}
\label {fund1}\left\{\BA {l}
\prt _tu-\Gd u+u^q=0\quad\mbox { in }\,Q_\infty\\[2mm]
\phantom{.\Gd u+u}
u(x,0)=\gm_n\quad\mbox {in }\,\BBR^N,
\EA\right.\end {equation}
converges to $u$  locally uniformly in $Q_\infty$.
\es

If $F$ is a closed subset of $\BBR^{N}$, we denote by $\overline u_{F}$ 
the maximal solution of $(\ref {mequ})$ with an initial trace vanishing 
on $F^c$, and by $\underline u_{F}$ the maximal $\gs$-moderate solution of $(\ref {mequ})$ with an initial trace vanishing 
on $F^c$. Thus $\underline u_{F}$ is defined by
\begin {equation}
\label {minimax}
\underline u_{F}=\sup\{u_{\gm}: \gm\in W^{-2/q,q}(\BBR^N)\cap \frak M^b_+(\BBR^N),\gm 
(F^c)=0\},
\end {equation}
(and clearly $W^{-2/q,q}(\BBR^N)\cap \frak M^b_+(\BBR^N)$ can be replaced by $\GTM_{+}^q(\BBR^{N})$).
One of the main goal of this article is to prove that $\overline u_{F}$ is $\gs$-moderate and more precisely,
\bth{th1} For any $q>1$ and any closed subset $F$ of $\BBR^N$, 
$\overline u_{F}=\underline u_{F}$.
\es

We define below a set function which will play a fundamental role in the sequel.
\bdef {cappot} Let $F$ be a closed subset of $\BBR^N$. The Bessel parabolic capacitary potential $W_F$ of $F$ is defined by 
\begin {equation}
\label {pot1}
 W_F(x,t)=\frac{1}{t^{\frac{N}{2}}}\mysum{n=0}{\infty}d_{n+1}^{N-\frac{2}{q-1}}e^{-\frac{n}{4}} C_{2/q,q'}\left(\myfrac {F_{n}}{d_{n+1}}\right)\qquad\forall (x,t)\in Q_\infty,
\end {equation}
where $C_{2/q,q'}$ is the $N$-dimensional Bessel capacity, $d_n=\sqrt{nt}$ and $F_n=\left\{y\in F:d_n\leq\abs{x-y}\leq d_{n+1}\right\}$.
\es

In our study, it is useful to introduce a variant of $W_F$ with the help of the Besov capacity: if $\Gw\subset\BBR^N$ is a bounded domain, we set
\begin{equation}\label{AS}
\norm{\gf}_{B_{2/q,q'}}=\left(\dint_{\!\!\!\Gw\ti\Gw}\frac{\abs{\gf(x)-\gf(y)}^{q'}}{\abs{x-y}^{N+\frac{2}{q-1}}}dx dy\right)^{1/q'},
\end{equation}
if $1<2/q<1$, and $\norm{\gf}_{B_{1,2}}=\norm{\nabla \gf}_{L^2}$ if $2/q=1$ (i.e. $N=2$ and $q=2$). The Besov capacity of a compact set $K\subset\Gw$ relative to $\Gw$ is expressed by
\begin{equation}\label{AS1}
R^\Gw_{2/q,q'}=\inf\left\{\norm{\gf}_{B_{2/q,q'}}^{q'}:\gf\in C^{\infty}_0(\Gw), 0\leq \gf\leq 1, \eta=1\text{ on K}\right\}.
\end{equation}
The Besov-parabolic capacitary potential $\tilde W_F$ of $F$ is defined by 
\begin {equation}
\label {pott2}
 \tilde W_F(x,t)=t^{-\frac{N}{2}}\mysum{n=0}{\infty}d_{n+1}^{N-\frac{2}{q-1}}e^{-\frac{n}{4}}
 R^{\Gamma_n}_{2/q,q'}\left(\myfrac {F_{n}}{d_{n+1}}\right)\qquad\forall (x,t)\in Q_\infty,
\end {equation}
where $\Gamma_n=B_{d_{n+1}}\setminus\overline{B_{d_{n}}}$. 
The Besov-parabolic capacitary potential is equivariant with respect to the same scaling transformation which let $(\ref{mequ})$ invariant in the sense that, for any $\ell>0$,
\begin {equation}
\label {inv}
\ell^{\frac{1}{q-1}} \tilde W_F(\sqrt\ell x,\ell t)= \tilde W_{F/\sqrt\ell}( x,t)\qquad\forall (x,t)\in Q_\infty.
\end {equation}
and we prove that there exists $c=c(N,q)>0$ such that 
\begin {equation}
\label {pot3}
c^{-1} \tilde W_F(x,t)\leq W_F(x,t)\leq c\tilde W_F(x,t) \qquad\forall (x,t)\in Q_\infty.
\end {equation}

One of the tool for proving \rth{th1} is the following bilateral estimate which is only meaningful in the supercritical case, otherwhile it reduces to $(\ref{bila})$;
\bth{th2} For any $q\geq q_c$ there exist two positive constants $C_1\geq C_2>0$, depending only on $N$ and $q$ such that for any closed subset $F$ of $\BBR^N$, there holds
\begin {equation}
\label {pot2}
C_2 W_F(x,t)\leq \underline u_F(x,t)\leq \overline u_F(x,t)\leq C_1 W_F(x,t)\qquad\forall (x,t)\in Q_\infty.
\end {equation}
\es

Actually our result is more general since the upper estimate in $(\ref{pot2})$ is valid for {\it any} positive solution of
\begin {equation}
\label {ineq1}
\prt_tu-\Gd u+u^q\leq 0\qquad\text{in }Q_T
\end{equation}
satisfying
\begin {equation}
\label {ineq1-1}
\lim_{t\to 0}u(x,t)=0\qquad\text{locally uniformly in }F^c.
\end {equation}
Extension to positive solutions of
\begin {equation}
\label {ineq2}
\prt_tu-\Gd u+f(u)= 0\qquad\text{in }Q_T
\end{equation}
where $f$ is continuous from $\BBR^+$ to $\BBR^+$ and satisfies
\begin {equation}\label {ineq3}
c_2r^q\leq f(r)\leq c_1r^q\qquad\forall r\geq 0
\end{equation}
for some $0<c_2\leq c_1$ is straightforward.\smallskip

This {\it quasi representation}, up to uniformly upper and lower bounded functions, is also interesting in the sense that it indicates precisely {\it how to characterize the blow-up points} of $\overline u_F=\underline u_F:=u_F$. Introducing an integral expression comparable to $ W_F$, we show in  particular the following results
\begin {equation}
\label {intpot-a}
\lim_{\gt\to 0} \gt^{\frac{2}{q-1}-N}C_{2/q,q'}\left(F\cap B_\gt(x)\right)=\gg\in[0,\infty)
 \Longrightarrow
\lim_{t\to 0}t^{\frac{1}{q-1}} u_F(x,t)=C\gg
\end {equation}
for some $C_\gg=C(N,q,\gg)>0$, and
\begin {equation}
\label {intpot-b}
\limsup_{\gt\to 0}\gt^{\frac{2}{q-1}} C_{2/q,q'}\left(\myfrac{F}{\gt}\cap B_1(x)\right)<\infty\Longrightarrow
\limsup_{t\to 0} u_F(x,t)<\infty.
\end {equation}

\medskip

 Our paper is organized as follows. In Section1 we recall some properties of the Besov spaces with fractional derivatives $B^{s,p}$ and their links with heat equation.  In Section 2 we obtain estimates 
 from above on $\overline u_{F}$. In Section 3 we give 
 estimates from below on $\underline u_{F}$. In Section 4 we prove the main theorems and expose various consequences. In Appendix we derive a series of sharp integral inequalities.\medskip 
 
 \noindent {\bf Aknowledgements} The authors are grateful to the European RTN Contract 
N$^\circ$ HPRN-CT-2002-00274 for the support 
 provided in the realization of this work. The authors are grateful to Luc Tartar for providing them the proof of the sharp Poincar\'e inequality \rprop{Poinca1} and related references.
  \mysection {Estimates from above}
 {\it Some notations.}  Let $\Gw$ be a domain in $\mathbb R^N$ with a 
 compact $C^{2}$ boundary and $T>0$. Set $B_{r}(a)$ the open ball 
 of radius $r>0$ and center $a$ (and $B_{r}(0):=B_{r}$) and
 $$
Q_{T}^\Gw:=\Gw\ti (0,T),\quad \prt_{\ell}Q_{T}^\Gw=\prt\Gw\ti (0,T), 
\quad  Q_{T}:=Q_{T}^{\BBR^{N}}, \quad Q_\infty:=Q_{\infty}^{\BBR^{N}}.
 $$
Let $\BBH^\Gw[.]$ (resp. $\BBH[.]$) denote the  heat 
potential in $\Gw$ with zero lateral boundary data (resp. the heat potential in
$\BBR^{N}$) with 
corresponding kernel 
$$
(x,y,t)\mapsto H^\Gw (x,y,t) \quad \mbox {(resp.} (x,y,t)\mapsto 
H(x,y,t)=(4\gp t)^{-\frac{N}{2}}e^{-\frac{\abs{x-y}^{2}}{4t}}\rm).$$ 
We denote by $q_c:=1+\frac{N}{2}$, the Brezis-Friedman critical exponent.

\bth {abovth}Let $q\geq q_c$. Then there exists a positive constant $C_1=C_1(N,q)$ such that for any closed subset $F$ of $\BBR^N$ and any 
$u\in C^{2}(Q_{\infty})\cap C(\overline {Q_{\infty}}\setminus F)$ satisfying 
\begin {equation}
\label {mainE} \BA {l}
\prt_{t}u-\Gd u+u^q=0\quad \mbox {in } Q_{\infty}\\[2mm]
\quad\;\;\;\displaystyle\lim_{t\to 0}u(x,t)=0\quad\mbox {locally uniformly in }F^c,
\EA
\end {equation}
there holds
\begin {equation}
\label {pott3}
u(x,t)\leq C_1  W_F(x,t)\qquad\forall (x,t)\in Q_\infty,
\end {equation}
where $ W_F$ is the $(2/q,q')$-parabolic capacitary potential of $F$ defined by $(\ref{pot1})$.
\es

 First we  consider the case where $F=K$ is compact and
      \begin {eqnarray}
\label {r}
 K\subset B_{r}\subset \overline 
 B_{r},
  \end {eqnarray}
  and then we  extend to the general case by a covering argument.
\subsection{Capacities and Besov spaces}

\subsubsection{$L^p$ regularity}

Throughout this paper $C$ will denote a generic positive constant, depending only on $N$, $q$ and sometimes $T$, the value of which may vary from one occurrence to another.  We  also use sometimes the notation $A\approx B$ for meaning that there exists a constant $C>0$ independent of the data such that
  $C^{-1}A\leq B\leq CA$.\smallskip

We recall some classical results dealing with $L^p$ capacities as they are developed in \cite{BP2}: if $1<p<\infty$
we denote 
\begin{equation}\label {W}
W^{2,1}_p(\BBR^{N+1}):=\{\gf\in L^p(\BBR^{N+1}):\prt_{t}\gf,\nabla \gf,D^2\gf\in L^p(\BBR^{N+1})\},
\end{equation}
with the associated norm
\begin{equation}\label {W1}
\norm{\gf}_{W^{2,1}_p}=\norm{\gf}_{L^p}+\norm{\nabla \gf}_{L^p}+\norm{\prt_t\gf}_{L^p}+\norm{D^2\gf}_{L^p}.
\end{equation}
We define a corresponding capacity on compact sets, that we extend it classicaly on capacitable sets.
\begin{equation}\label {Wt1}
C_{2,1,p}(E)=\inf\{\norm{\gf}_{W^{2,1}_p}:\gf\in C^{\infty}_0(\BBR^{N+1}):\gf\geq 1\text{ in a neighborhood of }E\},
\end{equation}
We extend the heat kernel $H$  in $\BBR^{N+1}=\{(x,t)\in\BBR^N\ti\BBR\}$ by assigning the value $0$ for $t<0$. Then, for any 
$\eta\in C_0(\BBR^N)$, 
\begin{equation}\label {W2}
\BBH[\eta](x,t)=\left\{\BA {ll}0\qquad&\text{if }t<0\\
H\ast(\eta\otimes\gd_{0})(x,t)\qquad&\text{if }t>0,
\EA\right.
\end{equation}
where $\gd_0$ has to be understood as the Dirac measure on $\BBR$ at $t=0$. For any subset $E\in\BBR^{N+1}$
\begin{equation}\label {W3}
C_{H,p}(E)=\inf\{\norm{f}_{L^p}:f\in L^p(\BBR^{N+1}),H\ast f\geq 1\text{ on }E\}.
\end{equation}
The following result is proved in \cite[Prop 2.1]{BP2}.

\bprop{BP1} For any $T>0$, there exists $c=c(T,p,N)$ such that
\begin{equation}\label {W4}
c^{-1}C_{H,p}(E)\leq C_{2,1,p}(E)\leq cC_{H,p}(E)\qquad\forall E\subset \BBR^N\ti]-T,T[,\,E\text{ Borel}.
\end{equation}
\es
We recall the Gagliardo Nirenberg inequality valid for any $\gf\in C^{\infty}_0(\BBR^{d})$
\begin{equation}\label {W4-1}
\norm{\nabla \gf}_{L^{2p}}^{2p}\leq c_{d,p}\norm{\gf}_{L^{\infty}}^{p}\norm{D^2\gf}_{L^{p}}^{p}.
\end{equation}
Furthermore, the trace at $t=0$ of functions in $W^{2,1}_p$ belongs to the Besov space $B^{2-\frac{2}{p},p}(\BBR^N)$. However, in our range of exponents  $B^{2-\frac{2}{p},p}(\BBR^N)=W^{2-\frac{2}{p},p}(\BBR^N)$. The reason for this is that $2-\frac{2}{p}$ is not an integer except if $p=2$, in which case equality holds also. If we set
\begin{equation}\label {W5}
c_{2-\frac{2}{p},p}(K)=\inf\{\norm{\gf}_{W^{2-\frac{2}{p},p}}:\gf\in C^\infty_0(\BBR^N),\gf\geq 1\text{ in a neighborhood of }K\}.
\end{equation}
then \cite[Prop 2.3]{BP2}.

\bprop{BP2} There exist $c=c(N,p)>0$ such that 
\begin{equation}\label {Wt5}
c^{-1}c_{2-\frac{2}{p},p}(E)\leq C_{2,1,p}(E\times \{0\})\leq cc_{2-\frac{2}{p},p}(E)\qquad\forall E\subset \BBR^{N},\,E\text{ Borel}.
\end{equation}
\es

The $c_{2-\frac{2}{p},p}$-capacity is equivalent to the Bessel capacity $C_{2-\frac{2}{p},p}$ defined by
 \begin{equation}\label {W6}
C_{2-\frac{2}{p},p}(E)=\inf\{\norm{f}_{L^p}:f\in L^p(\BBR^{N}),G_{2-\frac{2}{p}}\ast f\geq 1\text{ on }E\}
\end{equation}
where $G_{2-\frac{2}{p}}=\CF[(1+|\xi|^2)^{\frac{1}{p}-1}]$ denotes the Bessel kernel associated to the operator $(-\Gd +I)^{1-\frac{1}{p}}$.\medskip

\subsubsection{The Aronszajn-Slobodeckij integral}
If $\Gw$ is a domain in $\BBR^N$ and $0<s<1$, we denote by $\norm{.}_{\dot B^{s,p}(\Gw)}$ the Aronszajn-Slobodeckij norm defined on $C^\infty_0(\Gw)$ by
\begin{equation}\label {Be1}
\norm{\eta}_{\dot B^{s,p}}=\left(\dint_{\!\!\!\Gw\ti\Gw}\myfrac{|\eta(x)-\eta(y)|^p}{|x-y|^{N+sp}}dx dy\right)^{1/p}\qquad\forall \eta\in C^{\infty}_0(\Gw).
\end{equation}
In the case $1<s<2$, all the results which are presented still holds by replacing the function by its gradient.
We also consider the case $s=1$, but in our range of exponents the corresponding exponent for $p$ is $2$, in which case the space under consideration is just $H^1_0(\Gw)$.  Since the imbedding of $W^{1,p}(\Gw)$ is compact, it follows
  the imbedding of $B^{s,p}(\Gw)$ into $L^p(\Gw)$ is compact too. Therefore 
the following Poincar\'e type inequality holds \cite[p. 134]{Tar1}. Actually, the proof, obtained by contradiction,   is given with $W^{1,p}(\Gw)$ instead of $B^{s,p}(\Gw)$, but it depends only on the  compactness of the imbedding.  

\bprop{Poinca} Let $\Gw$ be a bounded domain and, $p\in (1,\infty)$ and $0<s\leq 1$ such that $sp\leq N$. Then there exists $\gl=\gl(\Gw,N,p)>0$ such that
\begin{equation}\label {poin1}
\dint_{\!\!\!\Gw\ti\Gw}\myfrac{|\eta(x)-\eta(y)|^p}{|x-y|^{N+sp}}dx dy\geq \gl\int_\Gw|\eta(x)|^p dx\quad\forall \eta\in C^{\infty}_0(\Gw).
\end{equation}
\es
\medskip

\noindent\Remark If $sp>N$, the  same proof re holds for all $\eta\in C^{\infty}_0(\Gw)$ (see the proof of \cite[Th 8.2]{DPV})
\begin{equation}\label {poin5}
\left(\dint_{\!\!\!\Gw\ti\Gw}\myfrac{|\eta(x)-\eta(y)|^p}{|x-y|^{N+s p}}dx dy\right)^{1/p}\geq C\frac{|\eta(z)-\eta(z')|}{|z-z'|^\ga}\quad\forall (z,z')\in\Gw\ti\Gw,\,z\neq z',
\end{equation}
with $\ga=s-N/p$ and $C=C(s,N,p)$. This estimate implies 
\begin{equation}\label {poin6}
\left(\dint_{\!\!\!\Gw\ti\Gw}\myfrac{|\eta(x)-\eta(y)|^p}{|x-y|^{N+s p}}dx dy\right)^{1/p}\geq Cd^{-\ga}\norm u_{L^\infty},
\end{equation}
where $d$ is the width of $\Gw$, i.e. the smallest of $\gd>0$ such that there exists an isometry $\CR$ such that  $\CR(\Gw)\subset D_\gd:=\{x=(x_1,x'):0<x_1<\gd\}$. \medskip

The related unpublished result due to L. Tartar \cite{Tar} will be useful in the sequel. We reproduce its proof for the sake of completeness.

\bprop{Poinca1} Assume $b>a$ and $\Gw\subset \Gg_{a,b}:=\{x=(x_1,x'):a<x_1<b\}$ is a domain. If $sp\leq N$ there exists $C=C(s,p,N,b/a)>0$ such that that 
\begin{equation}\label {poin7}
\dint_{\!\!\!\Gw\ti\Gw}\myfrac{|\eta(x)-\eta(y)|^p}{|x-y|^{N+sp}}dx dy\geq \gl\left(b-a\right)^{sp}\myint{\Gw}{}|\eta(x)|^p dx\quad\forall \eta\in C^{\infty}_0(\Gw).
\end{equation}
\es
\Proof Using the notation of \cite{LP}, $W^{s,p}(\BBR^N)$ is the interpolation space $[W^{1,p}(\BBR^N),L^p(\BBR^N)]_{s,p}$ and subset of $L^p(\BBR^{N-1};[W^{1,p}(\BBR),L^p{\BBR}]_{s,p})=L^p(\BBR^{N-1};W^{s,p}(\BBR))$, with continuous imbedding. Thus there exist $C>0$ such that
\begin{equation}\label {poin7-1}\BA {l}
\norm \eta_{L^p}^p+\myint{\BBR^{N-1}}{}\dint_{\BBR\ti\BBR}\myfrac{|\eta(x_1,x')-\eta(y_1,x')|^p}{|x_1-y_1|^{1+sp}}dx_1 dy_1dx'
\\[4mm]
\phantom{------}
\leq C\left(\norm \eta_{L^p}^p+\dint_{\!\!\!\BBR^N\ti\BBR^N}\myfrac{|\eta(x)-\eta(y)|^p}{|x-y|^{N+sp}}dx dy\right)
\EA\end{equation}
for all $\eta\in C^\infty_0(\BBR^N)$. This inequality is valid if $\eta$ is replaced by $\eta_\gt$ where $\eta_\gt(x)=\eta(\gt x)$ and  $\gt>0$. This gives
$$\BA {l}\norm \eta_{L^p}^p+\gt^{sp-N}\myint{\BBR^{N-1}}{}\dint_{\BBR\ti\BBR}\myfrac{|\eta(x_1,x')-\eta(y_1,x')|^p}{|x_1-y_1|^{1+sp}}dx_1 dy_1dx'\\[4mm]
\phantom{------}
\leq C\left(\norm \eta_{L^p}^p+\gt^{sp-N}\dint_{\!\!\!\BBR^N\ti\BBR^N}\myfrac{|\eta(x)-\eta(y)|^p}{|x-y|^{N+sp}}dx dy\right).
\EA$$
Letting $\gt\to 0$, we obtained
\begin{equation}\label {poin7-2}
\int_{\BBR^{N-1}}\dint_{\BBR\ti\BBR}\myfrac{|\eta(x_1)-\eta(y_1)|^p}{|x_1-y_1|^{1+sp}}dx_1 dy_1dx'
\leq C\dint_{\!\!\!\BBR^N\ti\BBR^N}\myfrac{|\eta(x)-\eta(y)|^p}{|x-y|^{N+sp}}dx dy\qquad\forall\eta\in C^{\infty}_0(\BBR^N).
\end{equation}
Using \rprop{Poinca} with $N=1$ we get
$$\myint{0}{1}\myint{0}{1}\myfrac{|\eta(x_1,x')-\eta(y_1,x')|^p}{|x_1-y_1|^{1+sp}}|dx_1 dy_1\geq
\gl \myint{0}{1}|\eta(x_1,x')|^pdx_1 \quad \forall \eta\in C^{\infty}_0 \left((0,1)\ti\BBR^{N-1}\right)
$$
for all $x'\in\BBR^{N-1}$. Using a standard change of scale, it transforms into
$$\myint{a}{b}\myint{a}{b}\myfrac{|\eta(x_1,x')-\eta(y_1,x')|^p}{|x_1-y_1|^{1+sp}}|dx_1 dy_1\geq
\gl (b-a)^{sp}\myint{a}{b}|\eta(x_1,x')|^pdx_1 \quad \forall \eta\in C^{\infty}_0 \left((a,b)\ti\BBR^{N-1}\right)
$$
Integrating over $\BBR^{N-1}$ and using $(\ref{poin7-2})$, we derive $(\ref{poin7})$.\qeda

\bdef {Besov} Assume $s\in (0,1)$ and $sp<1$ or $s=1$ and $p=2$. If $\Gw$ is any domain in $\BBR^N$, the Besov space $B_0^{s,p}(\Gw)$ is the closure of $C^{\infty}_0(\Gw)$ with respect to the norm
\begin{equation}\label {poin8}
\norm{\eta}_{B^{s,p}}=\norm{\eta}_{\dot B^{s,p}}+\norm{\eta}_{L^{p}}.
\end{equation}
\es

The following result is derived from \rprop{Poinca1}.

\bcor{Poinca2} Let $b>a>0$ and $\Gw$ be an open domain of $\BBR^N$ such that $\Gw\subset B_b\setminus\overline {B_a}$.
Then there exists a constant $C=C(s,p,N)>0$  such that for any $\eta\in C^\infty_0(\Gw)$
\begin{equation}\label {point8}
\norm{\eta}_{\dot B^{s,p}}\leq \norm{\eta}_{B^{s,p}}\leq C(b-a)^{sp}\norm{\eta}_{\dot B^{s,p}}.
\end{equation}
\es
\subsubsection{Heat potential and Besov space}

If $\eta\in C^{\infty}_0(\Gw)$, we extend it by $0$ outside $\Gw$ and set
    \begin {eqnarray}
    \label {Trt3}
\norm{\eta}_{\tilde B^{s,p}}=\left(\dint_{Q_\infty}\abs{t^{1-s/2}\prt_{t}\BBH[\eta]}^{p}dx\,\frac{dt}{t}\right)^{1/p}  
\end {eqnarray} 
  
  It is well known (see e.g. \cite{BeBu}) that the Besov space $B^{s,p}(\Gw)$ can be defined directly as the space of $\eta\in L^p(\Gw)$ functions such that $\norm{\eta}_{\dot B^{s,p}}<\infty$ or or such that $\norm{\eta}_{\tilde B^{s,p}}<\infty$. It coincides with the the interpolation space
  $\left[W^{2,p}(\Gw),L^p(\Gw)\right]_{s/2,p}$ (see \cite{LP}). Furthermore, there exists $C=C(s,p,N)>0$ such that 
     \begin {eqnarray}
    \label {Tr4-1}
C^{-1}\left(\norm \eta_{L^p}+\norm{\eta}_{ \dot B^{s,p}}\right)  \leq 
\norm \eta_{L^p}+\norm{\eta}_{\tilde B^{s,p}} \leq C\left(\norm \eta_{L^p}+\norm{\eta}_{ \dot B^{s,p}}\right)\quad\forall\eta\in B^{s,p}(\Gw).
\end {eqnarray}

  \blemma {Equiv}Assume $0<s<1$ and $1<p<\infty$ or $s=1$ and $p=2$. Then there exists a positive constant $C$, depending only on 
  $s,p,N$,  such that for any domain $\Gw$, there holds
     \begin {eqnarray}
    \label {Tr5}
C^{-1}\norm{\eta}_{ \dot B^{s,p}}  \leq 
\norm{\eta}_{\tilde B^{s,p}} \leq C\norm{\eta}_{ \dot B^{s,p}}\quad\forall\eta\in C^{\infty}_0(\Gw).
\end {eqnarray}  
    \es
\Proof Let $\eta\in C^{\infty}_0(\BBR^N)$ and $\gt>0$. Set $\eta_\gt(x)=\eta(\gt x)$, then $(\ref{Tr5})$ applied to $\eta_\gt$ yields to
$$C^{-1}\left(\norm \eta_{L^p}+\gt^{s}\norm{\eta}_{ \dot B^{s,p}}\right)  \leq 
\left(\norm \eta_{L^p}+\gt^{s}\norm{\eta}_{\tilde B^{s,p}}\right) \leq C\left(\norm \eta_{L^p}+\gt^{s}\norm{\eta}_{ \dot B^{s,p}}\right).
$$
Since it holds for any arbitrary large $\gt$ and  $\eta\in C^{\infty}_0(\BBR^N)$, $(\ref{Tr5})$ follows.\qeda\medskip

We denote by $\CT_{\Gw}(K)$ the set of functions $\eta\in C^\infty_0(\Gw)$ such that $0\leq\eta\leq 1$ and  $\eta=1\text{ on } K$. If $\Gw$ is a bounded subset of $\BBR^N$, we define the Besov capacity of a compact set $K\subset\Gw\subset\BBR^N$ by
       \begin {eqnarray}\label {Bt1}
R^\Gw_{s,p}(K)=\inf\{\norm{\eta}^p_{ \dot B^{s,p}}:\eta\in \CT_{\Gw}(K)\},       
  \end {eqnarray} 
  and the Bessel capacity relative to $\Gw$ by
   \begin {eqnarray}\label {B3}
C^\Gw_{s,p}(K)=\inf\{\norm{\eta}^p_{ B^{s,p}}:\eta\in \CT_{\Gw}(K)\}.       
  \end {eqnarray} 
We extend classicaly this capacity to any capacitable set $K\subset\Gw$. This capacity has the following scaling property.
\blemma{scal} For any $\gt>0$ and any capacitable set $K\subset\Gw$, there holds
       \begin {eqnarray}\label {B2}
R^{\Gw}_{s,p}(K)=\gt^{N-sp}R^{\gt^{-1}\Gw}_{s,p}(\gt^{-1}K).
  \end {eqnarray} 
  Furthermore, if $\Gw\subset B_b\setminus\overline{B_a}$, there exists $c=c(b-a,b/a, N,s,p)>0$ such that 
         \begin {eqnarray}\label {Bt2}
c^{-1}C^{\Gw}_{s,p}(K)\leq R^{\Gw}_{s,p}(K)\leq cC^{\Gw}_{s,p}(K).
  \end {eqnarray} 
  Finally, if $K\subset \Gw'\subset\overline{\Gw'}\subset\Gw$, there exists $c=c(N,s,p,\dist (\Gw',\Gw^c))$ such that
   \begin {eqnarray}\label {B4}
C_{s,p}(K)\leq C^{\Gw}_{s,p}(K)\leq cC_{s,p}(K).     
  \end {eqnarray} 
 \es
\Proof The scaling property $(\ref{B2})$ is clear by change of variable. Estimate $(\ref{Bt2})$ is a consequence of \rdef{Besov} and \rprop{Poinca1}. For the last statement, the left-hand side is obvious. For the right-hand side, consider a  smooth nonnegative cut-off function $\gz$ which is $1$ on $\overline{\Gw'}$, has value between $0$ and $1$ and has compact support in $\Gw$. If $\eta\in \CT_{\BBR^N}(K)$, $\gz\eta\in \CT_{\Gw}(K)$ and 
$$\BA {ll}\norm{\gz\eta}^p_{B^{s,p}}=\norm{\gz\eta}^p_{L^{p}(\Gw)}+\norm{\gz\eta}^p_{\dot B^{s,p}}
\\[2mm]
\phantom{\norm{\gz\eta}^p_{B^{s,p}}}\leq \norm{\eta}^p_{L^{p}(\Gw)}+\norm{\eta}^p_{\dot B^{s,p}}+\norm{\gz}^p_{\dot B^{s,\infty}}
\norm{\eta}^p_{L^{p}}
\\[2mm]
\phantom{\norm{\gz\eta}^p_{B^{s,p}}}
\leq c\norm{\eta}^p_{B^{s,p}},
\EA$$
where 
$$\norm{\gz}_{\dot B^{s,\infty}}=\sup_{x\neq y}\myfrac{|\gz(x)-\gz(y)|}{|x-y|^s}
$$
and $c\approx 1+(\dist (\Gw',\Gw^c))^{-s}$. The proof follows.
\qeda\medskip

In the sequel we assume that $q\geq q_c$ and we take $p=q'$ and $s=2/q$. If $K\subset \Gw$, $\Gw$ is bounded and $\eta\in \CT_{\Gw}(K)$, we set
       \begin {eqnarray}
\label {Tr2} R[\eta]=\abs {\nabla\BBH[\eta]}^{2}+\abs{\prt_{t}\BBH[\eta]}.
  \end {eqnarray}

\blemma {Lem2-1} There exists $C=C(N,q)>0$ such that for every $\eta\in \CT_{\Gw}(K)$
         \begin {eqnarray}
\label {Tr3}
\norm{\eta}^{q'}_{\tilde B^{2/q,q'}}\leq \dint_{Q_\infty}\left(R[\eta]\right)^{q'}dx\,dt:=\norm{R[\eta]}_{L^{q'}}^{q'}\leq C\norm{\eta}^{q'}_{\tilde B^{2/q,q'}}
  \end {eqnarray} 
  \es
 \Proof Using $(\ref{Trt3})$ and \rlemma {Equiv}, it follows from \rcor{Poinca2}  that
 $$\norm{\eta}^{q'}_{\tilde B^{2/q,q'}}\approx\iint_{Q_\infty}\abs{\prt_{t}\BBH[\eta]}^{q'}dxdt.$$
  Using the Gagliardo-Nirenberg inequality in $\BBR^N$, an elementary elliptic estimate and the fact that $0\leq \BBH[\eta]\leq 1$, we see that
  \begin {eqnarray}
\label {interpol2}
  \int_{\BBR^N}\abs{\nabla(\BBH[\eta](.,t))}^{2q'}dx
  \leq C\norm {D^2\BBH[\eta](.,t)}^{q'}_{L^{q'}}\norm {\BBH[\eta](.,t)}^{q'}_{L^\ity}\leq
  C\norm {\Gd\BBH[\eta](.,t)}^{q'}_{L^{q'}},
  \end {eqnarray} 
 for all $t>0$. Since  $\prt_{t}\BBH[\eta]=\Gd \BBH[\eta]$, it implies $(\ref{Tr3})$. 
  \qeda \medskip
  
  The dual space $B^{-2/q,q}(\Gw)$ of $B^{2/q,q'}(\Gw)$ is naturally endowed with the norm
  $$\norm{\gm}_{B^{-2/q,q}}=\sup\left\{\gm(\eta):\eta\in B^{2/q,q'}(\Gw),\norm{\eta}_{B^{2/q,q'}}\leq 1\right\}.
  $$
  
  The following result is may be already known, but we have not found it in the literature.  If $\gm$ is a bounded measure in $\BBR^N$, we denote by $\BBH[\gm]$ the solution of heat equation in $Q_\infty$ with initial data $\gm$. 

\blemma {meas} Assume $q\geq q_c$. For any $T>0$, there exist  a
constant $c>0$ such that,  for any bounded measure $\gm$ belonging to $B^{-2/q,q}(\BBR^{N})$, there holds 
 \begin {equation}
\label{mesest}
 c^{-1}{\norm{\gm}}_{B^{-2/q,q}(\BBR^{N})}\leq {\norm {\BBH[\gm]}}_{L^q(Q_{T})}
 \leq c{\norm{\gm}}_{B^{-2/q,q}(\BBR^{N})}.
 \end {equation}
Furthermore, if $q>q_c$ there holds
 \begin {equation}
\label{mesest0}
 c^{-1}{\norm{\gm}}_{B^{-2/q,q}(\BBR^{N})}\leq {\norm {\BBH[\gm]}}_{L^q(Q_{\infty})}
 \leq c{\norm{\gm}}_{B^{-2/q,q}(\BBR^{N})}+c{\norm{\gm}}_{\mathfrak M(\BBR^{N})}.
 \end {equation}.
\es
\Proof If $\gm\in B^{-2/q,q}(\BBR^{N})$, there exists a unique $\gw\in B^{2-2/q,q}(\BBR^{N})$ such that $\gm=(I-\Gd)\gw$, and
${\norm{\gm}}_{B^{-2/q,q}}\approx{\norm{\gw}}_{B^{2-2/q,q}}$. Applying standard interpolation methods to the analytic semi-group
$e^{-t(I-\Gd)}=e^{-t}e^{t\Gd}$ (see e.g. \cite {BeBu}, \cite{Tr}) we obtain, 
 \begin {equation}\label{mesest1}\BA {l}
 \left(\dint_{Q_\infty}\abs{t^{1/q}(I-\Gd)\BBH[\gw]}^qdx\myfrac{e^{-qt}dt}{t}\right)^{1/q}=\left(\dint_{Q_\infty}\abs{t^{1/q}\BBH[\gm]}^qdx\myfrac{e^{-qt}dt}{t}\right)^{1/q}\\
\phantom{\left(\dint_{Q_\infty}\abs{t^{1/q}(I-\Gd)\BBH[\gw]}^qdx\myfrac{e^{-qt}dt}{t}\right)^{1/q}}
\approx {\norm{\gw}}_{B^{2-2/q,q}}\\
\phantom{\left(\dint_{Q_\infty}\abs{t^{1/q}(I-\Gd)\BBH[\gw]}^qdx\myfrac{e^{-qt}dt}{t}\right)^{1/q}}
\approx{\norm{\gm}}_{B^{-2/q,q}}.
\EA \end {equation}
Clearly
$$e^{-qT}\dint_{Q_T}\abs{t^{1/q}\BBH[\gm]}^qdx\myfrac{dt}{t}\leq \dint_{Q_\infty}\abs{t^{1/q}\BBH[\gm]}^qdx\myfrac{e^{-qt}dt}{t},
$$
and
$$\BA {l}
\dint_{Q_\infty}\abs{t^{1/q}\BBH[\gm]}^qdx\myfrac{e^{-qt}dt}{t}=\displaystyle\sum_{n=0}^\infty\dint_{Q_{T+n+1}\setminus Q_{T+n}}\abs{t^{1/q}\BBH[\gm]}^qdx\myfrac{e^{-qt}dt}{t}
\\[4mm]\phantom{\dint_{Q_\infty}\abs{t^{1/q}\BBH[\gm]}^qdx\myfrac{e^{-qt}dt}{t}}
=\displaystyle\sum_{n=0}^\infty\dint_{Q_T}\abs{\BBH[\gm](s+n)}^qe^{-q(s+n)}ds
\\[4mm]\phantom{\dint_{Q_\infty}\abs{t^{1/q}\BBH[\gm]}^qdx\myfrac{e^{-qt}dt}{t}}
\leq \displaystyle\left(\sum_{n=0}^\infty e^{-qn}\right)\dint_{Q_T}\abs{t^{1/q}\BBH[\gm]}^q\myfrac{dt}{t}.
\EA
$$
This implies $(\ref{mesest})$. Furthermore, $\norm{\abs{\BBH[\gm](.,t)}}^q_{L^q}\leq ct^{-N(q-1)/2}\norm{\gm}^q_{\mathfrak M}$, thus 
$\BBH[\gm]\in L^q(Q_\infty)$ if $q>q_c$ (but this does not hold if $q=q_c$). If $q>q_c$ (equivalently $N(q-1)/2>1$), 

$$\BA{l}
\dint_{Q_\infty}\abs{t^{1/q}\BBH[\gm]}^qdx\myfrac{dt}{t}=\displaystyle\sum_{n=0}^\infty\dint_{Q_{T+n+1}\setminus Q_{T+n}}\abs{t^{1/q}\BBH[\gm]}^qdx\myfrac{dt}{t}\\[4mm]
\phantom{\dint_{Q_\infty}\abs{t^{1/q}\BBH[\gm]}^qdx\myfrac{dt}{t}}
=\dint_{Q_T}\abs{t^{1/q}\BBH[\gm]}^qdx\myfrac{dt}{t}+\displaystyle\dint_{Q_T}\sum_{n=1}^\infty\abs{\BBH[\gm](s+n)}^qdxds
\\[4mm]
\phantom{\dint_{Q_\infty}\abs{t^{1/q}\BBH[\gm]}^qdx\myfrac{dt}{t}}
\leq \dint_{Q_T}\abs{t^{1/q}\BBH[\gm]}^qdx\myfrac{dt}{t}+C\left(\displaystyle\sum_{n=1}^\infty n^{-N(q-1)/2}\right)\norm\gm^q_{\mathfrak M}.
\EA$$
Thus we obtain  $(\ref{mesest0})$.\qeda\medskip



 \subsection {Global $L^q$-estimates}
 Let $\gr>0$, we assume $(\ref {r})$ holds. With the previous notations, $\CT_{r,{r+\gr}}(K)$ denotes the set of functions $\eta\in C_{0}^\ity(B_{r+\gr})$, such that $0\leq \eta\leq 1$ and value 1 on $K$.  If $\eta\in \CT_{r,\gr}(K)$, we set 
 $$\eta^{*}=1-\eta\;\text{ and }\,
 \gz=\BBH[\eta^{*}]^{2q'}.$$

 \blemma {Lem2-2} Assume $u$ is a positive solution of $(\ref{mainE})$ in $Q_\infty$. There exists $C=C(N,q)>0$ such that for every $T>0$ and every compact set $K\subset B_r$, 
         \begin {eqnarray}
\label {Tr4}
\dint_{Q_T}u^q\gz dx\,dt+\int_{\BBR^{N}}(u\gz ) (x,T)dx\leq 
C{\norm 
{R[\eta]}}^{q'}_{L^{q'}}\qquad\forall\eta\in \CT_{r,\gr}(K).
  \end {eqnarray}
    \es
\Proof We recall that there always holds
  \begin {eqnarray}
\label {OK1} 0\leq u(x,t)\leq \left(\myfrac 
{1}{t(q-1)}\right)^{\frac{1}{q-1}}\qquad\forall (x,t)\in Q_{\infty},
  \end {eqnarray}
  and
  \begin {eqnarray}
\label {BF} 0\leq u(x,t)\leq \left(\myfrac 
{C}{t+(\abs x-r)^2}\right)^{\frac{1}{q-1}}\qquad\forall (x,t)\in Q_{\infty}\setminus 
B_{r}\ti \BBR,
  \end {eqnarray}  
by the Brezis-Friedman estimate \cite {BF}. Since $\eta^{*}$ vanishes in an open neighborhood $\CN_{1}$, for any 
open subset $\CN_{2}$ such that 
$K\subset \CN_{2}\subset \overline \CN_{2}\subset \CN_{1}$ there 
exist $c_2=c_{_{\CN_{2}}}>0$  and $C_2=C_{_{\CN_{2}}}>0$ such that 
$$
\BBH[\eta^{*}](x,t)\leq C_2e^{-\frac{c_2}{t}},\qquad\forall 
(x,t)\in Q_{T}^{\CN_{2}}. $$
Therefore $$
\lim_{t\to 0}\int_{\BBR^{N}}(u\gz)(x,t)dx=0. $$
Thus $\gz$ is an admissible test function and one has
  \begin {eqnarray}
\label {OK2}
\dint_{Q_{T}}u^q\gz dx\,dt+\int_{\BBR^{N}}(u\gz ) (x,T)dx =\dint_{Q_{T}}u(\prt_{t}\gz
+\Gd\gz)dx\,dt.
  \end {eqnarray}
  Notice that the two terms on the left-hand side are nonnegative. 
  Put $\mathbb H_{\eta^{*}}=\mathbb H[\eta^{*}]$, then
\begin {eqnarray*}
\prt_{t}\gz+\Gd\gz 
&=&2q'\mathbb H_{\eta^{*}}^{2q'-1}\left(\prt_{t}\mathbb H_{\eta^{*}}+\Gd\mathbb
H_{\eta^{*}}
\right)+2q'(2q'-1)\mathbb H_{\eta^{*}}^{2q'-2} {\abs{\nabla{\mathbb H_{\eta^{*}}}}}^{2},\\
&=&2q'\mathbb H_{\eta^{*}}^{2q'-1}\left(\prt_{t}\mathbb H_{\eta}+
\Gd\mathbb H_{\eta}
\right)+2q'(2q'-1)\mathbb H_{\eta}^{2q'-2} {\abs{\nabla{\mathbb H_{\eta}}}}^{2},\\
\end{eqnarray*}
because $\mathbb H_{\eta^{*}}=1-\mathbb H_{\eta}$, hence $$
u(\prt_{t}\gz+\Gd\gz ) =u\mathbb H_{\eta^{*}}^{2q'/q}
\left[2q'(2q'-1)\mathbb H_{\eta^{*}}^{2q'-2-2q'/q}{\abs{\nabla{\mathbb H_{\eta}}}}^{2}
-2q'\mathbb H_{\eta^{*}}^{2q'-1-2q'/q} (\Gd\mathbb H_{\eta}+\prt_{t}\mathbb
H_{\eta})\right] .$$
Finally, since $2q'-2-2q'/q=0$ and $0\leq \mathbb H_{\eta^{*}}\leq 1$, there holds
$$
\abs{\dint_{Q_{T}}u(\prt_{t}\gz+\Gd\gz)dx\,dt}
\leq C(q)\left(\dint_{Q_{T}}u^q\gz dx\,dt\right)^{1/q} 
\left(\dint_{Q_{T}}R^{q'}(\eta)dx\,dt\right)^{1/q'}, $$
where $$
R(\eta)=\abs{\nabla{\mathbb H_{\eta}}}^{2}+
\abs{\Gd\mathbb H_{\eta}+\prt_{t}\mathbb H_{\eta}}. $$
Using \rlemma {Lem2-1} one obtains $(\ref {Tr4})$.\qeda \medskip

\bprop {Prop2-3} Under the assumptions of \rlemma{Lem2-2}, let $r>0$, $\gr>0$, $T\geq(r+\gr)^{2}$
$$
\CE_{r+\gr}:=\{(x,t):{\abs x}^{2}+t\leq (r+\gr)^{2}\} $$
and $Q_{r+\gr,T}=Q_{T}\setminus \CE_{r+\gr}$. There exists $C=C(N,q,T)>0$ such that 
  \begin {eqnarray}
\label {OK4}
\dint_{Q_{r+\gr,T}}u^q dx\,dt+\int_{\BBR^{N}}u (x,T)dx
\leq C{\norm 
{R[\eta]}}^{q'}_{L^{q'}}\qquad\forall\eta\in \CT_{r,\gr}(K).
  \end {eqnarray}
\es
\Proof In view of \rlemma{Lem2-2} we only have to show that there exists a positive constant 
$c(N,q)$ such that, for $\eta$ as above and $T\geq (r+\gr)^2$, 
$$\gz=\BBH{\eta^*}^{2q'}>c(N,q).$$
Since, by assumption $K\subset B_{r}$, $\eta^{*}\equiv 1$  outside 
 $B_{r+\gr}$ and $0\leq\eta^*\leq 1$,
\begin {eqnarray*}
\BBH[\eta^{*}](x,t)\geq 
\BBH[1-\chi_{_{B_{r+\gr}}}](x,t)&=&
\left(\myfrac {1}{4\gp t}\right)^{\frac{N}{2}}\int_{\abs y\geq r+\gr}e^{-\frac{{\abs 
{x-y}}^{2}}{4t}}dy,\\
&=&1-\left(\myfrac 
{1}{4\gp t}\right)^{\frac{N}{2}}\int_{\abs y\leq r+\gr}e^{-\frac{{\abs 
{x-y}}^{2}}{4t}}dy.
\end {eqnarray*}
For $(x,t)\in Q_{r+\gr,T}$, put $x=(r+\gr)\xi$, $y=(r+\gr)\gu$ and $t=(r+\gr)^2\gt$. Then 
$(\xi,\gt)\in Q_{1,\frac{T}{(r+\gr)^2}}$ and
\begin {eqnarray*}
\left(\myfrac 
{1}{4\gp t}\right)^{\frac{N}{2}}\int_{\abs y\leq r+\gr}e^{-\frac{{\abs 
{x-y}}^{2}}{4t}}dy=\left(\myfrac 
{1}{4\gp \gt}\right)^{\frac{N}{2}}\int_{\abs {\gu}\leq 1}e^{-\frac{{\abs 
{\xi-\gu}}^{2}}{4\gt}}d\gu.
\end {eqnarray*}
We claim that
\begin {eqnarray}
\label {H3}
\max \left\{
\left(\myfrac {1}{4\gp \gt}\right)^{\frac{N}{2}}
\int_{\abs {\gu}\leq 1}e^{-\frac{{\abs 
{\xi-\gu}}^{2}}{4\gt}}d\gu:(\xi,\gt)\in Q_{1,\frac{T}{(r+\gr)^2}}\right\}=\ell ,
\end {eqnarray}
for some $\ell=\ell (N,\frac{T}{(r+\gr)^2})\in (0,1]$, and $\ell$ is actually independent of $\frac{T}{(r+\gr)^2}$Êif this quantity is larger than $1$. 
We recall that
\begin {eqnarray}
\label {H2}
\left(\myfrac 
{1}{4\gp \gt}\right)^{\frac{N}{2}}\int_{\abs {\gu}\leq 1}e^{-\frac{{\abs 
{\xi-\gu}}^{2}}{4\gt}}d\gu<1\qquad\forall \gt>0.
\end {eqnarray}
If the maximum is achieved for 
some $(\bar\xi,\bar\gt)\in Q_{1,\frac{T}{(r+\gr)^2}}$, it is smaller than $1$ and
\begin {eqnarray}
\label {H4}
\BBH[\eta^{*}](x,t)\geq 
\BBH[1-\chi_{_{B_{r+\gr}}}](x,t)\geq 1-\ell>0,\qquad\forall (x,t)\in  
Q_{r+\gr,T}.
\end {eqnarray}
Let us assume that the maximum is achieved following a sequence $\{(\xi_{n},\gt_{n})\}$ 
with $\gt_{n}\to 0$ and $\abs{\xi_n}\to \ga\geq 1$. Then
$$
\left(\myfrac {1}{4\gp \gt_n}\right)^{\frac{N}{2}}
\int_{\abs {\gu}\leq 1}e^{-\frac{\abs 
{\xi_n-\gu}^{2}}{4\gt_n}}d\gu=
\left(\myfrac {1}{4\gp \gt_n}\right)^{\frac{N}{2}}\int_{B_1(\xi_n)}e^{-\frac{\abs 
{\gu}^{2}}{4\gt_n}}d\gu\leq \myfrac{1}{2}.
$$
To verify this, note that $B_1(\xi_n)\cap B_1(-\xi_n)=\emptyset$, so that
$$\int_{B_1(\xi_n)}e^{-\frac{\abs 
{\gu}^{2}}{4\gt_n}}d\gu+\int_{B_1(-\xi_n)}e^{-\frac{\abs 
{\gu}^{2}}{4\gt_n}}d\gu<\int_{\BBR^N}e^{-\frac{\abs 
{\gu}^{2}}{4\gt_n}}d\gu<1
$$
and
$$\int_{B_1(\xi_n)}e^{-\frac{\abs 
{\gu}^{2}}{4\gt_n}}d\gu=\int_{B_1(-\xi_n)}e^{-\frac{\abs 
{\gu}^{2}}{4\gt_n}}d\gu.
$$
If the supremum is achieved with a sequence 
$\{(\xi_{n},\gt_{n})\}$ such that $|\xi_n|\to \infty$, the same argument applies. Finally if $\{\gx_n\}$ is bounded but $\tau_n\to \infty$ then the expression in $(\ref{H2})$ tends to zero.
Therefore $(\ref {H2})$ holds. Put $C=(1-\ell)^{-1}$, then
\begin {equation}
\label {capest3}
\dint_{Q_{r,T}}u^qdx\,dt+\myint{\BBR^{N}}{}u(.,T)dx
\leq C\norm{R[\eta]}_{L^{q'}}^{q'},
\end {equation}
and 
$(\ref {OK4})$ follows.\qeda 

\subsection {Pointwise estimates} 
In this subsection $u$ is a positive solution of $(\ref{mainE})$ in $Q_\infty$ and the assumptions of  \rlemma{Lem2-2} hold. We first derive a rough pointwise estimate.

\blemma{RPW} There exists a constant $C=C(N,q)>0$ such that, for any $\eta\in \CT_{r,\gr}(K)$,
\begin {equation}
\BA {r}\label {rpw1} u(x,(r+2\gr)^{2})\leq\myfrac {C
\norm{R[\eta]}_{L^{q'}}^{q'}}{(\gr(r+\gr))^{\frac{N}{2}}},
\qquad\forall 
x\in \BBR^{N}.
\EA 
\end {equation}
\es 
\Proof We recall that
\begin {equation}
\label {rpw2}
\int_{s}^{T}\int_{\BBR^{N}}u^qdx\,dt+\int_{\BBR^{N}}u(x,T)dx =
\int_{\BBR^{N}}u(x,s)dx\qquad\forall T> s>0,
\end {equation}
and
\begin {equation}
\label {rpw3}
\int_{\BBR^{N}}u(.,s)dx\leq C\norm{R[\eta]}_{L^{q'}}^{q'}\qquad\forall T> s\geq (r+\gr)^2,
\end {equation}
by \rprop{Prop2-3}. Using the fact that 
$$
u(x,\gt+s)\leq \BBH[u(.,s)](x,\gt)\leq \left(\frac {1}{4\gp 
\gt}\right)^{\frac{N}{2}}\int_{\BBR^{N}}u(.,s)dx, $$
$(\ref {rpw1})$ follows from $(\ref{rpw3})$ with $s=(r+\gr)^2$ and $\gt=(r+2\gr)^2-(r+\gr)^2\approx \gr(r+\gr)$.\qeda \medskip

The above estimate does not take into account the 
fact that $u(x,0)=0$ if $\abs x\geq r$. It is mainly interesting if 
$\abs x\leq r$. In order to derive a sharper estimate which takes this fact into account, we need some lateral boundary 
estimates.
\blemma {LATESLE} Let $\gg\geq r+2\gr$ and $c>0$ and either $N=1$ or $2$ and
$0\leq t\leq 
c\gg^{2}$ for some $c>0$, or $N\geq 3$ and $t>0$. Then, for any $\eta\in\CT_{r,\gr}(K)$, there holds 
\begin {equation}
\label {latest}
\int_{0}^{t}\int_{\prt B_{\gg}}udSd\gt\leq C_{5} \gg \norm{R[\eta]}_{L^{q'}}^{q'}.
\end {equation}
where $C>0$ depends on $N$, $q$ and $c$ if $N=1,\,2$ or depends only on $N$ and $q$ if $N\geq 3$.\es 
\Proof First we assume $N=1$ or $2$. Put $G^{\gg }:= B^c_{\gg }\ti 
(-\ity,0)$ and $\prt_{\ell}G^{\gg }=\prt B_{\gg }\ti 
(-\ity,0)$.
 We set 
 $$
h_{\gg}(x)=1-\frac {\gg }{\abs x},$$ 
 and let $\psi_{\gg}$ be 
the solution of
\begin {equation}
\label {testpsi}\BA {rc}
\prt_{\gt}\psi_{\gg}+\Gd\psi_{\gg}=0\,\;&\quad\mbox {in }G^{\gg }, \\
[2mm]
\psi_{\gg}=0\,\;&\quad\mbox {on }\prt_{\ell}G^{\gg }, \\
[2mm]
\psi_{\gg}(.,0)=h_{\gg}&\quad\mbox {in }B^c_{\gg }.
\EA 
\end {equation}
Thus the function $$
\tilde\psi (x,\gt)=\psi_{\gg}(\gg x,\gg^{2}\gt) $$
satisfies
\begin {equation}
\label {testpsi*}\BA {rc}
\prt_{t}\tilde\psi+\Gd\tilde\psi=0\,&\quad\mbox {in }G^{1} \\
[2mm]
\tilde\psi=0\,&\quad\mbox {on }\prt_{\ell}G^{1} \\
[2mm]
\tilde\psi(.,0)=\tilde h&\quad\mbox {in }B^c_{1},
\EA 
\end {equation}
and $\tilde h (x)=1-{\abs x}^{-1}$. By the maximum principle $0\leq 
\tilde \psi\leq 1$, and by Hopf Lemma
\begin {equation}
\label {testpsi*2} -\myfrac {\prt\tilde \psi}{\prt\bf 
n}\vline_{\prt B_{1}\ti [-c,0]}\geq \gth>0,
\end {equation}
where $\gth=\gth (N,c)$. Then $0\leq\psi_{\gg}\leq 1$ and
\begin {equation}
\label {testpsi*3} -\myfrac {\prt \psi_{\gg}}{\prt\bf n}\vline_{\prt B_{\gg }\ti 
 [-\gg^{2},0]}\geq 
\gth/\gg.
\end {equation}
Multiplying $(\ref {mequ})$ by $\psi_{\gg}(x,\gt-t)=\psi_{\gg}^{*}(x,\gt)$ and integrating on
$B^c_{\gg}\ti 
 (0,t)$ 
yields to
\begin {equation}
\BA {r}\label {rpw4}
\myint{0}{t}\myint {B^c_{\gg }}{}u^q\psi^{*}_{r}dxd\gt+\myint{B^c_{\gg 
}}{}(uh_{\gg})(x,t)dx -\myint{0}{t}\myint {\prt B_{\gg }}{}
\myfrac {\prt u}{\prt\bf n}\psi^{*}_{\gg}dSd\gt= 
-\myint{0}{t}\myint{\prt B_{\gg }}{}\myfrac {\prt \psi^{*}_{\gg}}{\prt\bf 
n}ud\gs d\gt.
\EA
\end {equation}
Since $\psi^{*}_{\gg}$ is bounded from above by $1$, estimate $(\ref {latest})$ follows from 
$(\ref {testpsi*3})$ and \rprop {Prop2-3} (notice that $B^c_{\gg }\ti 
 (0,t)\subset \CE^c_{\gg}$), first by taking $t=T= \gg^2\geq (r+2\gr)^2$, and then for any $t\leq \gg^2$.\medskip
 
 \noindent  If $N\geq 3$, we proceed as above except that we take
 $$h_\gg(x)=1-\left(\myfrac{\gg}{\abs x}\right)^{N-2}.
 $$
 Then $\psi_\gg(x,t)=h_\gg(x)$ and $\gth=N-2$ is independent of the length of the time interval. This leads to the conclusion.
 \qeda 
\blemma {decaylem1}I- Let $M,\,a>0$ and $\eta\in L^{\ity}(\BBR^{N})$ such that
\begin {equation}
\label {decay} 0\leq \eta(x)\leq Me^{-a{\abs x}^{2}}\qquad a.e.\mbox { in }\BBR^{N}.
\end {equation}
Then, for any $t> 0$, 
\begin {equation}
\label {decay1} 0\leq \BBH[\eta](x,t)\leq \frac {M}{(4at+1)^{\frac{N}{2}}} e^{-\frac{a{\abs
x}^{2}}{4at+1}}\qquad\forall x \in \BBR^{N}.
\end {equation}
II- Let $M,\,a,\,b>0$ and $\eta\in L^{\ity}(\BBR^{N})$ such that
\begin {equation}
\label {decay2} 0\leq \eta(x)\leq Me^{-a(\abs {x}-b)_{_{+}}^{2}}\qquad a.e.\mbox { in
}\BBR^{N}.
\end {equation}
Then, for any $t> 0$, 
\begin {equation}
\label {decay3} 0\leq \BBH[\eta](x,t)\leq \frac {Me^{-\frac{a(\abs
{x}-b)_{_{+}}^{2}}{4at+1}}}{(4at+1)^{\frac{N}{2}}} \qquad\forall x \in \BBR^{N},\,\forall t>0.
\end {equation}
\es 
\Proof For the first statement, put $a=\frac{1}{4}s$. Then 
$$
0\leq \eta (x)\leq M(4\gp s)^{\frac{N}{2}}\frac {1}{(4\gp s)^{\frac{N}{2}}} e^{-\frac{\abs x^{2}}{4s}}=C(4\gp
s)^{\frac{N}{2}}\BBH[\gd_{0}](x,s). $$
By the order property of the heat kernel, $$
0\leq \BBH[\eta](x,t)\leq M(4\gp 
s)^{\frac{N}{2}}\BBH[\gd_{0}](x,t+s)=M\left(\myfrac {s}{t+s}\right)^{\frac{N}{2}} e^{-\frac{\abs x^{2}}{4(t+s)}}, $$
and $(\ref {decay1})$ follows by replacing $s$ by $\frac{1}{4}a$.\medskip

\noindent For the second statement, let $\tilde a<a$ and 
$R=\max\{e^{-a(r-b)_{_{+}}^{2}+\tilde a r^{2}}:r\geq 0\}$. A direct computation gives $
R=e^{\frac{a\tilde a b^{2}}{a-\tilde a}}$, and $(\ref {decay3})$ implies $$
0\leq \eta(x)\leq Me^{\frac{a\tilde a b^{2}}{a-\tilde a}}e^{-\tilde a 
{\abs x}^{2}}. $$
Applying the statement I, we derive
\begin {equation}
\label {decay4} 0\leq \BBH[\eta](x,t)\leq \frac {Ce^{\frac{a\tilde a b^{2}}{a-\tilde a}}}{(4\tilde at+1)^{\frac{N}{2}}} e^{-\frac{\tilde a{\abs x}^{2}}{4\tilde at+1}}\qquad\forall x \in
\BBR^{N}, \;
\forall t>0.
\end {equation}
Since for any $x\in\BBR^{N}$ and $t>0$, 
$$
(4\tilde at+1)^{-\frac{N}{2}}e^{-\frac{\tilde a{\abs x}^{2}}{4\tilde at+1}}\leq 
e^{-\frac{a\tilde a b^{2}}{a-\tilde a}}(4 at+1)^{-\frac{N}{2}}e^{-\frac{ a(\abs x-b)^{2}}{4 
at+1}}, $$
$(\ref {decay3})$ follows from $(\ref{decay4})$.\qeda 
\blemma {decaylem3} There exists a 
constant $C=C(N,q)>0$ such that, for any $\eta\in\CT_{r,\gr}(K)$, there holds
\begin {equation}
\label {decay5} u(x,(r+2\gr)^{2})\leq 
C
\max\left\{\myfrac{r+\gr}{(\abs x-r-2\gr)^{N+1}},\myfrac{\abs x-r-2\gr}{(r+\gr)^{N+1}}\right\}
e^{-\frac{(\abs x-(r+2\gr))^{2}}{4(r+2\gr)^{2}}}\norm{R[\eta]}_{L^{q'}}^{q'},
\end {equation}
for any $ x \in \BBR^{N}\setminus B_{r+3\gr}$.
\es
\Proof It is classical that the Dirichlet 
heat kernel $H^{B^c_{1}}$ in the 
complement of $B_{1}$ satisfies, for some $C=C(N)>0$, 
\begin {equation}\label {dhk} 
H^{B^c_{1}}(x',y',t',s')\leq C_{7}(t'-s')^{-(N+2)/2}(|x'|-1)e^{-{\frac{|x'-y'|^2}{4(t'-s')}}},
\end {equation}
for $t'>s'$. 
By performing the change of variable $x'\mapsto (r+2\gr)x'$, $t'\mapsto 
 (r+2\gr)^{2}t'$, 
for any $x\in \BBR^{N}\setminus B_{ r+2\gr}$ and $0\leq t\leq T$, one obtains
\begin {equation}
\label {dhk1} u(x,t)\leq C(\abs x-r-2\gr)\int_{0}^t\int_{\prt 
B_{r+2\gr}}\myfrac {e^{-\frac{\abs 
{x-y}^{2}}{4(t-s)}}}{(t-s)^{1+\frac{N}{2}}}u(y,s)d\gs (y)ds.
\end {equation}
The right-hand side term in $(\ref {dhk1})$ is smaller than $$
\max\left\{\myfrac {C(\abs x-r-2\gr)}{(t-s)^{1+\frac{N}{2}}}e^{-\frac{(\abs x-r-2\gr)^{2}}{4(t-s)}}: s\in
(0,t)\right\}
\int_{0}^t\int_{\prt B_{r+2\gr}}u(y,s)d\gs (y)ds. $$
We fix $t=(r+2\gr)^{2}$ and $\abs x\geq r+3\gr$. Since 
\begin {eqnarray*}
&&\max\left\{\myfrac {e^{-\frac{(\abs x-r-2\gr)^{2}}{4s}}}{s^{1+\frac{N}{2}}}: s\in 
\left(0,(r+2\gr)^{2}\right)\right\}\\
&&\qquad\qquad\qquad= (\abs x-r-2\gr)^{-2-N}\max\left\{\myfrac {e^{-\frac{1}{4\gs}}}{\gs^{1+\frac{N}{2}}}: 
0<\gs<\left(\myfrac{r+2\gr}{\abs x-r-2\gr}\right)^{2}\right\},
\end {eqnarray*}
a direct computation gives
\begin {eqnarray*}
&&\max\left\{\myfrac {e^{-\frac{1}{4}\gs}}{\gs^{1+\frac{N}{2}}}: 
0<\gs<\left(\myfrac{r+2\gr}{\abs x-r-2\gr}\right)^{2}\right\}\\[2mm]
&&\phantom{---}=\left\{\BA {ll} (2N+4)^{1+\frac{N}{2}}e^{-(N+2)/2}&\mbox { if }r+3\gr\leq \abs x\leq
(r+2\gr)(1+\sqrt {4+2N}) ,\\
[2mm]
\left(\myfrac {\abs x-r-2\gr}{r+2\gr}\right)^{2+N}e^{-\left(\frac{\abs x-r-2\gr}{2r+4\gr}\right)^{2}}&\mbox { if } \abs x\geq
(r+2\gr)(1+\sqrt {4+2N}).
\EA\right.
\end {eqnarray*}
Thus there exists a constant $C(N)>0$ such that 
\begin {eqnarray}
\label {decay6*}
\max\left\{\myfrac {e^{-\frac{(\abs x-r-2\gr)^{2}}{4s}}}{s^{1+\frac{N}{2}}}: s\in 
\left(0,(r+2\gr)^{2}\right)\right\}\leq C(N)\gr^{-2-N}e^{-\left(\frac{\abs x-(r+2\gr)}{2r+4\gr}\right)^{2}}.
\end{eqnarray}
Combining this estimate with $(\ref {latest})$ with $\gg=r+2\gr$ and $(\ref {dhk1})$, 
one derives $(\ref {decay5})$.\qeda \\
\blemma {decaylem4} There exists a 
constant $C=C(N,q)>0$ such that
\begin {equation}
\label {decay7a} 0\leq u(x,(r+2\gr)^{2})\leq C
\max\left\{\myfrac{(r+\gr)^3}{\gr(\abs x-r-2\gr)^{N+1}},\myfrac{1}{(r+\gr)^{N-1}\gr}\right\}
e^{-\left(\frac{\abs x-r-3\gr}{2r+4\gr}\right)^{2}}
\norm{R[\eta]}_{L^{q'}}^{q'},
\end {equation}
for every $ x \in \BBR^{N}\setminus B_{r+3\gr}$.
\es
\Proof This is a direct consequence of the inequality
\begin {equation}\label{x}
(\abs x-r-2\gr)e^{-\left(\frac{\abs x-r-2\gr}{2r+4\gr}\right)^{2}}\leq \frac 
{C(r+\gr)^2}{\gr}e^{-\left(\frac{\abs x-r-3\gr}{2r+4\gr}\right)^{2}}, 
\qquad\forall x\in B^c_{r+2\gr}, \end {equation}
and \rlemma {decaylem3}.
\qeda 
\blemma {decaylem5} 
There exists a constant $C=C(N,q)>0$ such that, for any $\eta\in \CT_{r,\gr}(K)$, the following estimate holds
\begin {equation}
\label {decay6} u(x,t)\leq 
\myfrac {C\tilde Me^{-\frac{(\abs x-r-3\gr)_{_{+}}^{2}}{4t}}}{t^{\frac{N}{2}}} \norm{R[\eta]}_{L^{q'}}^{q'},
\qquad\forall x \in \BBR^{N},\,\forall t\geq 
(r+2\gr)^{2},
\end {equation}
where
\begin {equation}
\label {decay6''} \tilde M=\tilde M(x,r,\gr)=\left\{\BA{ll}
\left(1+\frac{r}{\gr}\right)^{\frac{N}{2}}\qquad\qquad\qquad&\mbox { if }\abs x<r+3\gr\\[2mm]
\frac{(r+\gr)^{N+3}}{\gr(\abs x-r-2\gr)^{N+2}}\quad&\mbox { if }r+3\gr\leq \abs x\leq c^*_N(r+2\gr)\\[2mm]
1+\frac{r}{\gr}\qquad\qquad\qquad&\mbox { if }\abs x\geq c^*_N(r+2\gr)
\EA\right.\end {equation}
with $c^*_N=1+\sqrt {4+2N}$.
\es
\Proof It follows by the maximum principle 
$$
u(x,t)\leq \BBH[u(.,(r+2\gr)^{2})](x,t-(r+2\gr)^{2}). 
$$
for $t\geq (r+2\gr)^{2}$ and $x\in\BBR^{N}$. By \rlemma {RPW} and \rlemma {decaylem4}
$$
u(x,(r+2\gr)^{2})\leq C_{10}\tilde Me^{-\frac{(\abs x-r-3\gr)^{2}}{4(r+2\gr)^{2}}}\norm{R[\eta]}_{L^{q'}}^{q'}, $$
where
$$\tilde M=\left\{\BA{ll}
((r+\gr)\gr)^{-\frac{N}{2}}\qquad\qquad\qquad&\mbox { if }\abs x<r+3\gr\\[2mm]
\frac{(r+\gr)^3}{\gr}\left(\abs x-r-2\gr)\right)^{N+2}\qquad\quad&\mbox { if }r+3\gr\leq \abs x\leq c^*_N(r+2\gr)\\[2mm]
\frac{1}{(r+\gr)^{N-1}\gr}\quad\qquad\qquad&\mbox { if }\abs x\geq c^*_N(r+2\gr)
\EA\right.$$
Applying \rlemma {decaylem1} with $a=(2r+4\gr)^{-2}$, $b=r+3\gr$ and $t$ 
replaced by $t-(r+2\gr)^{2}$ implies 
\begin {equation}
\label {decay6'} u(x,t)\leq C \myfrac {(r+2\gr)^N\tilde M}{t^{\frac{N}{2}}}
e^{-\frac{(\abs x-r-3\gr)^{2}}{4t}} \norm{R[\eta]}_{L^{q'}}^{q'},
\end {equation}
for all $ x \in B_{r+3\gr}^c$ and $ t\geq (r+2\gr)^{2}$, which is
$(\ref {decay6})$.\qeda\medskip

The next estimate gives a precise upper bound for $u$ when $t$ is not 
bounded from below.
\blemma {decaylem6} Assume that $0< t\leq (r+2\gr)^2$, then there exists a constant $C=C(N,q)>0$ such that the following
estimate holds
\begin {equation}
\label {decay7} u(x,t)\leq C(r+\gr)
\max\left\{\myfrac{1}{(\abs x-r-2\gr)^{N+1}},\myfrac{1}{\gr t^{\frac{N}{2}}}\right\} 
e^{-\frac{(\abs x-r-3\gr)^{2}}{4t}} \norm{R[\eta]}_{L^{q'}}^{q'},
\end {equation}
for any $ (x,t) \in \BBR^{N}\setminus
B_{r+3\gr}\ti (0,(r+2\gr)^{2}]$.\es
\Proof Thanks to $(\ref {latest})$ the following estimate is a 
straightforward variant of $(\ref {decay5})$ for any $|x|\geq r+2\gr$,
\begin {eqnarray}
\label {decay5*} u(x,t)\leq C_{8}(\abs x-r-2\gr)(r+2\gr)
\max\left\{\myfrac{e^{-\frac{(\abs x-r-2\gr)^{2}}{4s}} }{s^{1+\frac{N}{2}}}:0<s\leq t\right\}
\norm{R[\eta]}_{L^{q'}}^{q'}.
\end {eqnarray}
Clearly
$$\BA{l}\max\left\{\myfrac{e^{-\frac{(\abs x-r-2\gr)^{2}}{4s}} }{s^{1+\frac{N}{2}}}:0<s\leq t\right\}\\[6mm]\qquad
\phantom{-}
=\left\{\BA {ll}(2N+4)^{1+\frac{N}{2}}(\abs x-r-2\gr)^{-N-2}e^{-\frac{N+2}{2}}&\mbox {if }0<\abs x\leq r+2\gr+\sqrt{2t(N+2)}\\[2mm]
\myfrac{e^{-\frac{(\abs x-r-2\gr)^2}{4t}}}{t^{1+\frac{N}{2}}}\qquad&\mbox {if }\abs x>r+ 2\gr+\sqrt{2t(N+2)}.
\EA\right.\EA
$$
By elementary analysis, if $x\in B_{r+3\gr}^c$,
$$(\abs x-r-2\gr)
e^{-\frac{(\abs x-r-2\gr)^{2}}{4t}} \leq e^{-\frac{(\abs x-r-3\gr)^{2}}{4t}} \left\{\BA {ll}\gr e^{-\frac{\gr^2}{4t}}\,\qquad&\mbox {if } 2t< \gr^2\\[2mm]
\myfrac{2t}{\gr}e^{-1+\frac{\gr^2}{4t}}\quad&\mbox {if } \gr^2\leq 2t\leq 2(r+2\gr)^2.
\EA\right.$$
However, since
$$\myfrac{\gr}{t}e^{-\frac{\gr^2}{4t}}\leq \myfrac{4}{\gr},
$$
we derive
$$(\abs x-r-2\gr)
e^{-\frac{(\abs x-r-2\gr)^2}{4t}}\leq \myfrac{Ct}{\gr}e^{-\frac{(\abs x-r-3\gr)^2}{4t}},
$$
and $(\ref{decay7})$ follows.\qeda 

\medskip

\noindent \Remark  In the subcritical case $1<q<q_c$,  it is easy to show by using \rlemma{decaylem6},  that any positive solution $u$ of $(\ref{mainE})$, such that $u(x,0)=0$ for $x\neq0$, satisfies
\begin{equation}\label{up1}
u(x,t)\leq Ct^{-\frac{1}{q-1}}\min\left\{1,\left(\myfrac{\abs {x}}{\sqrt t}\right)^{\frac{2}{q-1}-N}e^{-\frac{\abs {x}^2}{4t}}\right\}\qquad\forall (x,t)\in Q_\infty.
\end {equation}
This upper estimate corresponds to the one obtained in \cite{BPT}. If $F=\overline B_r$ the upper estimate is less esthetic. However, it is proved in \cite{MV2} by a barrier method that, if the initial trace of positive solution $u$ of $(\ref {mainE})$, vanishes outside F, and if $1<q<3$, there holds
\begin{equation}\label{up3}
u(x,t)\leq t^{-\frac{1}{q-1}}f_1((\abs x-r)/\sqrt t)\qquad\forall (x,t)\in Q_\infty,\;\abs x\geq r,
\end {equation}
where $f=f_1$ is the unique positive (and radial) solution of 
\begin{equation}
\label {bila''}\left\{\BA{l}
 f''+\myfrac{y}{2}f'+\myfrac{1}{q-1}f-f^q=0\quad\mbox {in }(0,\infty)\\[2mm]
f'(0)=0\,,\;\lim_{ y\to\infty}\abs y^{\frac{2}{q-1}}f(y)=0.
\EA\right.
\end {equation}
Notice that the existence of $f_1$ follows from \cite {BPT} since $q$ belongs to the subcritical range on exponents in dimension one. Furthermore $f_1$ has the following asymptotic expansion
$$f_1(y)=Cy^{(3-q)/(q-1)}e^{-y^2/4t}(1+\circ (1)))\quad\mbox {as }y\to\infty.
$$



\subsection {The upper Wiener test} 

\noindent \bdef {paradist} \rm {We define on $\mathbb R^{N}\ti\mathbb R$ the two 
{\it parabolic distances} $\gd_{2}$ and 
$\gd_{\ity}$ by
\begin {equation}
\label {d2}
\gd_{2}[(x,t),(y,s)]:=\sqrt {{\abs {x-y}}^{2}+{\abs {t-s}}},
\end {equation}
and
\begin {equation}
\label {dinfini}
\gd_{\ity}[(x,t),(y,s)]:=\max\{{\abs {x-y}},\sqrt {\abs {t-s}}\}.
\end {equation}
}\es
If $K\subset\BBR^{N}$ and $i=2,\infty$, 
$$
\gd_{i}[(x,t),K]=\inf \{\gd_{i}[(x,t),(y,0)]:y\in K\}=\left\{\BA 
{lc}\max\left\{\dist (x,K),\sqrt {\abs t}\right\}&\mbox { if 
}i=\infty,\\
[2mm]
\sqrt{\dist^{2} (x,K)+\abs t}&\mbox { if }i=2.\EA\right. $$
For $\gb>0$ and $i=2,\infty$, we denote by $\CB_{\gb}^i(m)$ the parabolic ball of 
center $m=(x,t)$ and radius $\gb$ in the parabolic distance 
$\gd_{i}$. \medskip

Let $K$ be {\it any} compact subset of
$\BBR^{N}$ and $\overline u_{K}$ the maximal solution of $(\ref {mequ})$ which blows up on 
$K$. The function $\overline u_{K}$ is constructed in \cite {MV2} as being the decreasing limit of the $\overline u_{K_{\ge}}$
($\ge>0$) when $\ge \to 
0$, where 
$$
K_{\ge}=\{x\in\BBR^{N}:\dist (x,K)\leq\ge\} $$
and $\overline u_{K_{\ge}}=\lim_{k\to\infty}u_{k,K_{\ge}}=\overline 
u_{K}$, where $u_{k}$ is the 
solution of the classical problem, 
\begin {equation}
\label {CD-k}\left\{\BA{rll}
\prt_{t} u_{k}-\Gd u_{k}+u_{k}^q=0\phantom{\chi_{_{K_{\ge}}}}\;&\quad\mbox {in }Q_{T},\\
[2mm] u_{k}=0\phantom{\chi_{_{K_{\ge}}}}\;&\quad\mbox {on }\prt_{\ell}Q_{T},\\
[2mm] u_{k}(.,0)=k\chi_{_{K_{\ge}}}\;&\quad\mbox {in }\BBR^{N}.
\EA\right.
\end {equation}
If $(x,t)=m \in \BBR^{N}\ti (0,T]$, we set $ d_{K}=\dist (x,K)$,
$D_{K}=\max\{\abs{x-y}:y\in K\}$ and $\gl=\sqrt { d^2_{K}+t}=\gd_{2}[m,K]$. We define a 
 slicing of $K$, by setting $d_n=d_n(K,t):=\sqrt {nt}$ ($n\in \BBN $), $d^{\pm}_{n}=\left(\sqrt {nt}\pm\myfrac{\sqrt t}{\sqrt {n}}\right)_+$ (the positive part is only needed when $n=0$) and
 $$
T^*_{n}=\overline B_{d^+_{n+1}}(x)\setminus 
B_{d^-_{n}}(x)\,,\;T_{n}=\overline B_{d_{n+1}}(x)\setminus 
B_{d_{n}}(x), \qquad\forall n\in\BBN,$$
thus $T^*_{0}=B_{2\sqrt{t}}(x)$, $T_{0}=B_{\sqrt{t}}(x)$, and 
$$
K_{n}(x,t)=K\cap T_{n}(x,t)\;\mbox { for }n\in\BBN 
\mbox { and }\CQ_{n}(x,t)=K\cap B_{d_{n+1}}(x,t). $$
When there is no ambiguity, we will skip the $(x,t)$ variable in the above sets. 
The main result of this section is the following discrete upper Wiener-type 
estimate.

\bth {upperW} Assume  $q\geq q_c$. Then there exists $C=C(N,q,T)>0$ such that
\begin {equation}
\label {uwe}
\overline u_{K}(x,t)\leq \myfrac {C}{t^{\frac{N}{2}}}\mysum{n=0}{a_{t}}
d_{n+1}^{N-\frac{2}{q-1}}e^{-\frac{n}{4}} C_{2/q,q'}\left(\myfrac {K_{n}}{d_{n+1}} \right)\qquad\forall (x,t)\in Q_T,
\end {equation}
where $a_{t}$ is the largest integer $j$ such that
$K_j\neq\emptyset$.
\es 
\medskip 

\noindent With no loss of generality, we can assume that $x=0$. Furthermore, in considering the scaling transformation $u_\ell(y,t)=\ell^{\frac{1}{q-1}}u(\sqrt\ell y,\ell t)$, with $\ell>0$, we can assume $t=1$. Thus the new compact singular set of the initial trace becomes $K/\sqrt\ell$, that we still denote $K$. We also set $a_{_K}=a_{_{K,1}}$
%
For $n\in\BBN_*$ set  $\gd_n=d_{n+1}-d_{n}$, then $\frac{1}{2\sqrt {n+1}}\leq \gd_n\leq \frac{1}{2\sqrt {n}}$. By convention $\gd_0=1$.
 It is possible to exhibit a collection $\Gth_{n}$ of points $a_{n,j}$ with center on
the sphere 
$\Gs_{n}=\{y\in\BBR^N:\abs {y}=(d_{n+1}+d_{n})/2\}$, such that 
$$
T_{n}\subset \bigcup_{a_{n,j} \in\Gth_{n} }B_{\gd_n}(a_{n,j}), \quad
\abs {a_{n,j}- a_{n,k}}\geq \gd_n\,\,\mbox { and 
}\,\,\#\Gth_{n}\leq Cn^{N-1}, $$
for some constant $C=C(N)$. If  $K_{n,j}=K_{n}\cap B_{\gd_n}(a_{n,j})$,
there holds $$
K=\bigcup_{ 0\leq n\leq a_{_K}}\bigcup_{a_{n,j} \in\Gth_{n} }K_{n,j}.$$ 

The first intermediate step is based on the {\it quasi-additivity} property of capacities developed in \cite {AB}.

\blemma {QA}Let $q\geq q_c $. There exists a constant $C=C(N,q)$ such that 
\begin {eqnarray}
\label {quasiadd0}
\sum_{a_{n,j} \in\Gth_{n} }R^{B_{2\gd_n}(a_{n,j})}_{2/q,q'}(K_{n,j})
\leq Cd_{n+1}^{N-\frac{2}{q-1}} C_{2/q,q'}\left(\myfrac{K_{n}}{d_{n+1}}\right)\qquad\forall n\in\BBN_* .
\end {eqnarray}
\es
\Proof 
 The following result is proved in
\cite [Th 3]{AB}:  if the spheres $B_{\gr_{j}^\gth}(b_{j})$, $\theta=1-2/N(q-1)$, are disjoint in $\BBR^N$ and
$G$ is an analytic subset of $\bigcup B_{\gr_{j}}(b_j)$ where the $\gr_j $ are positive
and smaller than some $\gr^*>0 $, there holds 
\begin {eqnarray}
\label {quasiadd0-t} 
C_{2/q,q'}(G)\leq \sum_jC_{2/q,q'}(G\cap B_{\gr_{j}}(b_j))\leq AC_{2/q,q'}(G),
\end {eqnarray}
for some $A$ depending on $N$, $q$ and $\gr^* $. This property is
called  {\it quasi-additivity}. We define for $n\in\BBN_*$, $$
\tilde T_{n}=d_{n+1}T_{n},\quad \tilde K_{n}=
d_{n+1} K_{n}\quad\mbox {and }\;\tilde \CQ_{n}=d_{n+1}\CQ_{n}.$$
Since $K_{n,j}\subset B_{\gd_n}(a_{n,j})$, it follows that 
$$\tilde K_{n,j}:=d_{n+1}K_{n,j}\subset B_{d_{n+1}\gd_n}(\tilde a_{n,j}).$$
Note that  by \rlemma{scal}
\begin {equation}
\label {quasiadd1} \BA {ll}
R^{B_{2\gd_n}(a_{n,j})}_{2/q,q'}(K_{n,j})= d_{n+1}^{\frac{2}{q-1}-N}R^{B_{2\gd_nd_{n+1}}(d_{n+1}a_{n,j})}_{2/q,q'}(\tilde K_{n,j})\\[2mm]\phantom{R^{B_{2\gd_n}(a_{n,j})}_{2/q,q'}(K_{n,j})}
\approx d_{n+1}^{\frac{2}{q-1}-N}C^{B_{2\gd_nd_{n+1}}(d_{n+1}a_{n,j})}_{2/q,q'}(\tilde K_{n,j})
\\[2mm]\phantom{R^{B_{2\gd_n}(a_{n,j})}_{2/q,q'}(K_{n,j})}
\approx d_{n+1}^{\frac{2}{q-1}-N}C_{2/q,q'}(\tilde K_{n,j})
\EA
\end {equation}
where $\tilde K_{n,j}=d_{n+1}K_{n,j}$. 
For a fixed $n>0$ and each repartition $\Gl$ of points $\tilde 
a_{n,j}= d_{n+1}\,a_{n,j}$ such that the balls $B_{2^\gth}(\tilde a_{n,j})$ are
disjoint, the quasi-additivity property 
holds: if we set $$
K_{n,\Gl}=\bigcup_{a_{n,j}\in\Gl}K_{n,j}\;,\quad
\tilde K_{n,\Gl}=d_{n+1}\,K_{n,\Gl}=\bigcup_{a_{n,j}\in\Gl}\tilde K_{n,j}
\quad\mbox {and }\;\tilde K_{n}=d_{n+1}\,K_{n},$$
then
\begin {eqnarray}
\label {quasiadd3}
\sum_{a_{n,j}\in\Gl}C_{2/q,q'}(\tilde K_{n,j}) \approx C_{2/q,q'}(\tilde 
K_{n,\Gl}).
\end {eqnarray}
The maximal cardinal of any such repartition $\Gl$ is of the order 
of $Cn^{N-1}$ for some positive constant $C=C(N)$, 
therefore, the number of repartitions needed for a full covering 
of the set $\tilde T_{n}$ is of finite order depending upon the 
dimension. 
Because $\tilde K_{n}$ is the union of the $\tilde K_{n,\Gl}$,
\begin {eqnarray}
\label {quasiadd4}
\sum_{a_{n,j}\in\Gth_n}C_{2/q,q'}(\tilde K_{n,j})=\sum_{\Gl}\sum_{a_{n,j}\in\Gl}C_{2/q,q'}(\tilde K_{n,j}) \approx C_{2/q,q'}(\tilde 
K_{n}).
\end {eqnarray}
By \rlemma{scal}, 
$$C_{2/q,q'}(\tilde K_{n})\leq C^{B_{2d_{n+1}}}_{2/q,q'}(\tilde K_{n})\approx d_{n+1}^{N-\frac{1}{q-1}}C^{B_{2}}_{2/q,q'}\left( \myfrac{K_{n}}{d_{n+1}}\right)\approx d_{n+1}^{N-\frac{1}{q-1}}C_{2/q,q'}\left( \myfrac{K_{n}}{d_{n+1}}\right),
$$
we obtain $(\ref {quasiadd0})$ by combining this last inequality with $(\ref {quasiadd1})$ and $(\ref{quasiadd4})$.\qeda
\medskip 

\noindent {\it Proof of \rth {upperW}.} {\it Step 1.} We first notice that
\begin {eqnarray}
\label {super0} \overline u_{K}\leq 
\sum_{0\leq n\leq a_{_K}}\sum_{a_{n,j}\in\Gth_{n}}\overline u_{K_{n,j}}.
\end {eqnarray}
Actually, since $K=\bigcup_{n}\bigcup_{a_{n,j}}K_{n,j}$, for any $0<\ge'<\ge$, 
there holds $\overline {K_{\ge'}}\subset \bigcup_{n}\bigcup_{a_{n,j}}K_{n,j\,\ge}$.
Because a finite sum of 
positive solutions of $(\ref {mequ})$ is a super solution,
\begin {eqnarray}
\label {super1'} \overline u_{K_{\ge'}}\leq 
\sum_{0\leq n\leq a_{_K}}\sum_{a_{n,j}\in\Gth_{n}}\overline u_{K_{n,j\,\ge}}.
\end {eqnarray}
Letting successively $\ge'$ and $\ge$ go to $0$ implies $(\ref {super0})$.\medskip


\noindent {\it Step 2. } Let $n\in\BBN $. Since $K_{n,j} \subset B_{\gd_{n}}(a_{n,j})$ and $\abs {x-a_{n,j}}=(d_ n +d_{n+1})/2$, we can apply the previous lemmas with $r=\gd_{n}$ and $\gr=r$. For $n\geq n_N$, there holds 
$t=1\geq (r+2\gr)^2=9/(n+1)$ and $\abs{x- a_{n,j}}=(\sqrt{n+1}-\sqrt n)/2\geq (2+C_N)(3/\sqrt{n+1})$ (notice that $n_N\geq 8$). Thus
\begin{equation}\label{new1}\BA {l}
u_{K_{n,j}}(0,1)\leq Ce^{\left(\sqrt{n}-3/\sqrt{n+1}\right)^2\!\!/4}\,R_{2/q,q'}^{B_{2\gd_n}(a_{n,j})}(K_{n,j})\leq Ce^{3/2}e^{-\frac{n}{4}}R_{2/q,q'}^{B_{2\gd_n}(a_{n,j})}(K_{n,j}).\EA
\end{equation}
Using \rlemma {QA} we obtain, with $d_n=d_n(1)=\sqrt{n+1\,}$
\begin{equation}\label{S1}
\sum_{n=n_{_N}}^{a_{_K}}\sum_{a_{n,j}\in\Gth_n}u_{K_{n,j}}(0,1)\leq 
C\sum_{n=n_{_N}}^{a_{_K}}d_{n+1}^{N-\frac{2}{q-1}}e^{-\frac{n}{4}}C_{2/q,q'}\left(\myfrac
{K_{n}}{d_{n+1}}\right). 
\end{equation}
Finally, we apply \rlemma{RPW} if $1\leq n<n_{_N}$ and get
\begin{equation}\label{S2}\BA {l}
\mysum{1}{n_{_N}-1}\mysum{a_{n,j}\in\Gth_n}{}u_{K_{n,j}}(0,1)\leq 
C\mysum{1}{n_{_N}-1}C_{2/q,q'}\left(\myfrac
{K_{n}}{d_{n+1}}\right)\\
\phantom{\mysum{1}{n_{_N}-1}\mysum{a_{n,j}\in\Gth_n}{}u_{K_{n,j}}(0,1)}\leq
C'\mysum{1}{n_{_N}-1}d_{n+1}^{N-\frac{2}{q-1}}e^{-\frac{n}{4}}C_{2/q,q'}\left(\myfrac
{K_{n}}{d_{n+1}}\right).
\EA\end{equation}
For $n=0$, we proceed similarly, in splitting $K_1$ in a finite number of $K_{1,i}$, depending only on the dimension, such that diam$\,K_{1,i}<1/3$. Combining $(\ref{S1})$ and $(\ref{S2})$, we derive
\begin{equation}\label{S3}
\overline u_K(0,1)\leq 
C\sum_{n=0}^{a_{_K}}d_{n+1}^{N-\frac{2}{q-1}}e^{-\frac{n}{4}}C_{2/q,q'}\left(\myfrac
{K_{n}}{d_{n+1}}\right). 
\end{equation}
In order to derive the same result for any $t>0$, we notice that 
$$\overline u_K(y,t)=t^{-\frac{1}{q-1}}\overline u_{K/\sqrt t}(y/\sqrt t,1).$$
Going back to the definition of $d_n=d_n(K,t)=\sqrt{nt}=d_n(K\sqrt t,1)$,
we derive from $(\ref{S3})$ and the fact that $a_{_{K,t}}=a_{_{K\sqrt t,1}}$
\begin{equation}\label{S4}
\overline u_K(0,t)\leq 
Ct^{-\frac{N}{2}}\mysum{n=0}{a{_K}}d_{n+1}^{N-\frac{2}{q-1}}e^{-\frac{n}{4}}C_{2/q,q'}\left(\myfrac
{K_{n}}{d_{n+1}}\right), 
\end{equation}
with $d_n=d_n(t)=\sqrt{t(n+1)\,}$. This is $(\ref{uwe})$ with $x=0$, and a space translation leads to the final result. \qeda \medskip 


\noindent {\it Proof of \rth {abovth}}. Let $m>0$ and $F_m=F\cap \overline B_m$. We denote by 
$U_{B_m^c}$ the maximal solution of $(\ref{mequ})$ in $Q_\infty$ the initial trace of which vanishes on $B_m$. Such a solution is actually the unique solution of $(\ref{mainE})$ which satisfies
$$\lim_{t\to 0}u(x,t)=\infty
$$
uniformly on $B_{m'}^c$, for any $m'>m$: this can be easily proved by noticing that 
$$U_{B_m^c\,\ell}(y,t)=\ell^{\frac{1}{q-1}}U_{B_m^c}(\sqrt\ell y,\ell t)=U_{B^c_{m/\sqrt\ell}}(y,t).
$$
Furthermore
$$\lim_{m\to\infty}U_{B^c_{m}}(y,t)=\lim_{m\to\infty}m^{-\frac{2}{q-1}}U_{B_1^c}(y/m, t/m^2)=0
$$
uniformly on any compact subset of $\overline Q_\infty$.
Since $\overline u_{F_m}+U_{B_m^c}$ is a super-solution, it is larger that $\overline u_{F}$ and therefore $\overline u_{F_m}\uparrow \overline u_{F}$. Because $W_{F_m}(x,t)\leq W_{F}(x,t)$ and 
$\overline u_{F_m}\leq C_1W_{F_m}(x,t)$, the result follows.\qeda\medskip

\noindent\Remark It is clear that \rth {abovth} still holds if $u$ is a positive subsolution of $(\ref{mequ})$ satisfying the initial trace condition $(\ref{ineq1-1})$.\medskip

\rth {abovth} admits the following integral expression.

\bth {upperWint} Assume  $q\geq q_c$. Then there exists a positive 
constant $C^*_{1}=C^*(N,q,T)$  such that, for any closed subset $F$ of $\BBR^N$, there holds
\begin {equation}
\BA{l}
\label {uwe2}
\overline u_{F}(x,t)\leq \myfrac {
C^*_{1}}{t^{1+\frac{N}{2}}}\myint{\sqrt t}{\sqrt {t(a_t+2)}}e^{-\frac{s^{2}}{4t}}s^{N-\frac{2}{q-1}}
C_{2/q,q'}\left(\myfrac {1}{s} F\cap B_{1}(x)\right)s\,ds,
\EA
\end {equation}
where $a_t=\min\{n:F\subset B_{\sqrt{n+1)t}}(x)\}$.
\es 
\Proof 
We first use
$$
C_{2/q,q'}\left(\myfrac {F_{n}}{d_{n+1}} \right)
\leq C_{2/q,q'}\left(\myfrac {F}{d_{n+1}} \cap B_1\right), $$
and we denote
\begin {equation}
\label {cap}
\Gf (s)=C_{2/q,q'}\left(\myfrac {F}{s} \cap B_1\right)\qquad\forall s>0.
\end {equation}

\noindent {\it Step 1. } The following inequality holds
\begin {equation}
\label {cap'1} c_1\Gf (\ga s)\leq\Gf (s)\leq c_2\Gf (\gb s)\qquad\forall s>0,\!\!\!\!\qquad\forall
1/2\leq
\ga\leq 1\leq \gb\leq 2,
\end {equation}
for some positive constants $c_1$, $c_2$ depending on $N$ and $q$. See \cite {AH} and \cite {MV6}. If $\gb\in [1,2]$, 
$$
\Gf (\gb s)= C_{2/q,q'}\left(\myfrac {1}{\gb} \left(\myfrac {F}{s}\cap B_\gb\right)\right)
\approx C_{2/q,q'} \left(\myfrac {F}{s}\cap B_\gb\right)\geq c_{1}\Gf (s). $$
 If $\ga\in [1/2,1]$, 
 $$
\Gf (\ga s)= C_{2/q,q'}\left(\myfrac {1}{\ga} \left(\myfrac {F}{s}\cap B_\ga\right)\right)
\approx C_{2/q,q'} \left(\myfrac {F}{s}\cap B_\ga\right)\leq c_{2}\Gf (s). $$

\noindent {\it Step 2. } By $(\ref{cap'1})$ 
$$
C_{2/q,q'}\left(\myfrac {F}{d_{n+1}} \cap B_1\right)\leq c_2C_{2/q,q'}\left(\myfrac {F}{s}
\cap B_{1}\right)\qquad\forall\;  
s\in [d_{n+1},d_{n+2}], $$
and $n\leq a_{_{t}}$. Then
\begin {eqnarray*}
c_2\myint{d_{n+1}}{d_ {n+2}}s^{N-\frac{2}{q-1}}e^{-s^2/4t}C_{2/q,q'}\left(\myfrac {F}{s} \cap
B_{1}\right)s\,ds&&
\\
[2mm]
\geq C_{2/q,q'}\left(\myfrac {F}{d_{n+1}} \cap B_1\right)&&\!\!\!\!\!\!\!\!\!\!\!\!
\myint{d_{n+1}}{d_{n+2}}s^{N-\frac{2}{q-1}}e^{-s^2/4t}s\,ds.
\end {eqnarray*}
Using the fact that $N-\frac{2}{q-1}\geq 0$, we get, 
\begin {eqnarray}
\label {uwe2a}
\myint{d_{n+1}}{d_{n+2}}s^{N-\frac{2}{q-1}}e^{-\frac{s^{2}}{4t}}s\,ds
\geq e^{-\frac{n+2}{4}}d_{n+1}^{N-\frac{2}{q-1}+1}(d_{n+2}-d_{n+1})\\
\geq \myfrac{t}{4e^2}d_{n+1}^{N-\frac{2}{q-1}}e^{-\frac{n}{4}}.\phantom{-------} 
\end {eqnarray}
Thus
\begin {eqnarray}
\label {uwe2b}
\overline u_F(x,t)\leq \myfrac {C}{t^{1+\frac{N}{2}}}
\int_{\sqrt t}^{\sqrt {t(a_t+2)}} s^{N-\frac{2}{q-1}}e^{-\frac{s^{2}}{4t}}C_{2/q,q'}\left(\myfrac {1}{s}
F\cap B_{1}\right)s\,ds,
\end {eqnarray}
which ends the proof.
\qeda \medskip 

\mysection {Estimate from below} If $\gm\in\GTM_{_{+}}^q(\BBR^{N})\cap \GTM^b(\BBR^{N})$,
we denote by
$u_{\gm}=u_{\gm,0}$ 
the solution of 
\begin {equation}
\label {sub1}\left\{\BA{rll}
\prt_{t}u_{\gm}-\Gd u_{\gm}+u_{\gm}^q=0\;&\quad\mbox {in }Q_{T},\\
[2mm] u_{\gm}(.,0)=\gm& \quad\mbox {in }\BBR^{N}.
\EA\right.
\end {equation}
The maximal $\gs$-moderate solution of $(\ref{mequ})$ which has an initial trace vanishing outside a closed set $F$ is defined by
\begin {equation}
\label {lwe}
\underline u_F=\sup\left\{u_\gm:\gm\in\GTM_{_{+}}^q(\BBR^{N})\cap \GTM^b(\BBR^{N})\,,\;
\gm(F^c)=0\right\}.
\end {equation}
The main result of this section is the next one
\bth {lowerW} Assume  $q\geq q_c$. There exists a constant $C_2=C_{2}(N,q,T)>0$ such that, for any closed subset $F\subset\BBR^N$, there holds
\begin {equation}
\label {lwe-t}
\underline u_{F}(x,t)\geq C_{2}W_F(x,t)\qquad\forall (x,t)\in Q_T.
\end {equation}
\es\smallskip

We first assume that $F$ is compact, and we denote it by $K$. The first observation is that if $\gm\in\mathfrak M_+^q(\BBR^N)$,  $u_{\gm}\in L^q(Q_{T})$ (see lemma below) and $0\leq u_{\gm}\leq 
\BBH[\gm]:=\BBH_{\gm}$. Therefore
\begin {equation}
\label {sub2} u_{\gm}\geq \BBH[\gm]-\BBG\left[\BBH[\gm]^q\right],
\end {equation}
where $\BBG$ is the parabolic Green potential in $Q_{T}$ defined by $$
\BBG[f](t)=\int_{0}^t\BBH [f(s)](t-s)ds=\int_{0}^t\int_{\BBR^{N}}H(.,y,t-s)f(y,s)dyds.
$$

The main idea of the proof is as follows. For any $(x,t)\in Q_T$, construct a measure $\gm=\gm(x,t)\in\mathfrak M_+^q(\BBR^N)$ such that there holds
\begin {equation}
\label {sub2*}
\BBH[\gm](x,t)\geq CW_K(x,t)\qquad\forall (x,t)\in Q_T,
\end {equation}
and
\begin {equation}
\label {sub2'} \BBG\left(\BBH[\gm]\right)^q \leq  C\,\BBH[\gm]\quad\mbox {in }Q_T,
\end {equation}
with constants $C$ depends only on $N$, $q$, and $T$. Then replace $\gm$ by $\gm_\ge=\ge\gm$ with $\ge=(2C)^{-\frac{1}{q-1}}$ in order to derive
\begin {equation}
\label {sub2''} u_{\gm_\ge}\geq 2^{-1}\BBH_{\gm_\ge}\geq 2^{-1}CW_K.
\end {equation}
From this follows 
\begin {equation}
\label {sub2'''} \underline u_{K}\geq 2^{-1}\BBH_{\gm_\ge}\geq 2^{-1}CW_K.
\end {equation}
and the proof of \rth{lowerW} with $C_2=2^{-1}C$. In the following sections we describe the construction of measures $\mu(x,t)$ satisfying $(\ref{sub2*})$ and $(\ref{sub2'})$. \medskip

%
%
\subsection {Estimate from below of the solution of the heat equation}
\noindent The purely spatial slicing used is the trace on $\BBR^N\ti\{0\}$ of an {\it extended slicing} in $Q_T$ which is constructed as follows: if $K$ is a 
compact subset of $\BBR^{N}$, $m=(x,t)$, we define $d_{K}$, $\gl$, $d_{n}$ and $a_{t}$ as
in 
Section 2.3. Let $\ga\in (0,1)$ to be fixed 
later on, we define $\CT_{n}$ for $n\in\BBZ$ 
by $$
\CT_{n}=\left\{\BA {lc}\CB^2_{\sqrt {t(n+1)}}(m)\setminus 
\CB^2_{\sqrt {tn}}(m)\;&\mbox { if }n\geq 1,\\
[4mm]
\CB^2_{\ga^{-n}\sqrt t}(m)\setminus 
\CB^2_{\ga^{1-n}\sqrt t}(m)\;&\mbox { if }n\leq 0,
\EA\right.$$
and put $$
\CT^{*}_{n}=\CT_{n}\cap\{s:0\leq s\leq t\}, \;\mbox { 
for }n\in\BBZ. $$
We recall that for $n\in\BBN_{*}$, $$
\CQ_{n}=K\cap \CB^2_{\sqrt{t(n+1)} }(m) =K\cap B_{d_{n}}(x)$$
and $$
K_{n}=K\cap \CT_{n+1}=K\cap \left(B_{d_{n+1}}(x)\setminus 
B_{d_{n}}(x)\right).$$
Let $\gn_{n}\in \GTM^b_{_{+}}(\BBR^{N})\cap 
W^{-2/q,q}(\BBR^{N})$ be the {\it $q$-capacitary measure} of the set $K_{n}/d_{n+1}$. 
See \cite [Sec. 2.2]{AH}. Such a measure has support in $K_{n}/d_{n+1}$ 
and
\begin {equation}
\label {capmes}
\gn_{n}(K_{n}/d_{n+1})=C_{2/q,q'}(K_{n}/d_{n+1})\mbox { and 
}\norm{\gn_{n}}_{W^{-2/q,q'}(\BBR^{N})}=\left(C_{2/q,q'}(K_{n}/d_{n+1})\right)^{1/q}.
\end {equation}
We define $\gm_{n}$ as follows
\begin {equation}
\label {capmes0}
\gm_{n}(A)=d_{n+1}^{N-\frac{2}{q-1}}\gn_{n} (A/d_{n+1})\qquad\forall A\subset 
K_{n},\;A\;\mbox { Borel },
\end {equation}
and set $$
\gm_{t,K}=\sum_{n=0}^{a_{t}}\gm_{n}, $$
and
\begin {equation}
\label {capmes1}
\BBH_{\gm_{t,K}}=\sum_{n=0}^{a_{t}}\BBH_{\gm_{n}}.
\end {equation}
\bprop {heatsub} Let $q\geq q_{c}$, then there holds
 \begin {equation}
\label {capmes2}
 \BA {l}
 \BBH_{\gm_{t,K}}(x,t)
 \geq \myfrac{1}{(4\gp t)^{\frac{N}{2}}}\mysum{n=0}{a_{t}}\;
 e^{-\frac{n+1}{4}}d^{N-\frac{2}{q-1}}_{n+1}C_{2/q,q'}\left(\myfrac 
 {K_{n}}{d_{n+1}}\right),
 \EA
\end {equation}
 in $\BBR^{N}\ti (0,T)$.
\es 
\Proof Since
\begin {equation}
\label {capmes3}
\BBH_{\gm_{n}}(x,t)= \myfrac{1}{(4\gp t)^{\frac{N}{2}}}\int_{K_{n}} e^{-\frac{\abs
{x-y}^{2}}{4t}}d\gm_{n},
\end {equation}
and
 $$
y\in K_{n}\Longrightarrow \abs {x-y}\leq d_{n+1},
 $$
 $(\ref {capmes2})$ follows because of $(\ref {capmes0})$ and $(\ref 
 {capmes1})$.\qeda \medskip
\subsection {Estimate from above of the nonlinear term} 
We write $(\ref {sub2})$ under the form
 \begin {equation}
\label {nln1}\BA {l} u_{\gm}(x,t)\geq
\mysum{n\in\BBZ}{}\BBH_{\gm_{n}}(x,t)-\myint{0}t\myint{\BBR^{N}}{}H(x,y,t-s)\left%
[\mysum{n\in A_{K}}{}\BBH_{\gm_{n}}(y,s)\right]^qdyds\\
[4mm]
\phantom {u_{\gm}(x,t)\geq} = I_{1}- I_{2}.\EA
 \end {equation}
since $\gm_n=0$ if $n\notin A_K=\BBN\cap[1,a_t]$, and
  \begin {equation}
\label {nln2}\BA {l} I_{2}= 
\myfrac {1}{(4\gp)^{\frac{N}{2}}}\myint{0}t\myint{\BBR^{N}}{}(t-s)^{-\frac{N}{2}}e^{-\frac{\abs{x-y}^{2}}{4(t-s)}}\left[\mysum{n\in A_{K}} 
{}\BBH_{\gm_{n}}(y,s)\right]^qdyds\\
[4mm]
\phantom {I_{2}}= \myfrac {1}{(4\gp)^{\frac{N}{2}}}(J_{\ell}+J_{\ell}'),
\EA
 \end {equation}
 for some $\ell\in\BBN^{*}$ to be fixed later on, where 
$$
\BA {l} J_{\ell}\!=\!\mysum{p\in\BBZ} 
{}\dint_{\CT^{*}_{p}}\!\!(t-s)^{-\frac{N}{2}}e^{-\frac{\abs{x-y}^{2}}{4(t-s)}}\left[\mysum{n<p+\ell} 
{}\BBH_{\gm_{n}}(y,s)\right]^qdyds,
\EA $$
and $$
\BA {l} J_{\ell}'=\!\!\mysum{p\in\BBZ} 
{}\dint_{\CT^{*}_{p}}\!\!(t-s)^{-\frac{N}{2}}e^{-\frac{\abs{x-y}^{2}}{4(t-s)}}\left[\mysum{n\geq p+\ell} 
{}\BBH_{\gm_{n}}(y,s)\right]^qdyds.
\EA $$
\medskip

The next estimate will be used several times in the sequel. 

\blemma {kernest}
 Let $0<a<b$ and $t>0$, then, 
\begin {eqnarray*}
\max\left\{\gs^{-\frac{N}{2}}e^{ -\frac{\gr^2}{4\gs}}:0\leq\gs\leq t,\; at\leq\gr^{2}+\gs\leq 
bt\right\}=e^{\frac{1}{4}}\left\{\BA{ll}t^{-\frac{N}{2}}e^{-\frac{a}{4}}\;\, &\mbox { if 
}\myfrac {a}{2N}>1,\\
[4mm]
\left(\myfrac {2N}{at}\right)^{\frac{N}{2}}\!\!\!e^{-\frac{N}{2}}\;\, &\mbox { if 
} \myfrac {a}{2N}\leq 1.
\EA\right.
\end  {eqnarray*}
 \es 
 \Proof Set
 $$
\CJ(\gr,\gs)=\gs^{-\frac{N}{2}}e^{-\frac{\gr^2}{4\gs}}
 $$
 and
 $$
\CK_{a,b,t}=\left\{(\gr,\gs)\in [0,\infty)\ti (0,t]: \; at\leq\gr^{2}+\gs\leq 
bt\right\}.
 $$
We first notice that, for fixed $\gs$, the maximum of $\CJ(.,\gs)$ 
is achieved for $\gr$ minimal. If $\gs\in [at,bt]$ the 
 minimal value of $\gr$ is $0$, while if $\gs\in (0,at)$, the minimum 
 of $\gr$ is $\sqrt {at-s}$. \smallskip
 
\nind - Assume first $ a\geq 1$, 
 then 
 $\CJ(\sqrt {at-\gs},\gs)=e^{\frac{1}{4}}\gs^{-\frac{N}{4}}e^{-\frac{at}{4\gs}}$. 
 Thus if $1\leq a/2N$, the minimal value of $\CJ(\sqrt {at-\gs},\gs)$
 is $e^{\frac{1-2N}{4}}\left(\frac{2N}{at}\right)^{\frac{N}{2}}$, while if $a/2N<1\leq 
 a$, the minimum is $e^{\frac{1}{4}}t^{-\frac{N}{2}}e^{-\frac{a}{4}}$.\smallskip
 
 \nind - Assume now $ a\leq 1$. Then
\begin {eqnarray*}
\max\{\CJ(\gr,\gs):\;(\gr,\gs)\in \CK_{a,b,t}\}
 =\max \left\{\max_{\gs\in (at,t]}\CJ(0,\gs),\,
 \max_{\gs\in (0,at]}\CJ(\sqrt {at-\gs},\gs)\right\}\\
 =\max\left\{(at)^{-\frac{N}{2}},\,e^{\frac{1-2N}{4}}\left(\frac{2N}{at}\right)^{\frac{N}{2}}\right\}
 \qquad\quad\;\,\\
 =e^{\frac{1-2N}{4}}\left(\frac{2N}{at}\right)^{\frac{N}{2}}.\phantom {------------}
\end {eqnarray*}
 Combining these two estimates, we derive the result.\qeda \medskip
 
 \nind \Remark The following variant of \rlemma {kernest} will be 
 useful in the sequel: 
 {\it For any $\gth\geq 1/2N$ there holds}
 \begin {eqnarray}
\label {var}
\max\{\CJ(\gr,\gs):\;(\gr,\gs)\in\CK(a,b,t)\} \leq 
e^{\frac{1}{4}}\left(\myfrac {2N\gth}{t}\right)^{\frac{N}{2}}e^{-\frac{a}{4}}\quad \mbox {if } \gth a\geq 1.
\end {eqnarray}
\blemma {LJ1}There exists a positive constant $C=C(N,\ell,q)$ such that 
\begin {eqnarray}
\label {J1-2}
 J_{\ell}\leq Ct^{-\frac{N}{2}}\mysum{n=1}{a_{t}}d^{N-\frac{2}{q-1}}_{n+1}
 e^{-(1+(n-\ell)_{_{+}})/4}\;C_{2/q,q'}\left(\myfrac {K_{n}}{d_{n+1}}\right).
\end {eqnarray}
\es
\Proof  The set of the $p$'s for the summation in $J_{\ell}$ is reduced to 
$\BBZ\cap [-\ell+2,\infty)$, thus we write 
$$
J_{\ell}=J_{1,\ell}+J_{2,\ell} $$
where $$
J_{1,\ell}=\mysum{p=2-\ell}{0}\dint_{\CT^{*}_{p}}\!\!(t-s)^{-\frac{N}{2}}e^{-\frac{\abs {x-y}^{2}}{4(t-s)}}\left[\mysum{n<p+\ell} 
{}\BBH_{\gm_{n}}(y,s)\right]^q $$
and $$
J_{2,\ell}=\mysum{p=1}{\infty}\dint_{\CT^{*}_{p}}\!\!(t-s)^{-\frac{N}{2}}e^{-\frac{\abs {x-y}^{2}}{4(t-s)}}\left[\mysum{n<p+\ell} 
{}\BBH_{\gm_{n}}(y,s)\right]^q. $$
If $p=2-\ell,\ldots,0,$ $$
(y,s)\in\CT_{p}^{*}\Longrightarrow t\ga^{2-2p}\leq 
\abs {x-y}^{2}+t-s\leq t\ga^{-2p}, 
$$
and, if $p\geq 1$ $$
(y,s)\in\CT_{p}^{*}\Longrightarrow pt\leq 
\abs {x-y}^{2}+t-s\leq (p+1)t. 
$$
By \rlemma {kernest} and $(\ref {var})$, there exists $C=C(N,\ell,\ga)>0$ such that 
\begin {eqnarray}
\label {p-neg}
\max\left\{(t-s)^{-\frac{N}{2}}e^{-\frac{\abs {x-y}^{2}}{4(t-s)}}:(y,s)\in\CT^{*}_{p}\right\}\leq C 
t^{-\frac{N}{2}}e^{-\ga^{2-2p}/4},
\end {eqnarray}
if $p=2-\ell,\ldots,0$, and
\begin {eqnarray}
\label {p-pos}
\max\left\{(t-s)^{-\frac{N}{2}}e^{-\frac{\abs {x-y}^{2}}{4(t-s)}}:(y,s)\in\CT^{*}_{p}\right\}\leq C 
t^{-\frac{N}{2}}e^{-p/4},
\end {eqnarray}
if $p\geq 1$. When $p=2-\ell,\ldots,0$
\begin {eqnarray}
\label 
{sum-H}\left[\mysum{1}{p+\ell-1}\BBH_{\gm_{n}}(y,s)\right]^q
\leq C
\mysum{1}{p+\ell-1}\BBH^q_{\gm_{n}}(y,s),
\end {eqnarray}
for some $C=C(\ell,q)>0$, thus
\begin {eqnarray}
\label {p-neg2} J_{1,\ell}\leq Ct^{-\frac{N}{2}}\mysum{p=2-\ell}{0}e^{-\frac{\ga^{2-2p}}{4}}
\mysum{n=1}{p+\ell-1}\norm 
{\BBH_{\gm_{n}}}^q_{L^q(Q_{t})}\;\notag\\
\leq Ct^{-\frac{N}{2}}\mysum{n=1}{\ell-1}\norm{\BBH_{\gm_{n}}}^q_{L^q(Q_{t})}
\mysum{p=n-\ell+1}{0}e^{-\frac{\ga^{2-2p}}{4}}\\
\leq Ct^{-\frac{N}{2}}e^{-\frac{\ga^{2\ell-2}}{4}}
\mysum{n=1}{\ell-1}\norm{\BBH_{\gm_{n}}}^q_{L^q(Q_{t})}.
\;\;\;\;\quad\quad\notag
\end {eqnarray}
 If the set of $p$'s is not upper bounded, we introduce some parameter $\gd>0$ to be made 
 precise later on. Then
 \begin {eqnarray}
\label 
{sum-H'}\left[\mysum{1}{p+\ell-1}\BBH_{\gm_{n}}(y,s)\right]^q
\leq \left[\mysum{1}{p+\ell-1}e^{\gd q'\frac{n}{4}}\right]^{q/q'}
\mysum{1}{p+\ell-1}e^{-\frac{\gd qn}{4}}\BBH^q_{\gm_{n}}(y,s),
\end {eqnarray} 
with $q'=q/(q-1)$. If, by convention $\gm_{n}=0$ whenever $n>a_{t}$, we obtain, for some
$C>0$
 which depends also on $\gd$, 
\begin {eqnarray}
\label {p-pos2} J_{2,\ell}\leq Ct^{-\frac{N}{2}}\mysum{p=1}{\infty}e^{\frac{\gd (p+\ell-1)q-p}{4}}
\mysum{n=1}{p+\ell -1} e^{-\frac{\gd qn}{4}}\norm 
{\BBH_{\gm_{n}}}^q_{L^q(Q_{t})}\;\;\;\;\qquad\notag\\
\leq Ct^{-\frac{N}{2}}\mysum{n=1}{\infty}\norm{\BBH_{\gm_{n}}}^q_{L^q(Q_{t})} e^{-\frac{\gd qn}{4}}
\mysum{p=(n-\ell+1)\vee 1}{\infty}e^{\frac{\gd (p+\ell-1)q-p}{4}}\\
\leq Ct^{-\frac{N}{2}}
\mysum{n=1}{\infty}e^{-\frac{1+(n-\ell)_{+}}{4}}\norm{\BBH_{\gm_{n}}}^q_{L^q(Q_{t})}.
\qquad\qquad\qquad\;\;\;\;\;\;\;\;\quad\notag
\end {eqnarray}
Notice that we choose $\gd$ such that $\gd\ell q<1$. Combining $(\ref 
{p-neg2})$ and $(\ref {p-pos2})$, we derive $(\ref {J1-2})$ from \rlemma {meas}, 
$(\ref {capmes})$ and $(\ref {capmes0})$.
\qeda\\

    

\medskip

The set of indices $p$ for which the 
$\gm_{n}$ terms  are not zero in $J'_{\ell}$ is $\BBZ\cap 
(-\infty,a_{t}-\ell ]$. We write 
$$
J'_{\ell}=J'_{1,\ell}+J'_{2,\ell},$$ 
where $$
J'_{1,\ell}=\mysum{p=-\ity}{0}
\dint_{\CT^{*}_{p}}\!\!(t-s)^{-\frac{N}{2}}e^{-\frac{\abs {x-y}^{2}}{4(t-s)}}\left[\mysum{n=1\vee p+\ell} {\infty}
\BBH_{\gm_{n}}(y,s)\right]^q dyds, $$
and 
$$
J'_{2,\ell}=\mysum{p=1}{a_{t}-\ell}
\dint_{\CT^{*}_{p}}(t-s)^{-\frac{N}{2}}e^{-\frac{\abs {x-y}^{2}}{4(t-s)}}\left[\mysum{n=p+\ell} {\infty}
\BBH_{\gm_{n}}(y,s)\right]^q \!\!\!\!dyds. $$
                                                        %

\blemma {LJ3-1} There exists a constant $C=C
(N,q,\ell)>0$ such that
\begin {eqnarray}
\label {J3-1} J'_{1,\ell}\leq C t^{1-\frac{Nq}{2}}
\mysum{n=0} 
{a_{t}}e^{-\frac{(1+\gb_{0})(n-h)_{+}}{4}}d_{n+1}^{Nq-2q'}C^q_{2/q,q'}
\left(\myfrac{K_{n}}{d_{n+1}}\right),
\end {eqnarray}
where $\gb_{0}=(q-1)/4$ and $h=2q(q+1)/(q-1)^{2}$.
\es 
\Proof Since
\begin {eqnarray}
\label {J3-2} (y,s)\in\CT_{p}^{*},\mbox { and } (z,0)\in K_{n}\Longrightarrow
\abs {y-z}\geq (\sqrt n-\ga^{-p})\sqrt t,
\end {eqnarray}
there holds $$
\BBH_{\gm_{n}}(y,s)\leq (4\gp s)^{-\frac{N}{2}} e^{-\frac{(\sqrt n-\ga^{-p})^{2}t}{4s}}\gm_{n}(K_{n})\leq
C t^{-\frac{N}{2}}e^{-\frac{(\sqrt n-\ga^{-p})^{2}}{4}}\gm_{n}(K_{n}), $$
by \rlemma {kernest}. Let $\{\ge_{n}\}$ be a sequence of positive numbers such that 
$$
A_{\ge}=\mysum{n=1}{\infty}\ge^{q'}_{n}<\infty, $$
then 
\begin {equation}
\label {J3-3}\BA {l} J'_{1,\ell}\leq CA_{\ge}^{q/q'}t^{-\frac{Nq}{2}}\mysum{p=-\ity}{0}
\dint_{\CT^{*}_{p}}(t-s)^{-\frac{N}{2}}e^{-\frac{\abs {x-y}^{2}}{4(t-s)}}\!\!\!\!\!\!\mysum{n=1\vee
(p+\ell)} {\infty}
\ge^{-q}_{n}e^{-q\frac{(\sqrt n-\ga^{-p})^{2}}{4}}\gm^q_{n}(K_{n}) ds\,dy\\
[4mm]
\phantom{J'_{1,\ell}}\leq CA_{\ge}^{q/q'}t^{-\frac{Nq}{2}}\mysum{n=1}{\infty}
\ge^{-q}_{n}\gm^q_{n}(K_{n})\mysum{-\infty} {p=0\wedge( n-\ell)} e^{-\frac{q(\sqrt n-\ga^{-p})^{2}}{4}}
\dint_{\CT^{*}_{p}}(t-s)^{-\frac{N}{2}}e^{-\frac{\abs {x-y}^{2}}{4(t-s)}}ds\,dy\\
\phantom{J'_{1,\ell}}\leq CA_{\ge}^{q/q'}t^{-\frac{Nq}{2}}
\mysum{n=1}{\infty}
\ge^{-q}_{n}\gm^q_{n}(K_{n})e^{-\frac{q(\sqrt n-1)^{2}}{4}}
\dint_{\{\cup_{p\leq 0}\CT^{*}_{p}\}}(t-s)^{-\frac{N}{2}}e^{-\frac{\abs {x-y}^{2}}{4(t-s)}}ds\,dy\\
\phantom{J'_{1,\ell}}\leq CA_{\ge}^{q/q'}t^{1-\frac{Nq}{2}}
\mysum{n=1}{\infty}
\ge^{-q}_{n}\gm^q_{n}(K_{n})e^{-\frac{q(\sqrt n-1)^{2}}{4}}.
\EA
\end {equation}
Set $h=2q(q+1)/(q-1)^{2}$ and $Q=(1+q)/2$, then 
$q(\sqrt n-1)^{2}\geq Q(n-h)_{+}$ for any $n\geq 1$. If we choose
$\ge_{n}=e^{-\frac{(q-1)(n-h)_{+}}{16q}}$, there holds $\ge^{-q}_{n}e^{-\frac{q(\sqrt n-1)^{2}}{4}}\leq 
e^{-\frac{(q+3)(n-h)_{+}}{16}}$. Finally $$
J'_{1,\ell}\leq Ct^{1-\frac{Nq}{2}}
\mysum{n=1}{\infty} e^{-\frac{(1+\gb_{0})(n-h)_{+}}{4}}\gm^q_{n}(K_{n}), $$
with $\gb_{0}=(q-1)/4$, which yields to $(\ref {J3-1})$ by the choice of 
the $\gm_{n}$.\qeda\\

 In order to make easier the obtention of the estimate of the term $ 
 J'_{2,\ell}$, we first give the proof in dimension $1$.
\blemma {LJ3-2} Assume $N=1$ and $\ell$ is an integer larger than $1$. 
There exists a positive constant $C=C(q,\ell)>0$ such that 
\begin {eqnarray}
\label {J3-8} J'_{2,\ell}\leq Ct^{-1/2}
\mysum{n=\ell}{a_t}e^{-\frac{n}{4}}d^{\frac{q-3}{q-1}}_{n+1}
C_{2/q,q'}\left(\myfrac{K_{n}}{d_{n+1}}\right).
\end {eqnarray}
\es 
\Proof If $(y,s)\in \CT^{*}_{p}$ and $z\in K_{n}$ ($p\geq 1$, $n\geq 
p=\ell$) , there holds 
$\abs {x-y}\geq \sqrt t\sqrt p $ and $\abs{y-z}\geq \sqrt t(\sqrt n-\sqrt {p+1})$. 
Therefore $$
J'_{2,\ell}\leq C\sqrt t
\mysum{p=1}{a_{t}-\ell}\myfrac {1}{\sqrt p}\int_{0}^te^{-\frac{pt}{4(t-s)}}
\left(\mysum{n=p+\ell} {a_{t}}s^{-1/2}e^{-\frac{(\sqrt n-\sqrt {p+1}\,)^{2}t}{4s}}\gm_{n}(K_{n})\right)^q. $$
If $\ge\in (0,q)$ is some positive 
parameter which will be made more precise later on, there holds 
$$
\BA{l}\left(\mysum{n=p+\ell} {a_{t}}s^{-1/2}e^{-\frac{(\sqrt n-\sqrt {p+1}\,)^{2}t}{4s}}\gm_{n}(K_{n})\right)^q\\
\phantom {--------}\leq 
\left(\mysum{n=p+\ell} {a_{t}}e^{-\ge q'\frac{(\sqrt n-\sqrt {p+1}\,)^{2}t}{4s}}\right)^{q/q'}
\mysum{n=p+\ell} {a_{t}}s^{-\frac{q}{2}} e^{-(q-\ge)\frac{(\sqrt n-\sqrt {p+1}\,)^{2}t}{4s}}\gm^q_{n}(K_{n}),
\EA$$
by H\"older's inequality. By comparison between series and integrals and using Gauss integral
$$\BA {l}\mysum{n=p+\ell} {a_{t}}e^{-\ge q'\frac{(\sqrt n-\sqrt {p+1})^{2}t}{4s}}
\leq \myint{p+\ell}{\infty}e^{-\ge q'\frac{(\sqrt x-\sqrt {p+1})^{2}t}{4s}}dx\\
\phantom{\mysum{n=p+\ell} {a_{t}}e^{-\ge q'\frac{(\sqrt n-\sqrt {p+1})^{2}t}{4s}}}
=2\myint{\sqrt{p+\ell}-\sqrt{p+1}}{\infty}e^{-\frac{\ge q'x^{2}t}{4s}}(x+\sqrt{p+1})dx\\
\phantom{\mysum{n=p+\ell} {a_{t}}e^{-\ge q'\frac{(\sqrt n-\sqrt {p+1})^{2}t}{4s}}}
\leq \myfrac{4s}{\ge q't}e^{-\ge q'\frac{(\sqrt {p+\ell}-\sqrt 
{p+1})^{2}t}{4s}}+2\sqrt{p+1}\myint{\sqrt{p+\ell}-\sqrt{p+1}}{\infty}e^{-\frac{\ge q'x^{2}t}{4s}}dx
\\
\phantom{\mysum{n=p+\ell} {a_{t}}e^{-\ge q'\frac{(\sqrt n-\sqrt {p+1})^{2}t}{4s}}}
\leq C\sqrt\myfrac{(p+1)s}{t}e^{-\ge q'\frac{(\sqrt {p+\ell}-\sqrt {p+1})^{2}t}{2s}}\\
\phantom{\mysum{n=p+\ell} {a_{t}}e^{-\ge q'\frac{(\sqrt n-\sqrt {p+1})^{2}t}{4s}}}
\leq C\sqrt\myfrac{(p+1)s}{t}.
\EA$$
If we set $q_{\ge}=q-\ge$, then 
$$
J'_{2,\ell}\leq C\ge^{-q'/q}t^{1-\frac{q}{2}}\mysum{n=\ell+1}{\infty}\gm^q_{n}(K_{n})
\mysum{p=1}{n-\ell}p^{\frac{q-2}{2}}\myint{0}{t}(t-s)^{-1/2}s^{-1/2}
e^{-\frac{pt}{4(t-s)}}e^{-q_{\ge}\frac{(\sqrt n-\sqrt {p+1}\,)^{2}t}{4s}}ds. $$
where $C=C(\ge,q)>0$. Since 
$$
\BA{l}\myint{0}{t}(t-s)^{-1/2}s^{-1/2}
e^{-\frac{pt}{4(t-s)}}e^{-q_{\ge}\frac{(\sqrt n-\sqrt {p+1\,})^{2}t}{4s}}ds\\
\phantom {---------} =\myint{0}{1}(1-s)^{-1/2}s^{-1/2} e^{-\frac{p}{4(1-s)}}e^{-q_{\ge}\frac{(\sqrt n-\sqrt {p+1\,})^{2}}{4s}}ds,
\EA$$
we can apply \rlemma{integral} with $a=1/2$, $b=1/2$, $A=\sqrt p$ and 
$B=\sqrt {q_{\ge}}(\sqrt n-\sqrt {p+1})$. In this range of indices 
$B\geq\sqrt {q_{\ge}}(\sqrt {p+\ell}-\sqrt {p+1})\geq 
\frac{\sqrt {q_{\ge}}(\ell-1)}{\sqrt p}$, thus $\gk=\sqrt {q_{\ge}}(\ell-1)$ and
$$\sqrt \myfrac {A}{A+B}\sqrt \myfrac {B}{A+B}
\leq p^{\frac{1}{4}}n^{-1/2}(\sqrt n-\sqrt p)^{1/2}.
$$
Therefore
\begin {equation}
\label {integralest2}
\myint{0}{t}(t-s)^{-1/2}s^{-\frac{q}{2}} e^{-\frac{pt}{4(t-s)}}e^{-q\frac{(\sqrt n-\sqrt
{p+1})^{2}t}{4s}}ds
\leq\myfrac{Cp^{\frac{1}{4}}(\sqrt n-\sqrt p)^{1/2}}{\sqrt n}e^{-\frac{(\sqrt p+\sqrt{q_{\ge}}(\sqrt n-\sqrt {p+1}))^{2}}{4}},
\end {equation}
which implies
\begin {equation}
\label {integralest3} J'_{2,\ell}\leq 
Ct^{1-\frac{q}{2}}\mysum{n=\ell+1}{a_{t}}\myfrac{\gm_{n}^q(K_{n})}{\sqrt n}\mysum
{p=1}{n-\ell}p^{\frac{2q-3}{4}}(\sqrt n-\sqrt p)^{1/2}e^{-\frac{(\sqrt p+\sqrt{q_{\ge}}(\sqrt n-\sqrt {p+1}))^{2}}{4}},
\end {equation}
where $C$ depends of $\ge$, $q$ and $\ell$. By \rlemma {A2} 
\begin {equation}
\label {integralest4} J'_{2,\ell}\leq 
Ct^{1-\frac{q}{2}}\mysum{n=\ell+1}{a_{t}}n^{\frac{q-3}{2}}e^{-\frac{n}{4}}\gm_{n}^q(K_{n})
\end {equation}
 Because 
$\gm_{n}(K_{n})=d^{\frac{q-3}{q-1}}_{n+1}C_{2/q,q'}\left(\myfrac{K_{n}}{d_{n+1}}\right)$ 
(remember $N=1$) and diam$\,\frac{K_{n}}{d_{n+1}}\leq n^{-1}$, there holds 
\begin {equation}
\label {integraltest8}\mu^q_n(K_n)\leq C\left(\frac{\sqrt t}{\sqrt  n}\right)^{q-3}\mu_n(K_n)=C\left(\frac{\sqrt t}{\sqrt n}\right)^{q-3}d^{\frac{q-3}{q-1}}_{n+1}C_{2/q,q'}(K_{n}/d_{n+1})
\end {equation}
and inequality $(\ref{J3-8})$ follows.\qeda \medskip

Next we give the general proof. For this task we will use again the quasi-additivity with separated partitions.

\blemma {LJ3-2-0} Assume $N\geq 2$ and $\ell$ is an integer larger than $1$. 
There exists a positive constant $C_{1}=C_{1}(q,N,\ell)>0$ such that 
\begin {eqnarray}
\label {J3-8N} J'_{2,\ell}\leq C_{1}t^{-\frac{N}{2}}
\mysum{n=\ell}{a_t}e^{-\frac{n}{4}}d^{N-\frac{2}{q-1}}_{n+1}
C_{2/q,q'}\left(\myfrac{K_{n}}{d_{n+1}}\right).
\end {eqnarray}
\es 
\Proof As in the proof of \rth {upperW}, we know that there exists a finite number $J$, depending only on the dimension $N$,  of  separated sub-partitions $\{\#\Gth^h_{t,n}\}_{h=1}^J$ of the rescaled  sets $\tilde T_{n}=\sqrt{\frac{n+1}{t}\,}T_n$ by the $N$-dim balls 
$B_2(\tilde a_{n,j})$ where $\tilde a_{n,j}=\sqrt{\frac{n+1}{t}\,}{a_{n,j}}$, 
$\abs {a_{n,j}}=\myfrac{d_{n+1}+d_{n}}{2}$ and $\abs {a_{n,j}-a_{n,k}}\geq \sqrt{\frac{4t}{n+1}\,}$. Furthermore
$\#\Gth^h_{t,n}\leq Cn^{N-1}$. We denote $K_{n,j}=K_n\cap B_{\sqrt{\frac{t}{n+1}\,}}(a_{n,j})$.  
We write 
$\gm_n=\!\!\!\mysum{h=1}{J}\gm_n^h$, and accordingly $J'_{2,\ell}=\!\!\!\mysum{h=1}{J}
J_{2,\ell}'\,^{\!\!\!\!\!h}\,$, where $\gm_n^h=\mysum{j\in\Gth^h_{t,n}}{}\gm_{n,j}$, 
and $\gm_{n,j}$ are the capacitary measures of $K_{n,j}$ relative to $B_{n,j}=B_{6t/5\sqrt n}(a_n,j)$, which means
\begin{equation}\label {CM}
\nu_{n,j}(K_{n,j})=C_{2/q,q'}^{B_{n,j}}(K_{n,j})\;\;\mbox { and }\;\;\norm{\gn_{n,j}}_{W^{-2/q,q'}(B_{n,j})}=\left(C_{2/q,q'}^{B_{n,j}}(K_{n,j})\right)^{1/q}.
\end {equation}
Thus
$$
J'_{2,\ell}=\mysum{p=1}{a_{t}-\ell}
\dint_{\CT^{*}_{p}}(t-s)^{-\frac{N}{2}}e^{-\frac{\abs {x-y}^{2}}{4(t-s)}}\left[\mysum{n=p+\ell} {\infty}\;\mysum{h=1}{J}\;\mysum{j\in\Gth^h_{{t,n}}}{}
\BBH_{\gm_{n,j}}(y,s)\right]^q \!\!\!\!dyds. $$
We denote
$$
J_{2,\ell}'\,^{\!\!\!\!\!h}=\mysum{p=1}{a_{t}-\ell}
\dint_{\CT^{*}_{p}}(t-s)^{-\frac{N}{2}}e^{-\frac{\abs {x-y}^{2}}{4(t-s)}}\left[\mysum{n=p+\ell} {\infty}\;\mysum{j\in\Gth^h_{{t,n}}}{}
\BBH_{\gm_{n,j}}(y,s)\right]^q \!\!\!\!dyds, $$
 and clearly
 \begin{equation}\label{Su}
 J'_{2,\ell}\leq C\mysum{h=1}{J}J_{2,\ell}'\,^{\!\!\!\!\!h},
 \end {equation}
where $C$ depends only on $N$ and $q$. For integers $n$ and $p$ such that $n\geq\ell+1$, we set
$$\gl_{n,j,y}=\inf\{\abs {y-z}: z\in B_{\sqrt t/\sqrt {n+1}}(a_{n,j})\}=\abs {y-a_{n,j}}-\frac{\sqrt t}{\sqrt {n+1}}.
$$
Therefore 
$$\BA{l}
\mysum{n=p+\ell}{a_{t}}\myint{K_{n}}{}e^{-\frac{\abs {y-z}^{2}}{4s}}d\gm^h_{n}(z)=
\mysum{n=p+\ell}{a_{t}}\,\mysum{j\in\Gth^h_{t,n}}{}
\myint{K_{n,j}}{}e^{-\frac{\abs {y-z}^{2}}{4s}}d\gm_{n,j}(z)\\[2mm]
\phantom {\mysum{n=p+\ell}{a_{t}}\myint{K_{n}}{}e^{-\abs {y-z}^{2}/4t}d}
\leq \left(\mysum{n=p+\ell}{a_{t}}\,\mysum{j\in\Gth^h_{t,n}}{}e^{-\ge q'\frac{\gl_{n,j,y}^{2}}{4s}}\right)^{1/q'}
\left(\mysum{n=p+\ell}{a_{t}}\,
\mysum{j\in\Gth^h_{t,n}}{}e^{-q\gl_{n,j,y}^{2}\frac{1-\ge}{4s}}\gm^q_{n,j}(K_{n,j})\right)^{1/q}
\EA$$
where $\ge>0$ will be made precise later on. \medskip 

\noindent {\it Step 1} We claim that

\begin {equation}\label {kepler}
\mysum{n=p+\ell}{a_{t}}\,\mysum{j\in\Gth_{t,n}}{}e^{-\ge q'\frac{\gl_{n,j,y}^{2}}{4s}}
\leq C\sqrt\myfrac{ps}{t}
\end {equation}
where $C$ depends on $\ge$, $q$ and $N$. If $y$ is fixed in $T_{p}$, we denote by  $z_y$ the point of $T_n$ which solves $\abs {y-z_y}=\dist (y,T_n)$. Thus 
$$ \sqrt t(\sqrt n-\sqrt {p+1})\leq \abs {y-z_y}\leq t(\sqrt n-\sqrt p).$$
Let $Y=y\sqrt {t(p+1)}/\abs y$. On the axis $ \overrightarrow {0Y}$ we set ${\bf e}=Y/\abs Y$,  consider the points $b_{k}=(k\sqrt t/\sqrt n){\bf e}$ where 
$-n\leq k\leq n$ and denote by $G_{n,k}$ the spherical shell obtained by intersecting the spherical shell $T_{n}$ with the domain $H_{n,k}$ which is the set of points in $\BBR^N$ limited by the hyperplanes orthogonal to $ \overrightarrow {0Y}$
going through $((k+1)\sqrt t/\sqrt n){\bf e}$ and $((k-1)\sqrt t/\sqrt n){\bf e}$. The number of points $a_{n,j}\in G_{n,k}$ is smaller than  $C(n+1-\abs k)^{N-2}$, where $C$ depends only on $N$, and we denote by $\Gl_{n,k}$ the set of $j\in \Gth_{t,n}$ such that $a_{n,j}\in G_{n,k}$. Furthermore, if $a_{n,j}\in G_{n,k}$ elementary geometric considerations (Pythagora's theorem) imply that $\gl^2_{n,j,y}$ is greater than
$ t(n+p+1-2k\sqrt {p+1}/\sqrt n) $. Therefore
\begin{equation}\label {FN}
\BA {l}
\mysum{n=p+\ell}{a_{t}}\,\mysum{j\in \Gth_{t,n}}{}e^{-\ge q'\frac{\gl_{n,j,y}^{2}}{4s}}\leq C
\mysum{n=p+\ell}{a_{t}}\,\mysum{k=-n}{n}(n+1-\abs k)^{N-2}
e^{-\frac{\ge q'\left(n+p+1-2k\sqrt {p+1}/\right) t}{4s\sqrt n}}.
\EA\end {equation}
{\it Case $N=2$.} Summing a geometric series and using the inequality 
$\frac{e^u}{e^u-1}\leq 1+u^{-1}$ for $u>0$, we obtain
\begin{equation}\label{ke2}
\BA {l}
\mysum{k=-n}{n}e^{\frac{\ge q'\left(k\sqrt {p+1}\right) t}{2s\sqrt n}}
\leq e^{\frac{\ge q' t\sqrt {n(p+1)}}{2s}}
\myfrac{e^{\frac{\ge q't\sqrt {p+1}}{2s\sqrt n}}}{e^{\frac{\ge q't\sqrt {p+1}}{2s\sqrt n}-1}}\\[2mm]
\phantom{\mysum{k=-n}{n}e^{\frac{\ge q'\left(k\sqrt {p+1}\right) t}{2s\sqrt n}}}
\leq e^{\frac{\ge q' t\sqrt {n(p+1)}}{2s}}\left(1+\myfrac{2s\sqrt n}{\ge q't\sqrt {p+1}}\right).
\EA
\end {equation}
Thus, by comparison between series and integrals, 
\begin{equation}\label{ke3}\BA {l}
\mysum{n=p+\ell}{a_{t}}\,\mysum{j\in \Gth_{t,n}}{}e^{-\frac{\ge q'\gl_{n,j,y}^{2}}{4s}}
\leq C\mysum{n=p+\ell}{a_{t}}\left(1+\myfrac{s\sqrt n}{t\sqrt p}\right)e^{-\frac{\ge q'(\sqrt n-\sqrt {p+1\,})^2}{4s}}\\
\phantom{\mysum{n=p+\ell}{a_{t}}\,\mysum{j\in \Gth_{t,n}}{}e^{-\frac{\ge q'\gl_{n,j,y}^{2}}{4s}}}
\leq C\myint{p+1}{\infty}e^{-\frac{\ge q'(\sqrt x-\sqrt {p+1\,})^2t}{4s}}dx\\
\phantom{------------------}+
\myfrac{Cs}{t\sqrt p}\myint{p+1}{\infty}\sqrt xe^{-\frac{\ge q'(\sqrt x-\sqrt {p+1\,})^2t}{4s}}dx.
\EA
\end {equation}
Next
\begin{equation}\label{E1}
\BA {l}\myint{p+1}{\infty}e^{-\frac{\ge q'(\sqrt x-\sqrt {p+1\,})^2t}{4s}}dx
=2\myint{\sqrt {p+1}}{\infty}e^{-\frac{\ge q'(y-\sqrt {p+1\,})^2t}{4s}}ydy\\[4mm]
\phantom{\myint{p+1}{\infty}e^{-\frac{\ge q'(\sqrt x-\sqrt {p+1\,})^2t}{4s}}dx}
=2\myint{0}{\infty}e^{-\frac{\ge q'y^2t}{4s}}ydy+2\sqrt {p+1}\myint{0}{\infty}e^{-\frac{\ge q'y^2t}{4s}}dy\\
\phantom{\myint{p+1}{\infty}e^{-\frac{\ge q'(\sqrt x-\sqrt {p+1\,})^2t}{4s}}dx}
=\myfrac{2s}{t}\myint{0}{\infty}e^{-\frac{\ge q'z^2}{4}}zdz+2\sqrt\myfrac{(p+1)s}{t}
\myint{0}{\infty}e^{-\frac{\ge q'z^2}{4}}dz,
\EA
\end {equation}
and
\begin{equation}\label{E2}\BA {l}
\myint{p+1}{\infty}\sqrt xe^{-\frac{\ge q'(\sqrt x-\sqrt {p+1\,})^2t}{4s}}dx
=2\myint{\sqrt {p+1}}{\infty}e^{-\frac{\ge q'(y-\sqrt {p+1\,})^2t}{4s}}y^2dy\\
\phantom{\myint{p+1}{\infty}\sqrt xe^{-\frac{\ge q'(\sqrt x-\sqrt {p+1\,})^2t}{4s}}dx}
=2\myint{0}{\infty}e^{-\frac{\ge q'y^2t}{4s}}(y+\sqrt {p+1})^2dy\\
\phantom{\myint{p+1}{\infty}\sqrt xe^{-\frac{\ge q'(\sqrt x-\sqrt {p+1\,})^2t}{4s}}dx}
\leq 4\myint{0}{\infty}e^{-\frac{\ge q'y^2t}{4s}}y^2dy+4(p+1)
\myint{0}{\infty}e^{-\frac{\ge q'y^2t}{4s}}dy\\
\phantom{\myint{p+1}{\infty}\sqrt xe^{-\frac{\ge q'(\sqrt x-\sqrt {p+1\,})^2t}{4s}}dx}
\leq 4\left(\myfrac{s}{t}\right)^{3/2}
\myint{0}{\infty}e^{-\frac{\ge q'z^2}{4}}z^2dz+4(p+1)\sqrt \myfrac{s}{t}
\myint{0}{\infty}e^{-\frac{\ge q'z^2}{4}}dz
\EA
\end {equation}
Jointly with $(\ref{ke3})$, these inequalities imply
\begin{equation}\label{ke4}
\mysum{n=p+\ell}{a_{t}}\,\mysum{j\in \Gth_{t,n}}{}e^{-\frac{\ge q'\gl_{n,j,y}^{2}}{4s}}
\leq C\sqrt\myfrac{ps}{t}.
\end {equation}
{\it Case $N>2$}. Because the value of the right-hand side of $(\ref {FN})$ is an increasing value of $N$, it is sufficient to prove $(\ref{kepler})$ when $N$ is even, say $(N-2)/2=d\in\BBN_*$. There holds
\begin{equation}\label{ke5}
\BA {l}
\mysum{k=-n}{n}(n+1-\abs k)^de^{\frac{\ge q'\left(k\sqrt {p+1}\right) t}{2s\sqrt n}}
\leq 2\mysum{k=0}{n}(n+1- k)^de^{\frac{\ge q'\left(k\sqrt {p+1}\right) t}{2s\sqrt n}}.\EA
\end {equation}
We set
$$\ga=\ge q'\frac{t\sqrt {p+1}}{2s\sqrt n\,}\quad \mbox { and }\;
I_{d}=\mysum{k=0}{n}(n+1- k)^de^{k\ga}.
$$
Since
$$e^{k\ga}=\myfrac{e^{(k+1)\ga}-e^{k\ga}}{e^{\ga}-1},$$
we use Abel's transform to obtain
$$\BA {l}
 I_{d}=\myfrac{1}{e^{\ga}-1}\left(e^{(n+1)\ga}-(n+1)^d+
 \mysum{k=1}{n}\left((n+2- k)^d-(n+1- k)^d\right)e^{k\ga}\right)\\
 \phantom{ I_{d}}
 \leq \myfrac{1}{e^{\ga}-1}\left((1-d)e^{(n+1)\ga}-(n+1)^d+
 de^{\ga}\mysum{k=1}{n}\left((n+1- k)^{d-1}\right)e^{k\ga}\right).
\EA$$
Therefore the following induction holds
\begin{equation}\label{ke6}
I_d\leq \myfrac{de^{\ga}}{e^{\ga}-1}I_{d-1}.
\end {equation}
In $(\ref{ke2})$, we have already used the fact that 
$$\myfrac{de^{\ga}}{e^{\ga}-1}\leq C\left(1+\myfrac{s\sqrt n}{t\sqrt p}\right),
$$
and
$$ I_d\leq C\left(1+\left(\myfrac{s\sqrt n}{t\sqrt p}\right)^{d+1}\right)I_0.
$$
Thus $(\ref{ke3})$ is replaced by
\begin{equation}\label{ke7}\BA {l}
\mysum{n=p+\ell}{a_{t}}\,\mysum{j\in \Gth_{t,n}}{}e^{-\frac{\ge q'\gl_{n,j,y}^{2}}{4s}}
\leq C\mysum{n=p+\ell}{a_{t}}\left(1+\left(\myfrac{s\sqrt n}{t\sqrt p}\right)^{d+1}\right)e^{-\frac{\ge q'(\sqrt n-\sqrt {p+1\,})^2t}{4s}}\\
\phantom{\mysum{n=p+\ell}{a_{t}}\,\mysum{j\in \Gth_{t,n}}{}e^{-\frac{\ge q'\gl_{n,j,y}^{2}}{4s}}}
\leq C\myint{p+1}{\infty}e^{-\frac{\ge q'(\sqrt x-\sqrt {p+1\,})^2t}{4s}}dx\\
\phantom{------------------}+
\left(\myfrac{Cs}{t\sqrt p}\right)^{d+1}\myint{p+1}{\infty}x^{(d+1)/2}e^{-\frac{\ge q'(\sqrt x-\sqrt {p+1\,})^2t}{4s}}dx.
\EA
\end {equation}
The first integral on the right-hand side has already been estimated in $(\ref{E1})$, for the second integral, there holds
\begin{equation}\label{E3}\BA {l}
\myint{p+1}{\infty}x^{(d+1)/2}e^{-\frac{\ge q'(\sqrt x-\sqrt {p+1\,})^2t}{4s}}dx
=
\myint{0}{\infty}(y+\sqrt{p+1}\,)^{d+2}e^{-\frac{\ge q'y^2t}{4s}}dx\\
\phantom{\myint{p+1}{\infty}x^{(d+1)/2}e^{-\frac{\ge q'(\sqrt x-\sqrt {p+1\,})^2t}{4s}}dx}
\leq C\myint{0}{\infty}y^{d+2}e^{-\frac{\ge q'y^2t}{4s}}dy+Cp^{1+\frac{d}{2}}
\myint{0}{\infty}e^{-\frac{\ge q'y^2t}{4s}}dy\\
\phantom{\myint{p+1}{\infty}x^{(d+1)/2}e^{-\frac{\ge q'(\sqrt x-\sqrt {p+1\,})^2t}{4s}}dx}
\leq
C\left(\myfrac{s}{t}\right)^{2+\frac{d}{2}}\myint{0}{\infty}z^{(d+1)/2}e^{-\frac{\ge q'z^2}{4}}dz\\
\phantom{---------------------}
+C\left(\myfrac{s}{t}\right)^{3/2}p^{1+\frac{d}{2}}\myint{0}{\infty}e^{-\frac{\ge q'z^2}{4}}dz.
\EA
\end {equation}
Combining $(\ref {E1})$, $(\ref {ke7})$ and $(\ref {E3})$, we derive $(\ref{kepler})$.
\\

\noindent {\it Step 2. } Since $\CT_{p}^{*}\subset \Gg_{p}\ti [0,t]$ where 
$\Gg_{p}=B_{d_{p+1}}(x)\setminus B_{d_{p-1}}(x)$, 
$(y,s)\in\CT_{p}^{*}$ implies that $\abs {x-y}^{2}\geq (p-1)t$, thus
$J_{2,\ell}'\,^{\!\!\!\!\!h}\,$ satisfies
\begin {equation}\label {L1}\BA {l}
J_{2,\ell}'\,^{\!\!\!\!\!h}\,\leq 
Ct^{\frac{1-q}{2}}\mysum{p=1}{\ity}p^{\frac{q-1}{2}}\myint{0}{t}\myint{\Gg_{p}}{}(t-s)^{-\frac{N}{2}}
s^{-(q(N-1)+1)/2}
e^{-\frac{\abs {x-y}^{2}}{4(t-s)}}
\\[2mm]
\phantom {C\mysum{p=1}{\ity}\myint{0}{t}\myint{\Gg_{p}}{}p^{(N-1)(q-1)}(t-s)^{-\frac{N}{2}}}
\ti\mysum{n=p+\ell}{a_{t}}\,
\mysum{j\in\Gth^h_{t,n}}{}e^{-\frac{q\gl_{n,j,y}^{2}(1-\ge)}{4s}}\gm^q_{n,j}(K_{n,j})dsdy\\[2mm]
\phantom {J_{2,\ell}'\,^{\!\!\!\!\!h}\,}\leq 
Ct^{\frac{1-q}{2}}\mysum{n=\ell+1}{a_{t}}\,
\mysum{j\in\Gth^h_{t,n}}{}\gm^q_{n,j}(K_{n,j})\\[2mm]
\phantom {\mysum{n=\ell+1}{a_{t}}\,}
\ti
\mysum{p=1}{n-\ell}p^{\frac{q-1}{2}}\myint{0}{t}\myint{\Gg_{p}}{}
(t-s)^{-\frac{N}{2}}s^{-(q(N-1)+1)/2}
e^{-\abs {x-y}^{2}/4(t-s)}e^{-\frac{q\gl_{n,j,y}^{2}(1-\ge)}{4s}}dsdy

\EA\end {equation}
and the constant $C$ depends on $N, q$ and $\ge$. Next we set $q_{\ge}=(1-\ge)q$. 
Writting 
$$\abs {y-a_{n,j}}^{2}=\abs {x-y}^{2}+\abs 
{x-a_{n,j}}^{2}-2\langle y-x,a_{n,j}-x\rangle \geq pt+\abs 
{x-a_{n,j}}^{2}-2\langle y-x,a_{n,j}-x\rangle ,$$
we get
$$\int_{\Gg_{p}}e^{-\frac{q_{\ge}\abs {y-a_{n,j}}^{2}}{4s}}dy
= e^{-\frac{q_{\ge}\abs {x-a_{n,j}}^{2}}{4s}}\int_{\sqrt {tp}}^{\sqrt {t(p+1)}}
e^{-\frac{q_{\ge}r^{2}}{4s}}\int _{\abs {x-y}=r}
e^{2q_{\ge}\langle y-x, a_{n,j}-x\rangle/4s}dS_{r}(y)dr.
$$
For estimating the value of the spherical integral, we can 
assume that $a_{n,j}-x=(0,\ldots,0,\abs {a_{n,j}-x})$, $y=(y_{1},\ldots,y_{N})$ 
and, using spherical coordinates with center at $x$,  that the unit sphere has the representation 
$S^{N-1}=\{(\sin\gf.\gs,\cos\gf)\in \BBR^{N-1}\ti\BBR:\gs\in 
S^{N-2},\,\gf\in [0,\gp]\}$. With this representation, 
$dS_{r}=r^{N-1}\sin^{N-2}\gf\, d\gf\, d\gs$ and $\langle 
y-x,a_{n,j}-x\rangle=\abs {a_{n,j}-x}\abs {y-x}\cos\gf$. Therefore
$$\int _{\abs {x-y}=r}e^{2q_{\ge}\frac{\langle y-x,a_{n,j}-x\rangle}{4s}}dS_{r}(y)
=r^{N-1}\abs {S^{N-2}}\myint {0}{\gp}e^{2q_{\ge}\frac{\abs {a_{n,j}-x}r\cos\gf}{4s}}
\sin^{N-2}\gf\, d\gf.
$$
By \rlemma {A3}
\begin {equation}\BA {l}\label {L2}
\myint {\abs {x-y}=r}{}e^{2q_{\ge}\frac{\langle y-x,a_{n,j}-x\rangle}{4s}}dS_{r}(y)
\leq C\myfrac {r^{N-1}e^{2q_{\ge}\frac{r\abs {a_{n,j}-x}}{4s}}}{\left(1+\frac{r\abs 
{a_{n,j}-x}}{s}\right)^{\frac{N-1}{2}}}\\
[2mm]
\phantom {\myint {\abs {x-y}=r}{}e^{2q_{\ge}\frac{\langle y-x,a_{n,j}-x\rangle}{4s}}dS_{r}(y)}
\leq Cs^{\frac{N-1}{2}}\left(\myfrac {r}{\abs 
{a_{n,j}-x}}\right)^{\frac{N-1}{2}}e^{2q_{\ge}\frac{r\abs {a_{n,j}-x}}{4s}}.
\EA \end {equation}
Therefore
\begin {equation}\BA {l}\label {L3}
\myint{\Gg_{p}}{}e^{-q_{\ge}\frac{\abs {y-a_{n,j}}^{2}}{4s}}dy
\leq Ct^{\frac{N-1}{4}}p^{\frac{N-3}{4}}\myfrac {s^{\frac{N-1}{2}}e^{-q_{\ge}\frac{(\abs {a_{n,j}-x}-\sqrt 
{t(p+1)}\,)^{2}}{4s}}}
{\abs {a_{n,j}-x}^{\frac{N-1}{2}}},
\EA \end {equation}
and, since $\abs {a_{n,j}-x}\geq \sqrt {tn}$,
\begin {equation}\BA {l}\label {L4}\myint{0}{t}\myint{\Gg_{p}}{}(t-s)^{-\frac{N}{2}}
s^{-(q(N-1)+1)/2}
e^{-\frac{\abs{x-y}^2}{4(t-s)}}e^{-q_{\ge}\frac{\gl_{n,j,y}^{2}}{4s}}dy\,ds\\[2mm]
\phantom {-----}
\leq C\myfrac {\sqrt t p^{\frac{N-3}{4}}}{n^{\frac{N-1}{4}}}
\myint {0}{t}(t-s)^{-\frac{N}{2}}s^{-\frac{(q-1)(N-1)+1}{2}}e^{-\frac{pt}{4(t-s)}}
e^{-q_{\ge}\frac{(\sqrt {tn}-\sqrt {t(p+1)}\,)^{2}}{4s}}ds\\[4mm]
\phantom {-----}
\leq C\myfrac {t^{\frac{1-q(N-1)}{2}}p^{\frac{N-3}{4}}}{n^{\frac{N-1}{4}}}
\myint {0}{1}(1-s)^{-\frac{N}{2}}s^{-\frac{(q-1)(N-1)+1}{2}}e^{-\frac{p}{4(1-s)}}
e^{-q_{\ge}\frac{(\sqrt {n}-\sqrt {p+1}\,)^{2}}{4s}}.
\EA \end {equation}
We apply \rlemma {integral}, with $A=\sqrt p$, $B=\sqrt {q_{\ge}}(\sqrt n-\sqrt {p+1})$, 
$b=\frac{(q-1)(N-1)+1}{2}$, $a=\frac{N}{2}$ and $\gk=\sqrt {q_{\ge}}(\ell -1)/8$ as in the case $N=1$, 
and noticing that, for these specific values, 
$$\BA{l}
A^{1-a}B^{1-b}(A+B)^{a+b-2}=p^{\frac{2-N}{4}}(\sqrt {q_{\ge}}(\sqrt n-\sqrt {p+1}))^{\frac{1-(q-1)(N-1)}{2}}\\
\phantom{--------------------}
\ti
(\sqrt p+\sqrt {q_{\ge}}(\sqrt n-\sqrt {p+1}))^{\frac{(q-1)(N-1)+N-3}{2}}
\\[4mm]
\phantom {A^{1-a}B^{1-b}(A+B)^{a+b-2}}
\leq
C\left(\myfrac{n}{p}\right)^{\frac{N}{4}-1/2}\left(\myfrac{\sqrt n-\sqrt {p}}{\sqrt {n}}\right)^{\frac{1-(q-1)(N-1)}{2}},
\EA$$
where $C$ depends on $N$, $q$ and $\gk$. Therefore
\begin {equation}\label {L5}\BA {l}
\myint{0}{t}\myint{\Gg_{p}}{}(t-s)^{-\frac{N}{2}}s^{-\frac{N}{2}}
e^{-\frac{\abs{x-y}^2}{4(t-s)}}e^{-q_{\ge}\abs {y-z}^{2}/4s}dy\,ds\\[4mm]
\phantom {--}
\leq C\myfrac {t^{(1-q(N-1))/2}p^{\frac{N-3}{4}}}{n^{\frac{N-1}{4}}}\left(\myfrac{n}{p}\right)^{\frac{N}{4}-1/2}\left(\myfrac{\sqrt n-\sqrt {p}}{\sqrt {n}}\right)^{\frac{1-(q-1)(N-1)}{2}}e^{-\frac{(\sqrt p+\sqrt {q_{\ge}}(\sqrt n-\sqrt {p+1}))^{2}}{4}}\\[4mm]
\phantom {--}
\leq Ct^{\frac{1-q(N-1)}{2}}p^{-\frac{1}{4}}n^{\frac{(q-1)(N-1)-2}{4}}(\sqrt n-\sqrt {p})^{\frac{1-(q-1)(N-1)}{2}}
e^{-\frac{(\sqrt p+\sqrt {q_{\ge}}(\sqrt n-\sqrt {p+1}))^{2}}{4}}
.
\EA \end {equation}
We derive from $(\ref{L1})$, $(\ref{L5})$,
\begin {equation}\label {L6}\BA {l}
J_{2,\ell}'\,^{\!\!\!\!\!h}\,\leq Ct^{1-\frac{Nq}{2}}\\
\ti
\mysum{n=\ell+1}{a_{t}}\,
\mysum{j\in\Gth^h_{t,n}}{}n^{\frac{(q-1)(N-1)-2}{4}}\gm^q_{n,j}(K_{n,j})
\mysum{p=1}{n-\ell}p^{\frac{2q-3}{4}}(\sqrt n-\sqrt {p})^{\frac{1-(q-1)(N-1)}{2}}e^{- \frac{(\sqrt p+\sqrt {q_{\ge}}(\sqrt 
n-\sqrt {p+1}\,))^{2}}{4}}.
\EA \end {equation}
By \rlemma {A2} with $\ga=\myfrac{2q-3}{4}$, $\gb=\frac{1-(q-1)(N-1)}{2}$, $\gd=\frac{1}{4}$ and $\gg=q_\ge$, we obtain
\begin {equation}\label {L7}\BA {l}
\mysum{p=1}{n-\ell}p^{\frac{2q-3}{4}}(\sqrt n-\sqrt {p})^{\frac{1-(q-1)(N-1)}{2}}e^{- \frac{(\sqrt p+\sqrt {q_{\ge}}(\sqrt 
n-\sqrt {p+1}\,))^{2}}{4}}\leq Cn^{\frac{N(q-1)+q-3}{4}}e^{-\frac{n}{4}},
\EA \end {equation}
thus
\begin {equation}\label {L8}\BA {l}
J_{2,\ell}'\,^{\!\!\!\!\!h}\,\leq Ct^{1-\frac{Nq}{2}}
\mysum{n=\ell+1}{a_{t}}\,n^{\frac{N(q-1)}{2}-1}e^{-\frac{n}{4}}
\mysum{j\in\Gth^h_{t,n}}{}\gm^q_{n,j}(K_{n,j}).
\EA \end {equation}
Because 
$$\gm_{n,j}(K_{n,j})=C^{B_{n,j}}_{2/q,q'}(K_{n,j}),$$ 
we use the rescaling procedure as in the proof of \rlemma{QA}, except that the scale factor is $\sqrt{(n+1)t}$ instead of $\sqrt{n+1}$ so that the sets $\tilde T_n$, $\tilde K_n$, $\tilde \CQ_n$ and $\tilde K_n$ remains unchanged
Using again the quasi-additivity and the fact that $J'_{2,\ell}=\mysum{h=1}{J}J_{2,\ell}'\,^{\!\!\!\!\!h}\,$, we deduce
\begin {equation}\label {L11}\BA {l}
J_{2,\ell}\leq C't^{-\frac{N}{2}}
\mysum{n=\ell+1}{a_{t}}\,d^{N-\frac{2}{q-1}}_{n+1}e^{-\frac{n}{4}}C_{2/q,q'}\left(\myfrac{K_{n}}{d_{n+1}}\right),
\EA \end {equation} 
which implies $(\ref {J3-8N})$.\qeda\medskip 

The proof of \rth {lowerW} follows from the previous estimates on $J_1$ and $J_2$. Furthermore 
 the following integral expression holds
 \bth {lowerWint} Assume  $q\geq q_c$. Then there exists a positive constants $C_2^*$ , depending on $N$,$q$ and $T$, such that for any closed set $F$, there holds
\begin {equation}
\BA{l}
\label {lwe2}
\underline u_{F}(x,t)\geq \myfrac {C_2^*}{t^{1+\frac{N}{2}}}
\myint{0}{\sqrt {ta_t}}e^{-\frac{s^{2}}{4t}}s^{N-\frac{2}{q-1}}
C_{2/q,q'}\left(\myfrac {F}{s} \cap B_{1}(x)\right)s\,ds,
\EA
\end {equation}
where $a_t$ is the smallest integer $j$ such that $F\subset B_{\sqrt {jt}}(x)$.
\es 
\Proof We distinguish according $q=q_c$, or $q>q_c$, and for simplicity we denote $B_r=B_r(x)$ for the various values of $r$. \smallskip

\noindent{\it{Case 1: $q=q_c\Longleftrightarrow N-\frac{2}{q-1}=0$}}. Because 
$F_n=F\cap (B_{d_{n+1}}\setminus B_{d_{n}})$ 
there holds
 $$\BA {l}C_{2/q,q'}\left(\myfrac {F_n}{d_{n+1}} \right) \geq 
 C_{2/q,q'}\left(\myfrac {F}{d_{n+1}}\cap B_{1} \right)-
 C_{2/q,q'}\left(\myfrac {F\cap  B_{d_{n}}}{d_{n+1}} \right),
 \EA$$
Furthermore, since $d_{n+1}\geq d_{n}$, 
 $$C_{2/q,q'}\left(\myfrac {F\cap  B_{d_{n}}}{d_{n+1}} \right)
 =C_{2/q,q'}\left(\myfrac{d_{n}}{d_{n+1}}\myfrac {F\cap B_{d_{n}}}{d_{n}} \right)
 \leq C_{2/q,q'}\left(\myfrac {F}{d_{n}}\cap  B_{1} \right),
 $$
thus
 $$C_{2/q,q'}\left(\myfrac {F_n}{d_{n+1}} \right) \geq
  C_{2/q,q'}\left(\myfrac {F}{d_{n+1}}\cap B_{1} \right)
  - C_{2/q,q'}\left(\myfrac {F}{d_{n}}\cap B_{1} \right),
 $$
it follows
 $$\BA {l}
 \mysum{n=1}{a_{t}}e^{-\frac{n}{4}}
 C_{2/q,q'}\left(\myfrac {F_n}{d_{n+1}} \right)\geq 
 \mysum{n=1}{a_{t}}e^{-\frac{n}{4}}
 C_{2/q,q'}\left(\myfrac {F}{d_{n+1}}\cap B_{1} \right) -\mysum{n=1}{a_t}e^{-\frac{n}{4}} 
 C_{2/q,q'}\left(\myfrac {F}{d_{n}} \cap B_{1}\right)\\
 \phantom{ \mysum{n=1}{a_{t}}e^{-\frac{n}{4}}
 C_{2/q,q'}\left(\myfrac {F_n}{d_{n+1}} \right)}
  \geq  \mysum{n=1}{a_{t}}e^{-\frac{n}{4}} C_{2/q,q'}\left(\myfrac {F}{d_{n+1}}\cap B_{1} \right) - e^{-\frac{1}{4}}\mysum{n=0}{a_{_{t}}-1}e^{-\frac{n}{4}}C_{2/q,q'}\left(\myfrac {F}{d_{n+1}}\cap B_{1} \right)\\
  \phantom{ \mysum{n=1}{a_{t}}e^{-\frac{n}{4}}
 C_{2/q,q'}\left(\myfrac {F_n}{d_{n+1}} \right)}
  \geq
(1-e^{-\frac{1}{4}})\mysum{n=1}{a_{_{t}}-1}e^{-\frac{n}{4}}
C_{2/q,q'}\left(\myfrac {F}{d_{n+1}}\cap B_{1} \right) - e^{-\frac{1}{4}}
   C_{2/q,q'}\left(\myfrac {F}{\sqrt t}\cap B_{1} \right).
\EA $$
Since, by $(\ref{cap'1})$,
 $$
C_{2/q,q'}\left(\myfrac {F}{s'} \cap  B_{1}\right)\geq C_{2/q,q'}\left(\myfrac {F}{d_{n+1}}\cap B_{1} \right)\geq 
 C_{2/q,q'}\left(\myfrac {F}{s} \cap  B_{1}\right),
 $$
 for any $ s'\in[d_{n+1},d_{n+2}]$ and $ s\in[d_n,d_{n+1}]$, there holds
 $$\BA {l}
 te^{-\frac{n}{4}}C_{2/q,q'}\left(\myfrac {F}{d_{n+1}}\cap B_{1} \right)
 \geq
 C_{2/q,q'}\left(\myfrac {F}{d_{n+1}}\cap B_{1} \right)
  \myint{d_n}{d_{n+1}}e^{-s^2/4t}s\,ds\\
  \phantom{t e^{-\frac{n}{4}}C_{2/q,q'}\left(\myfrac {F}{d_{n+1}}\cap B_{1} \right)}
 \geq 
 \myint{d_n}{d_{n+1}}e^{-s^2/4t}
  C_{2/q,q'}\left(\myfrac {F}{s}\cap  B_{1}\right)s\,ds.
 \EA$$
 This implies 
$$\BA {l}
W_F(x,t)\geq (1-e^{-\frac{1}{4}}) t^{-(1+\frac{N}{2})}\myint{0}{\sqrt{ta_t}}e^{-s^2/4t}
  C_{2/q,q'}\left(\myfrac {F}{s}\cap  B_{1}\right)s\,ds .
\EA$$
\smallskip

\noindent{\it{Case 2: $q>q_c\Longleftrightarrow N-\frac{2}{q-1}>0$}}. In that case it follows from \rlemma{scal} that
$$ C_{2/q,q'}\left(\myfrac {F_n}{d_{n+1}} \right)\approx d_{n+1}^{\frac{2}{q-1}-N}C_{2/q,q'}\left(F_n \right).
$$
Thus
$$W_F(x,t)\approx t^{-1-\frac{N}{2}}\mysum{n=0}{a_t}e^{-\frac{n}{4}}C_{2/q,q'}\left(F_n \right).
$$
Since
$$C_{2/q,q'}\left(F_n\right) \geq 
 C_{2/q,q'}\left(F\cap B_{d_{n+1}} \right)-
 C_{2/q,q'}\left(F\cap B_{d_{n}} \right),
 $$
and again
 $$\BA {l}
 t^{-\frac{N}{2}}\mysum{n=0}{a_{t}}e^{-\frac{n}{4}}
 C_{2/q,q'}\left(F_n\right)\geq 
(1-e^{-\frac{1}{4}}) t^{-\frac{N}{2}}\mysum{n=0}{a_{_{t}}-1}e^{-\frac{n}{4}}
C_{2/q,q'}\left(F\cap B_{d_{n+1}} \right)\\
\phantom{ t^{-\frac{N}{2}}\mysum{n=0}{a_{t}}e^{-\frac{n}{4}}
 C_{2/q,q'}\left(F_n\right)}
\geq (1-e^{-\frac{1}{4}}) t^{-(1+\frac{N}{2})}\myint{0}{\sqrt{ta_t}}e^{-\frac{s^2}{4t}}
  C_{2/q,q'}\left(F\cap  B_{s}\right)s\,ds.
\EA $$
Because $  C_{2/q,q'}\left(F\cap  B_{s}\right)\approx   
s^{N-\frac{2}{q-1}}C_{2/q,q'}\left(s^{-1}F\cap  B_{1}\right)$, $(\ref{lwe2})$ follows.
\qeda\medskip
\section{Applications}

The first result of this section is the following
\bth {BIG} Assume $N\geq 1$ and $q> 1$. Then $\overline u_K=\underline u_K$.
\es
\Proof If $1<q<q_c$, the result is already proved in \cite{MV2}. The proof in the super-critical case is an adaptation that we  recall, for the sake of completeness. By \rth{upperWint} and \rth {lowerWint} there exists a positive constant $C$, depending on $N$, $q$ and $T$ such that
$$\overline u_F(x,t)\leq C\underline u_F(x,t)\qquad\forall (x,t)\in Q_T.
$$
By convexity $\tilde u=\underline u_F-\myfrac{1}{2C}(\overline u_F-\underline u_F)$ is a super-solution, which is smaller than $\underline u_F$ if we assume that $\overline u_F\neq \underline u_F$. If we set $\gth:=1/2+1/(2C)$, then $u_\gth=\gth\overline u_F$ is a subsolution. Therefore there exists a solution $u_1$ of $(\ref{mequ})$ in $Q_{\infty}$ such that
$u_\gth\leq u_1\leq \tilde u<\underline u_F$. If $\gm\in\mathfrak M_+^q(\BBR^N)$ satisfies $\gm (F^c)=0$, then $ u_{\gth\gm}$ is the smallest solution of $(\ref{mequ})$ which is above the subsolution $\gth u_{\gm}$. Thus $u_{\gth \gm}\leq u_1<\underline u_F$ and finally
$\underline u_F \leq u_1<\underline u_F$, a contradiction.\qeda\medskip

If we combine \rth{upperWint} and \rth {lowerWint} 
we derive the following integral approximation of the parabolic capacitary potential
\bprop{intrep} Assume $q\geq q_c$. Then there exist two positive constants $C^\dag_1$, 
$C^\dag_2$, depending only on $N$, $q$ and $T$ such that
\begin{equation}\label {equivpot}\BA {l}
 C^\dag_2t^{-(1+\frac{N}{2})}\myint{0}{\sqrt{ta_t}}s^{N-\frac{2}{q-1}}e^{-\frac{s^2}{4t}}
 C_{2/q,q'}\left(\myfrac{F}{s}\cap B_1(x)\right)s\,ds
 \leq W_F(x,t)\\
 \phantom{--------}
 \leq  
 C^\dag_1t^{-(1+\frac{N}{2})}\myint{\sqrt t}{\sqrt{t(a_t+2)}}s^{N-\frac{2}{q-1}}e^{-\frac{s^2}{4t}}
 C_{2/q,q'}\left(\myfrac{F}{s}\cap B_1(x)\right)s\,ds
\EA \end {equation}
 for any $(x,t)\in Q_T$.
\es

\bdef{intcap} If $F$ is a closed subset of $\BBR^N$, we define the $(2/q,q')$-integral parabolic capacitary potential $\CW_F$ by
 \begin{equation}\label{intpot}
 \CW_F(x,t)=t^{-1-\frac{N}{2}}\myint{0}{D_F(x)}s^{N-\frac{2}{q-1}}e^{-s^2/4t} C_{2/q,q'}\left(\myfrac{F}{s}\cap B_1(x)\right)s\,ds\qquad\forall (x,t)\in Q_\infty,
\end {equation}
where $D_F(x)=\max\{\abs{x-y}:y\in F\}$.
\es
An easy computation shows that
 \begin{equation}\label{intpot1}\BA {l}
0\leq\CW_F(x,t)- t^{-(1+\frac{N}{2})}\myint{0}{\sqrt{ta_t}}s^{N-\frac{2}{q-1}}e^{-\frac{s^2}{4t}}
 C_{2/q,q'}\left(\myfrac{F}{s}\cap B_1(x)\right)s\,ds\\[4mm]
 \phantom{-----------------------}
 \leq 
 C\myfrac{t^{(q-3)/2(q-1)}}{D_F(x)}e^{-D^2_F(x)/4t},
\EA\end {equation}
and
 \begin{equation}\label{intpot2}\BA {l}
0\leq t^{-(1+\frac{N}{2})}\myint{0}{\sqrt{t(a_t}+2)}s^{N-\frac{2}{q-1}}e^{-\frac{s^2}{4t}}
 C_{2/q,q'}\left(\myfrac{F}{s}\cap B_1(x)\right)s\,ds-\CW_F(x,t)\\[4mm]
 \phantom{-----------------------}
 \leq 
 C\myfrac{t^{(q-3)/2(q-1)}}{D_F(x)}e^{-\frac{D^2_F(x)}{4t}},
\EA\end {equation}
 for some $C=C(N,q)>0$. Furthermore
  \begin{equation}\label{intpot3}
\CW_F(x,t)=t^{-\frac{1}{q-1}}\myint{0}{D_F(x)/\sqrt t}s^{N-\frac{2}{q-1}}e^{-\frac{s^2}{4}}C_{2/q,q'}\left(\myfrac{F}{s\sqrt t}\cap B_1(x)\right)s\,ds.
\end {equation}

The following result gives a sufficient condition in order that $\overline u_F$ does not have a strong blow-up at a point $x$. 
\bprop{b-u} Assume $q\geq q_c$ and $F$ is a closed subset of $\BBR^N$. If there exists $\gg\in [0,\infty)$ such that
  \begin{equation}\label{intpot4}
\lim_{\gt\to 0}C_{2/q,q'}\left(\myfrac{F}{\gt}\cap B_1(x)\right)=\gg,
\end {equation}
then
  \begin{equation}\label{intpot5}
\lim_{t\to 0}t^{\frac{1}{q-1}}\overline u_F(x,t)=C\gg,
\end {equation}
for some $C=C(N,q)>0$.
\es
\Proof Clearly, condition $(\ref{intpot4})$ implies 
$$\lim_{t\to 0}C_{2/q,q'}\left(\myfrac{F}{\sqrt{t}s}\cap B_1(x)\right)=\gg
$$
for any $s>0$. Then $(\ref {intpot5})$ follows by Lebesgue's theorem. Notice also that the set of 
$\gg$ is bounded from above by a constant depending on $N$ and $q$.\qeda\medskip

In the next result we give a condition in order that the solution remains bounded at a point $x$. The proof is similar to the previous one.
\bprop{bnd} Assume $q\geq q_c$ and $F$ is a closed subset of $\BBR^N$. If
  \begin{equation}\label{intpot6}
\limsup_{\gt\to 0}\gt^{-\frac{2}{q-1}}C_{2/q,q'}\left(\myfrac{F}{\gt}\cap B_1(x)\right)<\infty,
\end {equation}
then $\overline u_F(x,t)$ remains bounded when $t\to 0$.
\es
\Remark If we assume that $f$ is a convex function on $\BBR^+$ satisfying
\begin{equation}\label{A}
c_2 r^q\leq f(r)\leq c_1r^q\qquad\forall r\geq 0
\end{equation}
for some $0<c_2\leq c_1$ we can construct in the same way as for $(\ref{mequ})$ the solutions
$\underline u_F$ and $\overline u_F$ for equation
\begin{equation}\label{A-A}
\prt_t u-\Gd u+f(u)=0\qquad\text{in }Q_T.
\end{equation}
The bilateral estimate estimate $(\ref{pot2})$ is still valid (up to change of the $C_i$). Since only convexity of $f$ is used in the proof of \rth{BIG}, there still holds $\underline u_F=\overline u_F$. Similar extensions of \rprop{b-u} and \rprop{bnd} are also clear.

\appendix\mysection {Appendix}
The next estimate is crucial in our study of semilinear parabolic 
equations.
\blemma {integral} Let $a$ and $b$ be two real numbers, $a>0$ and $\gk>0$. Then there exists a constant 
$C=C(a,b,\gk)>0$  
such that for any $A>0$, $B>\gk/A$ there holds
\begin {equation}\label {Aintegralest}
\myint{0}{1}(1-x)^{-a}x^{-b}e^{-A^2/4(1-x)}e^{-B^2/4x}dx\leq 
Ce^{-(A+B)^2/4}A^{1-a}B^{1-b}(A+B)^{a+b-2}.
\end {equation}
\es 
\Proof We first notice that 
\begin {eqnarray}\label {supexp}
\label {max}
\max \{e^{-A^2/4(1-x)}e^{-B^2/4x}:0\leq x\leq 1\}=e^{-(A+B)^2/4},
\end {eqnarray}  
and it is achieved for $x_0=B/(A+B)$. Set
$\Gf(x)=(1-x)^{-a}x^{-b}e^{-A^2/4(1-x)}e^{-B^2/4x}$, 
thus 
$$
\myint{0}{1}\Gf(x)dx=\myint{0}{x_{0}}\Gf(x)dx
+\myint{x_{0}}{1}\Gf(x)dx=I_{a,b}+J_{a,b}. 
$$
Put 
\begin {equation}\label {A1}
u=\myfrac {A^2}{4(1-x)}+\myfrac {B^2}{4x},
\end {equation}
then
\begin {equation}\label {A2}
4ux^{2}-(4u+B^{2}-A^{2})x+B^{2}=0.
\end {equation}
If $0<x<x_{0}$ this equation admits the solution
$$x=x(u)=\myfrac {1}{8u}\left(4u+B^{2}-A^{2}-
\sqrt {16u^{2}-8u(A^{2}+B^{2})+(A^{2}-B^{2})^{2}}\right)
$$
$$\BA{l}
\myint{0}{x_{0}}(1-x)^{-a}x^{-b}e^{-A^{2}/4(1-x)-B^{2}/4x}dx=
-\myint{(A+B)^{2}/4}\infty (1-x(u))^{-a}x(u)^{-b}e^{-u}x'(u)du
\EA
$$
Putting $x'=x'(u)$ and differentiating $(\ref {A2})$,
$$4x^{2}+8uxx'-(4u+B^{2}-A^{2})x'-4x=0\Longrightarrow
-x'=\myfrac{4x(1-x)}{4u+B^{2}-A^{2}-8ux}.
$$
Thus
\begin {equation}\label {A3}
\myint{0}{x_{0}}\Gf(x)dx=4\myint{(A+B)^{2}/4}{\infty}
\myfrac {(1-x(u))^{-a+1}x(u)^{-b+1}e^{-u}du}{4u+B^{2}-A^{2}-8ux(u)}.
\end {equation}
Using the explicit value of the root $x(u)$, we finally get
\begin {equation}\label {A4}
\myint{0}{x_{0}}\Gf(x)dx=4\myint{(A+B)^{2}/4}{\infty}
\myfrac {(1-x(u))^{-a+1}x(u)^{-b+1}e^{-u}du}
{\sqrt {16u^{2}-8u(A^{2}+B^{2})+(A^{2}-B^{2})^{2}}},
\end {equation}
and the factorization below holds
$$16u^{2}-8u(A^{2}+B^{2})+(A^{2}-B^{2})^{2}
=16(u-(A+B)^{2}/4)(u-(A-B)^{2}/4).
$$
We set $u=\gu+(A+B)^{2}/4$ and obtain
$$x(u)=\myfrac {v+(AB+B^{2})/2-\sqrt {v(v+AB)}}{2\left(v+(A+B)^{2}/4\right)},
$$
and
$$1-x(u)=\myfrac {v+(A^2+AB)/2+\sqrt {v(v+AB)}}{2\left(v+(A+B)^{2}/4\right)}.
$$
We introduce the relation $\approx$ linking two positive quantities depending on $A$ and $B$. It means that the two sided-inequalities up to multiplicative constants independent of $A$ and $B$. Therefore
\begin {equation}\label {A5}\BA {c}
\myint{0}{x_{0}}\Gf(x)dx= 2^{a-b-4}e^{-(A+B)^{2}/4}\myint{0}{\infty}\tilde\Gf (v)dv\quad \mbox {where }\\
\tilde\Gf (v)=
\myfrac{\left(v+(AB+B^{2})/2-\sqrt {v(v+AB)}\right)^{1-b}
\left(v+(A^2+AB)/2+\sqrt {v(v+AB)}\right)^{1-a}}
{\left(v+(A+B)^{2}/4\right)^{2-a-b}\sqrt {v(v+AB)}}e^{-v}dv.
\EA\end {equation}
{\it Case 1: $a\geq 1$, $b\geq 1$}. First 
\begin {equation}\label {A6}\BA {l}
\myfrac{\left(v+(A+B)^{2}/4\right)^{a+b-2}}{\sqrt {v(v+AB)}}\leq 
\myfrac{\left(v+(A+B)^{2}/4\right)^{a+b-2}}{\sqrt {v(v+\gk)}}
\approx \myfrac{\left(v+(A+B)^2\right)^{a+b-2}}{\sqrt {v(v+\gk)}}
\EA\end {equation}
since $a+b-2\geq 0$ and $AB\geq\gk$. Next
\begin {equation}\label {A7}\BA {l}\left(v+(A^2+AB)/2+\sqrt {v(v+AB)}\right)^{1-a}
\approx\left(v+A(A+B)\right)^{1-a}.
\EA\end {equation}
Furthermore
\begin {equation}\label {A8}\BA {l}
v+(AB+B^{2})/2-\sqrt {v(v+AB)}=B^2\myfrac{v+(A+B)^2/4}{v+B(A+B)/2+\sqrt {v(v+AB)}}\\
\phantom 
{v+(AB+B^{2})/2-\sqrt {v(v+AB)}}\approx
B^2\myfrac{v+(A+B)^2}{v+B(A+B)}.
\EA\end {equation}
Then
\begin {equation}\label {A9}\BA {l}
\left(v+(AB+B^{2})/2-\sqrt {v(v+AB)}\right)^{1-b}
\approx B^{2-2b}\left(\myfrac{v+B(A+B)}{v+(A+B)^2}\right)^{b-1}\\
\EA\end {equation}
It follows
\begin{equation}\label {N1}\BA {l}
\tilde\Gf(v)\leq CB^{2-2b}
\left(\myfrac{v+(A+B)^2}{v+A(A+B)}\right)^{a-1}
\myfrac{\left(v+B(A+B)\right)^{b-1}}{\sqrt {v(v+\gk)}}\\
\phantom{\Gf(x)}
\leq CB^{2-2b}
\left(\myfrac{v+(A+B)^2}{v+A(A+B)}\right)^{a-1}
\myfrac{v^{b-1}+ (B^2+AB)^{b-1}}{\sqrt {v(v+\gk)}}
\EA\end {equation}
where $C$ depends on $a$, $b$ and $\gk$. The function 
$v\mapsto (v+(A+B)^2)/(v+A(A+B))$ is decreasing on $(0,\infty)$. If we set
$$C_1=\myint{0}{\infty}\myfrac{v^{b-1}e^{-v}dv}{\sqrt {v(v+\gk)}}\quad \mbox {and }\;\;
C_2=\myint{0}{\infty}\myfrac{e^{-v}dv}{\sqrt {v(v+\gk)}}
$$
then
$$C_1\leq K(B^2+AB)^{b-1}C_2
$$
with $K=C_1\gk^{1-b}/C_2$. Therefore
\begin{equation}\label {N2}\BA {l}
\myint{0}{x_0}\Gf(x)dx\leq Ce^{-(A+B)^2/4}B^{1-b}A^{1-a}(A+B)^{a+b-2}.
\EA\end {equation}
The estimate of $J_{a,b}$ is 
obtained by exchanging $(A,a)$ with $(B,b)$ and replacing $x$ by $1-x$. {\it Mutadis mutandis}, this 
yields directely to the same expression as in $\ref {N2}$ and finally
\begin{equation}\label {N3}
\myint{0}{1}\Gf(x)dx\leq Ce^{-(A+B)^2/4}A^{1-a}B^{1-b}(A+B)^{a+b-2}.
\end {equation}
{\it Case 2: $a\geq 1$, $b< 1$}. Estimates $(\ref{A5})$, $(\ref{A6})$, $(\ref{A7})$, $(\ref{A8})$
and $(\ref{A9})$ are valid. Because $v\mapsto (v+B(A+B))^{b-1}$ is decreasing,  $(\ref{N1})$ has to be replaced by\begin{equation}\label {N4}\BA {l}
\tilde\Gf(v)\leq CB^{2-2b}
\left(\myfrac{v+(A+B)^2}{v+A(A+B)}\right)^{a-1}
\myfrac{\left(AB+B^2\right)^{b-1}}{\sqrt {v(v+\gk)}}.
\EA\end {equation}
This implies $(\ref{N2})$ directly. The estimate of $J_{a,b}$ is performed by the change of variable $x\mapsto 1-x$. If $x_1=1-x_0$ , there holds
$$J_{a,b}=\myint{0}{x_1}x^{-a}(1-x)^{-b}e^{-A^2/4x}e^{-B^2/4(1-x)}dx
=\myint{0}{x_1}\Psi(x)dx.
$$
Then
\begin {equation}\label {A'5}\BA {c}
\myint{0}{x_{1}}\Psi(x)dx= 2^{b-a-4}e^{-(A+B)^{2}/4}\myint{0}{x_{1}}\tilde\Psi(v)dv
\quad \mbox{where }\\
\tilde\Psi(v)=
\myfrac{\left(v+(AB+A^{2})/2-\sqrt {v(v+AB)}\right)^{1-a}
\left(v+(B^2+AB)/2+\sqrt {v(v+AB)}\right)^{1-b}}
{\left(v+(A+B)^{2}/4\right)^{2-a-b}\sqrt {v(v+AB)}}e^{-v}dv.
\EA\end {equation}
Equivalence $(\ref{A6})$ is unchanged; $(\ref{A7})$ is replaced by
\begin {equation}\label {A'7}\BA {l}\left(v+(B^2+AB)/2+\sqrt {v(v+AB)}\right)^{1-b}
\approx\left(v+B(A+B)\right)^{1-b},
\EA\end {equation}
$(\ref{A8})$ by
\begin {equation}\label {A'8}\BA {l}
v+(AB+A^{2})/2-\sqrt {v(v+AB)}\approx
A^2\myfrac{v+(A+B)^2}{v+A(A+B)},
\EA\end {equation}
and $(\ref{A9})$ by
\begin {equation}\label {A'9}\BA {l}
\left(v+(AB+A^{2})/2-\sqrt {v(v+AB)}\right)^{1-a}
\approx A^{2-2a}\left(\myfrac{v+A(A+B)}{v+(A+B)^2}\right)^{a-1}.\\
\EA\end {equation}
Because $a>1$,  $(\ref{N1})$ turns into
\begin{equation}\label {N'1}\BA {l}
\tilde\Psi(v)\leq CA^{2-2b}
(v+(A+B)^2)^{b-1}
\myfrac{(v+A^2+AB)^{a-1}(v+B^2+AB)^{1-b}}{\sqrt {v(v+\gk)}}\\[2mm]
\phantom{\Psi(x)}
\leq Ce^{-(A+B)^2/4}A^{2-2b}
(A+B)^{2b-2}\\
\phantom{--------}
\ti \myfrac{v^{a-b}+(A^2+AB)^{a-1}v^{1-b}+(B^2+AB)^{1-b}v^{a-1}
+A^{a-1}B^{1-b}(A+B)^{a-b}
}{\sqrt {v(v+\gk)}}.
\EA\end {equation}
Because $AB\geq \gk$, there exists a positive constant $C$, depending on $\gk$, such that
\begin{equation}\label {N'2}\BA {l}
\myint{0}{\infty}
 \myfrac{v^{a-b}+(A^2+AB)^{a-1}v^{1-b}+(B^2+AB)^{1-b}v^{a-1}
}{\sqrt {v(v+\gk)}}e^{-v}dv\\
\phantom{------------------}
\leq C A^{a-1}B^{1-b}(A+B)^{a-b}\myint{0}{\infty}\myfrac {e^{-v}dv
}{\sqrt {v(v+\gk)}}.
\EA\end {equation}
Combining $(\ref{N'1})$ and $(\ref{N'2})$ yields to
\begin{equation}\label {N'3}\BA {l}
\myint{0}{x_1}\Psi(x)dx\leq Ce^{-(A+B)^2/4}A^{1-a}B^{1-b}(A+B)^{a+b-2}.
\EA\end {equation}
This, again, implies that $(\ref{Aintegralest})$ holds.\\

\noindent {\it Case 3: $\max\{a,b\}<1$}. Inequalities $(\ref{A5})$-$(\ref{A9})$ hold, but $(\ref {N1})$ has to be replaced by
\begin{equation}\label {N5}\BA {l}
\tilde\Gf(v)\leq CB^{2-2b}
\left(\myfrac{v+(A+B)^2}{v+A(A+B)}\right)^{a-1}
\myfrac{\left(v+B^2+AB\right)^{b-1}}{\sqrt {v(v+\gk)}}\\[2mm]
\phantom{\Gf(x)}
\leq 
CB^{1-b}(A+B)^{2a+b-3}
\myfrac{v^{1-a}+\left(A^2+AB\right)^{1-a}}{\sqrt {v(v+\gk)}}
\EA\end {equation}
Noticing that
$$\myint{0}{\infty}\myfrac{v^{1-a}e^{-v}dv}{\sqrt {v(v+\gk)}}
\leq C\left(A^2+AB\right)^{1-a}\myint{0}{\infty}\myfrac{e^{-v}dv}{\sqrt {v(v+\gk)}},
$$
it follows that $(\ref{N2})$ holds. Finally $(\ref{N3})$ holds by exchanging $(A,a)$ and $(B,b)$.
\qeda

\medskip

\blemma {A2}. Let $\ga$, $\gb$, $\gg$, $\gd$ be real numbers and $\ell$ an integer. We assume $\gg>1$, $\gd>0$ and $\ell\geq 2$. Then there exists a positive constant $C$ such that, for any integer $n>\ell$
\begin {equation}\label {ser1}
\mysum{p=1}{n-\ell}p^{\ga}(\sqrt n-\sqrt p\,)^{\gb}e^{-\gd(\sqrt p+\sqrt\gg(\sqrt n-\sqrt {p+1}))^2}\leq Cn^{\ga-\gb/2}e^{-\gd n}.
\end {equation}
\es
\Proof The function $x\mapsto (\sqrt x+\sqrt\gg(\sqrt n-\sqrt {x+1}))^2$ is decreasing on $[(\gg-1)^{-1},\infty)$. Furthermore there exists $C>0$ depending on $\ell$, $\ga$ and $\gb$ such that 
$p^{\ga}(\sqrt n-\sqrt p\,)^{\gb}\leq Cx^{\ga}(\sqrt n-\sqrt {x+1}\,)^{\gb}$ for $x\in [p,p+1]$
If we denote by $p_0$ the smallest integer larger than $(\gg-1)^{-1}$, we derive
$$\BA {l}
S=\mysum{p=1}{n-\ell}p^{\ga}(\sqrt n-\sqrt p\,)^{\gb}e^{-(\sqrt p+\sqrt\gg(\sqrt n-\sqrt {p+1}))^2/4}=
\mysum{p=1}{p_0-1}+\mysum{p_0}{n-\ell}p^{\ga}(\sqrt n-\sqrt p\,)^{\gb}e^{-\gd(\sqrt p+\sqrt\gg(\sqrt n-\sqrt {p+1}))^2}\\
\phantom{S}
\leq \mysum{p=1}{p_0-1}p^{\ga}(\sqrt n-\sqrt p\,)^{\gb}e^{-\gd(\sqrt p+\sqrt\gg(\sqrt n-\sqrt {p+1}))^2}\\
\phantom{---------------}
+C\myint{p_0}{n+1-\ell}x^{\ga}(\sqrt n-\sqrt {x}\,)^{\gb}
e^{-\gd(\sqrt x+\sqrt\gg(\sqrt n-\sqrt {x+1}))^2}dx,
\EA$$
(notice that $\sqrt n-\sqrt x\approx \sqrt n-\sqrt {x+1}$ for $x\leq n-\ell$). Clearly 
\begin{equation}\label{ser2}
\mysum{p=1}{p_0-1}p^{\ga}(\sqrt n-\sqrt p\,)^{\gb}e^{-\gd(\sqrt p+\sqrt\gg(\sqrt n-\sqrt {p+1}))^2}
\leq C_0n^{\ga}(\sqrt n-\sqrt {n-\ell}\,)^{\gb}e^{-\gd n}
\end {equation}
for some $C_0$ independent of $n$. We set 
$y=y(x)=\sqrt {x+1}-\sqrt x/\sqrt\gg$. Obviously
$$y'(x)=\myfrac{1}{2}\left(\myfrac{1}{\sqrt{x+1}}-\myfrac{1}{\sqrt\gg\sqrt{x}}\right)\qquad\forall x\geq p_0,
$$
and their exists $\ge=\ge (\gd,\gg)>0$ such that
$\sqrt {2}\sqrt x\geq y(x)\geq \ge\sqrt x$ and $y'(x)\geq \ge/\sqrt x$. Furthermore
$$\sqrt x=\myfrac{\sqrt\gg\left(y+\sqrt{\gg y^2+1-\gg}\right)}{\gg-1},$$
$$\BA {l}\sqrt n-\sqrt x=\myfrac{\sqrt n(\gg-1)-\sqrt\gg y-\sqrt\gg\sqrt{\gg y^2+1-\gg}}{\gg-1}\\[4mm]
\phantom{\sqrt n-\sqrt x}
=\myfrac{n(\gg-1)+\gg-2y\sqrt{\gg n}-\gg y^2}{\sqrt n(\gg-1)-\sqrt\gg y+\sqrt\gg\sqrt{\gg y^2+1-\gg}}\\[5mm]
\phantom{\sqrt n-\sqrt x}
\approx \myfrac{n(\gg-1)+\gg-2y\sqrt{\gg n}-\gg y^2}{\sqrt n}
\EA$$
since $y(x)\leq \sqrt n$. Furthermore
$$\BA {l}n(\gg-1)+\gg-2y\sqrt{\gg n}-\gg y^2=
\gg(\sqrt{n+1}+\sqrt n/\sqrt\gg+y)(\sqrt{n+1}-\sqrt n/\sqrt\gg-y)\\
\phantom{n(\gg-1)+\gg-2y\sqrt{\gg n}-\gg y^2}
\approx \sqrt n (\sqrt{n+1}-\sqrt n/\sqrt\gg-y),
\EA$$
because $y$ ranges between $\sqrt {n+2-\ell}-\sqrt {n+1-\ell}\sqrt \gg\approx \sqrt n$ and $\sqrt {p_0+1}-\sqrt {p_0}\sqrt \gg$. 
Thus
$$\BA {l}(\sqrt n-\sqrt {x}\,)^{\gb}
\approx
\left(\sqrt{n+1}-\sqrt n/\sqrt\gg-y\right)^{\gb}.
\EA$$
This implies
\begin{equation}\label{ser3}\BA {l}
\myint{p_0}{n+1-\ell}x^{\ga}(\sqrt n-\sqrt {x}\,)^{\gb}
e^{-\gd(\sqrt x+\gg(\sqrt n-\sqrt {x+1}))^2}dx\\
\phantom{----}
\leq C\myint{y(p_0)}{y(n+1-\ell)}y^{2\ga+1}
\left(\sqrt{n+1}-\sqrt n/\sqrt\gg-y\right)^{\gb}e^{-\gg\gd(\sqrt n-y)^2}dy\\
\phantom{----}
\leq Cn^{\ga+\gb/2+1}
\myint{1-y(n+1-\ell)/\sqrt n}{1-y(p_0)/\sqrt n}(1-z)^{2\ga+1}
(z+\sqrt {1+1/n}-1-1/\sqrt\gg)^{\gb}e^{-\gg\gd nz^2}dz.
\EA\end {equation}
Moreover
\begin{equation}\label{ser4}\BA {l}
\phantom{--;;;;}
1-\myfrac {y(p_0)}{\sqrt n}= 1-\myfrac {1}{\sqrt 
  n}\left(\sqrt {p_0+1}-\myfrac 
 {\sqrt {p_0}}{\sqrt {\gg}}\right),\\[4mm]
1-\myfrac {y(n-\ell+1)}{\sqrt n}=1-\myfrac {\sqrt {n-\ell+2}}{\sqrt {n}}
 +\myfrac {\sqrt {n-\ell+1}}{\sqrt {n\gg}}\\
 \phantom{1-\myfrac {y(n-\ell+1)}{\sqrt n}}
 =\myfrac {1}{\sqrt {\gg}}
 \left(1 +\myfrac {\sqrt \gg\,(\ell-2)-\ell+1}{2 n} + 
\myfrac {\sqrt \gg\,(\ell-2)^2-(\ell-1)^2}{8n^{2}}\right)+O(n^{-3}).
\EA\end {equation}
 Let $\gth$ fixed such that $1-\myfrac {y(n-\ell+1)}{\sqrt n}<\gth<1-\myfrac {y(p_0)}{\sqrt n} $ for any $n>p_0$. Then
 $$\BA {l}\myint{\gth}{1-y(p_0)/\sqrt n}\!\!\!\!\!\!(1-z)^{2\ga+1}
(z+\sqrt {1+1/n}-1-1/\sqrt\gg)^{\gb}e^{-\gg\gd nz^2}dz\leq C_\gth
\myint{\gth}{1-y(p_0)/\sqrt n}(1-z)^{2\ga+1}e^{-\gg\gd nz^2}dz\\
\phantom{\myint{\gth}{1-y(p_0)/\sqrt n}\!\!\!\!\!\!(1-z)^{2\ga+1}
(z+\sqrt {1+1/n}-1-1/\sqrt\gg)^{\gb}e^{-\gg\gd nz^2}dz}
\leq C_\gth\; e^{-\gg\gd n\gth^2}\myint{\gth}{1-y(p_0)/\sqrt n}\!\!\!\!\!\!(1-z)^{2\ga+1}dz\\
\phantom{\myint{\gth}{1-y(p_0)/\sqrt n}\!\!\!\!\!\!(1-z)^{2\ga+1}
(z+\sqrt {1+1/n}-1-1/\sqrt\gg)^{\gb}e^{-\gg\gd nz^2}dz}
\leq C\; e^{-\gg\gd n\gth^2}\max \{1,n^{-\ga-1/2}\}.
 \EA$$
Because $\gg\gth^2>1$ we derive
\begin{equation}\label{ser5}\BA{l}
\myint{\gth}{1-y(p_0)/\sqrt n}\!\!\!\!\!\!(1-z)^{2\ga+1}
(z+\sqrt {1+1/n}-1-1/\sqrt\gg)^{\gb}e^{-\gg\gd nz^2}dz\leq Cn^{-\gb}e^{-\gd n},
\EA\end {equation}
for some constant $C>0$.
On the other hand
  $$\BA {l}\myint{1-y(n+1-\ell)/\sqrt n}{\gth}\!\!(1-z)^{2\ga+1}
(z+\sqrt {1+1/n}-1-1/\sqrt\gg)^{\gb}e^{-\gg\gd nz^2}dz\\
\phantom{\myint{1-y(n+1-\ell)/\sqrt n}{\gth}\!\!(1-z)^{2\ga+1}
----}
\leq C'_\gth
\myint{1-y(n+1-\ell)/\sqrt n}{\gth}\!\!
(z+\sqrt {1+1/n}-1-1/\sqrt\gg)^{\gb}e^{-\gg\gd nz^2}dz.\\
 \EA$$
The minimum of $z\mapsto (z+\sqrt {1+1/n}-1-1/\sqrt\gg)^{\gb}$ is achieved at $1-y(n+1-\ell)$ with value
 $$\myfrac{\sqrt\gg(\ell+1)+1-\ell}{2n\sqrt\gg}+O(n^{-2}),
 $$
and the maximum of the exponential term is achieved at the same point with value
$$e^{-n\gd+((\ell-2)\sqrt\gg +1-\ell)/2}(1+\circ (1))=C_\gg e^{-n\gd}(1+\circ (1)).
$$
We denote
 $$z_{\gg,n}=1+1/\sqrt\gg-\sqrt {1+1/n}\quad\mbox {and }\;
 I_\gb=\myint{1-y(n+1-\ell)/\sqrt n}{\gth}\!\!
(z-z_{\gg,n})^{\gb}e^{-\gg\gd nz^2}dz.
 $$
Since $1-y(n+1-\ell)\geq 1/\sqrt{2\gg}$ for $n$ large enough, 
 $$\BA {l}I_\gb
\leq \sqrt{2\gg}\myint{1-y(n+1-\ell)/\sqrt n}{\gth}\!\!
(z-z_{\gg,n})^{\gb}ze^{-\gg\gd nz^2}dz\\[4mm]
\phantom{I_\gb}
\leq\myfrac{-\sqrt{2\gg}}{2n\gg\gd}\left[(z-z_{\gg,n})^{\gb}e^{-\gg\gd nz^2}\right]_{1-y(n+1-\ell)/\sqrt n}^\gth+\myfrac{\gb\sqrt{2\gg}}{2n\gg\gd}
\myint{1-y(n+1-\ell)/\sqrt n}{\gth}\!\!
(z-z_{\gg,n})^{\gb-1}ze^{-\gg\gd nz^2}dz
 \EA$$
 But $1-y(n+1-\ell)/\sqrt n-z_{\gg,n}=(\ell-1)(1-1/\sqrt\gg)/2n$, therefore
\begin{equation}\label{ser7} 
I_\gb\leq C_1n^{-\gb-1}e^{-\gd n}+\gb C'_1n^{-1}I_{\gb-1}.
\end {equation}
 If $\gb\leq 0$ , we derive 
 $$ I_\gb\leq C_1n^{-\gb-1}e^{-\gd n},$$
 which inequality, combined with $(\ref{ser3})$ and $(\ref {ser5})$, yields to $(\ref{ser1})$.
If $\gb>0$, we iterate and get
 $$
 I_\gb\leq C_1n^{-\gb-1}e^{-\gd n}+C'_1n^{-1}(C_1n^{-\gb}e^{-\gd n}+(\gb-1)C'_1n^{-1}I_{\gb-2})
$$
If $\gb-1\leq 0$ we derive
 $$
 I_\gb\leq C_1n^{-\gb-1}e^{-\gd n}+C_1C'_1n^{-1-\gb}e^{-\gd n}=C_2n^{-\gb-1}e^{-\gd n},$$
which again yields to $(\ref{ser1})$. If $\gb-1>0$, we continue up we find a positive integer $k$ such that $\gb-k\leq 0$, which again yields to 
$$ I_\gb\leq C_kn^{-\gb-1}e^{-\gd n}
$$
and to $(\ref{ser1})$.\qeda\\

\medskip

The next estimate is fundamental in deriving the $N$-dimensional estimate.

\blemma {A3} For any integer $N\geq 2$ there exists a constant 
$c_{N}>0$ such that
\begin {equation}\label {spherical}
\int_{0}^\gp e^{m\cos\gth}\sin^{N-2}\gth \,d\gth \leq c_{N}\myfrac 
{e^m}{(1+m)^{(N-1)/2}}\qquad\forall m>0.
\end{equation}
\es
\Proof Put $\CI_{N}(m)=\myint{0}{\gp} e^{m\cos\gth}\sin^{N-2}\gth \,d\gth $.
Then $\CI_{2}'(m)=\myint{0}{\gp} e^{m\cos\gth}\cos\gth\,d\gth $
and 
$$\BA{l}\CI_{2}''(m)=\myint{0}{\gp} e^{m\cos\gth}\cos^{2}\gth\,d\gth 
=\CI_{2}(m)-\myint{0}{\gp} e^{m\cos\gth}\sin^{2}\gth \,d\gth \\[2mm]
\phantom {\CI_{2}''(m)}
=\CI_{2}(m)-\myfrac {1}{m}\myint{0}{\gp}e^{m\cos\gth}\cos\gth\,d\gth\\[2mm]
\phantom {\CI_{2}''(m)}=\CI_{2}(m)-\myfrac {1}{m}\CI'_{2}(m).
\EA$$
Thus $\CI_{2}$ satisfies a Bessel equation of order $0$. 
Since $\CI_{2}(0)=\gp$ and $\CI'_{2}(0)=0$, $\gp^{-1}\CI_{2}$ is the modified Bessel function of index 
$0$ (usually denoted by $I_{0})$ the asymptotic behaviour of which is 
well known, thus $(\ref {spherical})$ holds. If $N=3$
$$\CI_{3}(m)=\int_{0}^\gp e^{m\cos\gth}\sin\gth \,d\gth =
\left[\myfrac{-e^{m\cos\gth}}{m}\right]_{0}^\gp=\myfrac {2\sinh m}{m}.$$
For $N>3$ arbitrary
\begin {equation}\label{I3}\CI_{N}(m)=\myint{0}{\gp}\myfrac {-1}{m}\myfrac 
{d}{d\gth}(e^{m\cos\gth})\sin^{N-3}\gth\,d\gth
=\myfrac {N-3}{m}\myint{0}{\gp}
e^{m\cos\gth}\cos\gth\sin^{N-4}\gth\,d\gth.
\end {equation}
Therefore,
$$\CI_{4}(m)=\myfrac {1}{m}\myint{0}{\gp}
e^{m\cos\gth}\cos\gth\,d\gth=\CI_{2}'(m),
$$
and, again $(\ref {spherical})$ holds since $I_{0}'(m)$ has the same 
behaviour as $I_{0}(m)$ at infinity. For $N\geq 5$
$$\CI_{N}(m)=\myfrac {3-N}{m^{2}}
\left[e^{m\cos\gth}\cos\gth\sin^{N-5}\gth\right]_{0}^\gp+
\myfrac {N-3}{m^{2}}\myint{0}{\gp}e^{m\cos\gth}\myfrac{d}{d\gth}
\left(\cos\gth\sin^{N-5}\gth\right)d\gth.
$$
Differentiating $\cos\gth\sin^{N-5}\gth$ and using $(\ref {I3})$, we 
obtain
$$\CI_{5}(m)=\myfrac {4\sinh m}{m^{2}}-\myfrac {4\sinh m}{m^{3}},
$$
while
\begin {equation}\label{IN}
\CI_{N}(m) =\myfrac {(N-3)(N-5)}{m^{2}}\left(\CI_{N-4}(m)-\CI_{N-2}(m)\right),
\end {equation}
for $N\geq 6$. Since the estimate $(\ref{spherical})$ for $\CI_{2}$, $\CI_{3}$, 
$\CI_{4}$ and $\CI_{5}$ has already been obtained, a straigthforward induction 
yields to the general result.\qeda \\

\nind\Remark Although it does not has any importance for our use, it must 
be noticed that $\CI_{N}$ can be expressed either with hyperbolic 
functions if $N$ is odd, or with Bessel functions if $N$ is even.

\end {document}